\documentclass[11pt]{article}
\usepackage[ansinew]{inputenc}
\usepackage{graphicx}
\usepackage{hyperref}
\usepackage{url}
\usepackage[numbers,square]{natbib}
\usepackage{multirow}
\usepackage{tikz}
\usepackage{tikz-cd}
\usetikzlibrary{arrows,matrix}
\usepackage{authblk}
\usepackage{subcaption}
\usepackage{diagbox}

\usepackage{stmaryrd}
\usepackage{comment}
\usepackage{soul}

\usetikzlibrary{calc,patterns,shapes.geometric}
\usetikzlibrary{arrows,shapes,positioning}
\usetikzlibrary{decorations.markings,decorations.pathmorphing,decorations.pathreplacing}
\usetikzlibrary{decorations.text,decorations.shapes}

\usepackage[bbgreekl]{mathbbol}
\usepackage{amsmath,amssymb}
\usepackage{amsthm}  
\usepackage{bm}
\usepackage{bbm}
\usepackage{upgreek}

\def\bbb{\phantom{\big|} }

\def\mint{\int \mspace{-21.mu}\raise .8ex\hbox{\rotatebox{-80}{$|$}\,}}
\def\smint{\smallint \mspace{-13.mu}\raise -.1ex\hbox{\rotatebox{10}{--}}}

\def\bproof{\bigskip\noindent{\em Proof.}~}
\def\eproof{\hfill$\square$\nl}

\def\beq{\begin{equation}}
\def\eeq{\end{equation}}

\def\nl{\newline}

\def\<{\langle}
\def\>{\rangle}

\def\t{\tilde}
\def\wt{\widetilde}

\def\Chi{\raise .3ex \hbox{\large $\chi$}} 

\newcommand{\ii}{\ensuremath{\mathrm{i}}}

\def\dd{\, {\rm d}} 

\def\tb{\hbox{$\|\kern -.09em |$}}
\def\bigtb{\hbox{$\big\|\kern -.09em \big|$}}
\def\Bigtb{\hbox{$\Big\|\kern -.09em \Big|$}}

\def\px{\bX}  
\def\pv{\bV}  

\def\gx{\bx}  

\newcommand{\test}[1]{\check{#1}}

\def\RR{\mathbbm{R}}
\def\SS{\mathbbm{S}}

\def\UU{\mathbbm{U}}
\def\VV{\mathbbm{V}}

\newcommand{\mat}[1]{\mathbb{#1}}
\def\matd{\mathbb{d}}
\def\matr{\mathbb{r}}

\def\matB{\mathbb{B}}
\def\matC{\mathbb{C}}
\def\matD{\mathbb{D}}

\def\matF{\mathbb{F}}
\def\matG{\mathbb{G}}
\def\matI{\mathbb{I}}
\def\matJ{\mathbb{J}}
\def\matK{\mathbb{K}}
\def\matM{\mathbb{M}}

\def\matO{\mathbb{O}}
\def\matR{\mathbb{R}}
\def\matS{\mathbb{S}}
\def\matT{\mathbb{T}}

\def\matW{\mathbb{W}}

\def\matLambda{\mathbb{\Lambda}}

\def\LaB{\mathbb{\Lambda}} 

\newcommand{\arr}[1]{{\bm{\mathsf{#1}}}}
\newcommand\arre{{\bm{\mathsf{e}}}}
\newcommand\arrg{{\bm{\mathsf{g}}}}
\newcommand\arrB{{\bm{\mathsf{B}}}}

\newcommand\arrE{{\bm{\mathsf{E}}}}
\newcommand\arrU{{\bm{\mathsf{U}}}}
\newcommand\arrV{{\bm{\mathsf{V}}}}

\newcommand\arrX{{\bm{\mathsf{X}}}}

\newcommand\arrLambda{{\bm{\mathsf{\Lambda}}}}
\newcommand\arrsigma{{\bm{\mathsf{\sigma}}}}

\newcommand\ttc{{\texttt{c}}}
\newcommand\tte{{\texttt{e}}}
\newcommand\ttf{{\texttt{f}}}

\newcommand\sfk{{\sf k}}

\newcommand\sfF{{\sf F}}
\newcommand\sfG{{\sf G}}
\newcommand\sfH{{\sf H}}

\def\Vsp{V}

\def\cA{\mathcal{A}}

\def\cC{\mathcal{C}}

\def\cF{\mathcal{F}}
\def\cG{\mathcal{G}}
\def\cH{\mathcal{H}}
\def\cI{\mathcal{I}}

\def\cL{\mathcal{L}}

\def\cP{\mathcal{P}}

\def\cS{\mathcal{S}}

\usepackage{mathrsfs}

\def\scrC{\mathscr{C}}
\def\scrE{\mathscr{E}}
\def\scrF{\mathscr{F}}
\def\scrP{\mathscr{P}}

\newcommand\uvec{{\be}}

\def\bb{{\bs b}}

\def\be{{\bs e}}

\def\bk{{\bs k}}
\def\bsm{{\bs m}}
\def\bn{{\bs n}}

\def\bv{{\bs v}}

\def\bx{{\bs x}}

\def\bA{{\bs A}}
\def\bB{{\bs B}}
\def\bC{{\bs C}}
\def\bE{{\bs E}}
\def\bF{{\bs F}}
\def\bG{{\bs G}}

\def\bJ{{\bs J}}
\def\bK{{\bs K}}

\def\bM{{\bs M}}

\def\bR{{\bs R}}

\def\bV{{\bs V}}

\def\bX{{\bs X}}

\def\bsigma{{\bs \sigma}}
\def\bLambda{{\bs \Lambda}}

\def\bs{\boldsymbol}
\def\vp{{\varphi}}

\newcommand{\ee}{\ensuremath{\mathrm{e}}}

\def\Ampere{{Amp\`{e}re}}

\def\dt{{\partial_{t}}}

\DeclareMathOperator{\Span}{Span}
\DeclareMathOperator{\diag}{diag}
\DeclareMathOperator{\sinc}{sinc}

\DeclareMathOperator{\Div}{div}

\DeclareMathOperator{\curl}{curl}

\DeclareMathOperator{\grad}{grad}

\providecommand{\abs}[1]{\lvert#1\rvert}

\providecommand{\sprod}[2]{\langle#1,#2\rangle}
\providecommand{\bigsprod}[2]{\big\langle#1,#2\big\rangle}
\providecommand{\Bigsprod}[2]{\Big\langle#1,#2\Big\rangle}

\providecommand{\diff}[2]{\frac{\delta #1}{\delta #2}}

\newcommand{\fdt}{\frac{\partial }{\partial t}}

\providecommand{\range}[2]{\llbracket #1, #2 \rrbracket}

\providecommand{\martin}[1]{\textcolor{purple}{(Martin: #1)}}

\newcommand{\Lab}{\ensuremath{\boldsymbol{\Lambda}}}

\newcommand{\cbra}[2]{\ensuremath{\left \{ #1 , #2 \right \}}}

\newtheorem{theorem}{Theorem}
\newtheorem{lemma}{Lemma}
\newtheorem{remark}{Remark}
\newtheorem{definition}{Definition}
\newtheorem{assumption}{Assumption}

\newcommand\extended[1]{$~$\newline {\bf $~ \sharp ~ >> $ extended version: } #1 $~$ \newline {\bf $~ \sharp ~ << $} \newline}
\renewcommand\extended[1]{}

\title{Variational Framework for Structure-Preserving Electromagnetic Particle-In-Cell Methods}

\author[1]{Martin Campos Pinto}
\author[1,2]{Katharina Kormann}
\author[1,2]{Eric Sonnendr\"ucker }
\affil[1]{Max-Planck-Institut f\"ur Plasmaphysik, Garching, Germany}
\affil[2]{Technische Universit\"at M\"unchen, Zentrum Mathematik, Garching, Germany}

\begin{document}

\maketitle

\begin{abstract}
 In this article we apply a discrete action principle for the Vlasov--Maxwell equations in a
 structure-preserving particle-field discretization framework.
 In this framework the finite-dimensional electromagnetic potentials and fields are represented in
 a discrete de Rham sequence involving general finite element spaces, and the particle-field coupling
 is represented by a set of projection operators that commute with the differential operators.
 With a minimal number of assumptions which allow for a variety of finite elements and shape functions
 for the particles, we show that the resulting variational scheme has a general
 discrete Poisson structure and thus leads to a semi-discrete Hamiltonian system. By introducing discrete
 interior products we derive a second type of space discretization which is momentum preserving,
 based on the same finite elements and shape functions.
 We illustrate our method by applying it to spline finite elements, and to a new spectral
 discretization where the particle-field coupling relies on discrete Fourier transforms.
\end{abstract}

\section{Introduction}

Since the early days of Particle-in-Cell (PIC) schemes, plasma physicists have devised variational algorithms based on least action principles to
preserve key invariants such as the total energy and Gauss's laws \cite{Lewis:1970, Lewis:1972, eastwood1991virtual}.
In parallel, a Hamiltonian structure of the Vlasov--Maxwell equations has been proposed,
that involves a non-canonical Poisson bracket \cite{Morrison:1980,Weinstein:1981,Marsden:1982}.
Although the first methods were developed for finite difference field solvers,
many improvements have been made and in the last decade several schemes have been proposed that
rely on the de Rham structure of the Maxwell equations \cite{bossavit1998computational, hiptmair2002finite} to guarantee
an exact preservation of proper discrete Gauss laws for general Finite Element PIC methods on general
meshes~\cite{Campos.Pinto.Jund.Salmon.Sonnendrucker.2014.crm}, later extended
to variational PIC schemes in e.g.~\cite{Squire:2012} and \cite{Evstatiev:2013,Shadwick:2014},
where it was shown that variational spectral methods also preserve the total momentum of the plasma.

Following these ideas a Geometric Electromagnetic PIC (GEMPIC) method based on spline
finite elements has been proposed in \cite{kraus2016gempic}, that possess a Hamiltonian 
structure relying on a discrete Poisson bracket. 
Coupled with Hamiltonian splitting methods \cite{HairerLubichWanner:2006, Crouseilles:2015, He:2015},
this approach leads to fully discrete schemes that preserve a modified energy, discrete Gauss laws,
and the Poisson structure of the semi-discrete problem, including its associated Casimir invariants \cite{kraus2016gempic}.

In this article, we extend these constructions to a flexible and general setting that allows for arbitrary
structure-preserving discretizations of the electromagnetic fields and a variety of particle-field coupling operators
which in particular includes almost arbitrary smoothing shape functions.
By applying a discrete action principle we rigorously derive a variational system of discrete Vlasov--Maxwell equations, and we show that
it has a non-canonical Poisson structure. This approach allows for instance to derive numerical Maxwell solvers with a strong \Ampere{} and Gauss
equation, and also extends the strong Faraday solver of \cite{kraus2016gempic} to more general particle-field coupling schemes.
Another direct application is the design of variational spectral particle methods, where the Maxwell equations are solved in discrete Fourier spaces.

The outline of the paper is as follows. In Section~\ref{sec:variational} we first present the commuting de Rham complex
that serves as the basis of our discrete derivation. This setting is now common in the structure-preserving (mimetic)
discretization of Maxwell equations, and has been thoroughly studied in the Finite Element Exterior Calculus (FEEC) literature.
For the Vlasov--Maxwell equations it describes how the particle-field coupling operators are connected with the differential
operators involved in the discrete Maxwell equations.
Then, we derive a variational particle discretization of the Vlasov--Maxwell system in a strong \Ampere{} formulation from
a discrete action principle, and analyze its main conservation properties together with its discrete Poisson structure.
In Section~\ref{sec:momentum}, we present a variant of our method that preserves exactly the Gauss laws and the total momentum.
In Section~\ref{sec:hamiltonian_matrix} a matrix form of the equations is carefully detailed, which also allows to
derive a matrix form of the discrete Poisson bracket.
In Section~\ref{sec:SF}, we show how our analysis extends to a more general setting, and easily applies to the case of strong Faraday solvers.
A detailed application to the case of structure-preserving Spline and Fourier discretizations is then presented in Section~\ref{sec:interpolation_histopolation}, with particle-field coupling operators based on geometric degrees of freedom
which amount to discrete Fourier transforms in the spectral case.
In Section~\ref{sec:numerical_experiments}, we conclude with preliminary numerical experiments that validate our approach
and we compare the results obtained by various space discretizations that fit into our general framework,
including different Maxwell solvers and different orders of particle smoothing.


\section{Variational particle-field discretization}\label{sec:variational}

\subsection{Maxwell equations and particle trajectories} 

A kinetic description of the dynamics of a plasma in an 
electromagnetic field $(\bE, \bB)$ models the particles of species $s$
by a distribution function $f_s$ in phase-space that evolves according to the Vlasov equation
\begin{equation}
	\dt f_s(t, \bx, \bv) + \bv \cdot \nabla_\bx f_s(t, \bx, \bv)
    + \frac{q_s}{m_s}\big( \bE(t,\bx) + \bv \times \bB(t,\bx)   \big) \cdot \nabla_\bv f_s(t, \bx, \bv) = 0
\end{equation}
where $m_s$ and $q_s$ denote the mass and charge of the particle species $s$.
The self-consistent fields evolve according to Maxwell's equations
\begin{equation} \label{Maxwell}
\begin{aligned}
&\dt \bE(t, \bx) = \curl \bB(t, \bx) - \bJ(t, \bx) \\
&\dt \bB(t, \bx) = - \curl \bE(t, \bx)\\
&\Div \bE(t, \bx) = \rho(t, \bx) \\
&\Div \bB(t, \bx) = 0
\end{aligned}
\end{equation}
which are coupled to the Vlasov equation through the charge and current densities, 
\begin{equation}
	\rho(t, \bx) = \sum_s q_s \int_{\RR^3} f_s(t, \bx, \bv) \dd \bv,
  \qquad \bJ(t, \bx) = \sum_s q_s \int_{\RR^3} \bv f_s(t, \bx, \bv) \dd \bv.
\end{equation}
We refer to e.g. \cite{Glassey.1996.siam, ACL.2018} for a detailed presentation of these equations.
As the Vlasov equation is a conservative transport equation,
the distribution function $f_s$ is constant over time along the characteristic trajectories
for that species, which are solution to the characteristic ODEs 
\begin{equation} \label{traj}
\frac{\dd}{\dd t} \bX(t) = \bV(t),
\qquad
\frac{\dd}{\dd t} \bV(t) =  \frac{q_s}{m_s}\big( \bE(t, \bX(t)) + \bV(t) \times \bB(t, \bX(t))   \big) .
\end{equation}
In particle methods the distribution function is often represented by a collection of $N$ macro-particles
with phase-space positions $(\px_p, \pv_p)(t)$ and weights $w_p$, of the form
\begin{equation}  \label{fN}
 f^S_{s,N}(t, \bx, \bv) = \sum_{p=1}^N w_p S(\bx - \px_p(t)) \delta(\bv- \pv_p(t)),
\end{equation}
where $S$ is a shape function that can either by the Dirac $\delta$ distribution or some smoothing kernel,
depending on the particular configuration of the particle method.
Starting from a collection of initial positions $(\px_p^0, \pv_p^0)$, $p = 1, \dots, N$,
the weights are initialized so as to provide a good approximation to the initial density $f^0_s$,
and the particle positions are evolved according to some discrete characteristic equation,
in order to approximate the trajectories \eqref{traj}.
For the solution of Maxwell's equations, a grid-based solver is commonly used.

\subsection{Structure of Maxwell's equations and finite element exterior calculus}

As has been evidenced by several key contributions in the last decades
\cite{Bossavit.1988.IEE-A,bossavit1998computational,hiptmair2002finite},
the Maxwell equations \eqref{Maxwell} 
possess a geometric structure where a central role is played by de Rham sequence
\begin{equation}\label{deRham}
H^1(\RR^3)
\xrightarrow{ \mbox{$~ \grad ~$} }
H(\curl;\RR^3)
\xrightarrow{ \mbox{$~ \curl ~$}}
H(\Div;\RR^3)
\xrightarrow{ \mbox{$~ \Div ~$}}
L^2(\RR^3).
\end{equation}
In order to derive structure-preserving schemes we will follow 
the framework of Finite Element Exterior Calculus (FEEC) developed in e.g.~\cite{Monk.1993.jcam, hiptmair2002finite,%
Arnold.Falk.Winther.2006.anum, Arnold.Falk.Winther.2010.bams,%
buffa2010isogeometric, campos_pinto_compatible_2017_I}.
A central feature of these approaches is to involve a discretization that preserves the
sequence \eqref{deRham} at the discrete level, and that admits a sequence of
projection operators $\Pi^0, \dots, \Pi^3$ mapping infinite-dimensional function
spaces into discrete ones:
\begin{equation}\label{CD}
\begin{tikzpicture}[baseline=(current  bounding  box.center)]
  \matrix (m) [matrix of math nodes,row sep=3em,column sep=4em,minimum width=2em] {
    \bbb \Vsp^0 \bbb
      & \bbb \Vsp^1 \bbb
          & \bbb \Vsp^2 \bbb%
              & \bbb \Vsp^3 \bbb
    \\
           V_h^0 & V_h^1 & V_h^2 & V_h^3
    \\
    };
  \path[-stealth]
    (m-1-1) edge node [left]  {$\Pi^0$} (m-2-1)
    (m-1-1) edge node [above] {$\grad$} (m-1-2)
    (m-1-2) edge node [left]  {$\Pi^1$} (m-2-2)
    (m-2-1) edge node [above] {$\grad$} (m-2-2)
    (m-1-3) edge node [left]  {$\Pi^2$} (m-2-3)
    (m-1-4) edge node [left]  {$\Pi^3$} (m-2-4)
    (m-1-2) edge node [above] {$\curl$} (m-1-3)
    (m-2-2) edge node [above] {$\curl$} (m-2-3)
    (m-1-3) edge node [above] {$\Div$}  (m-1-4)
    (m-2-3) edge node [above] {$\Div$}  (m-2-4)
    ;
\end{tikzpicture}
\end{equation}
In our framework, it is these operators $\Pi^\ell$, together with some shape (smoothing) functions $S$,
that will encode the coupling mechanism between the particles and the discrete fields.
Here the top row contains the infinite-dimensional domain spaces $\Vsp^\ell$ of the operators $\Pi^\ell$,
which are in general proper subsets of the natural Hilbert spaces involved in the sequence \eqref{deRham},
and the bottom row consists of general discrete spaces such as finite-element or spectral spaces,
see e.g. Section~\ref{sec:interpolation_histopolation}.

A key ingredient in our variational derivation will be that
the operators $\Pi^\ell$ make the diagram {\em commuting}.
In practice many choices can be made for these operators and the associated finite-element spaces
where the fields are discretized. Each choice will result in a different coupling mechanism
between the particles and the fields, but all of them will lead to Hamiltonian systems,
provided the following property holds.
\begin{assumption} \label{as:Pi_cd}
  The operators $\Pi^\ell : V^\ell \to V^\ell_h$ are such that:
  \begin{itemize}
    \item
      the diagram \eqref{CD} commutes, i.e., we have
      \begin{align} \label{CD-nabla}
      \Pi^1 \grad G &= \grad \Pi^0 G \mspace{-80mu } &&\text{ for all } G \in \Vsp^0
      \\
      \label{CD-curl}
      \Pi^2 \curl \bG &= \curl \Pi^1 \bG \mspace{-80mu }&&\text{ for all } \bG \in \Vsp^1
      \\
      \label{CD-Div}
      \Pi^3 \Div \bG &= \Div \Pi^2 \bG \mspace{-80mu }&&\text{ for all } \bG \in \Vsp^2
      \end{align}
    \item
      the domain spaces $\Vsp^\ell$ are translation invariant function (or distribution) spaces,
      in the sense that if $G \in \Vsp^0$, then $G(\cdot - \bx) \in \Vsp^0$ for all $\bx \in \RR^3$.
  \end{itemize}
\end{assumption}

Since the commuting projection operators will be applied to particle shape functions,
we also need to specify when these shapes are admissible.
\begin{definition}[admissible shape functions] \label{def:adm_S}
  A shape function $S$ is said to be {\em admissible} for a given sequence of operators $\Pi^\ell$
  if it belongs to the domain spaces $\Vsp^0$ and $\Vsp^3$ of $\Pi^0$ and $\Pi^3$,
  and if for any $\be \in \RR^3$, $\be S$ belongs to the domains $\Vsp^1$ and
  $\Vsp^2$ of $\Pi^1$ and $\Pi^2$.
\end{definition}

\begin{remark}
  In practice, the translation invariance assumption corresponds to defining the projection operators
  on domain spaces $\Vsp^\ell$ characterized by some homogeneous regularity
  over $\RR^3$, which simplifies the notion of an admissible shape function $S$.
  In special cases where one works with localized or heterogeneous domain spaces,
  some additional care may need to be taken to guarantee that the projection operators
  can be applied on the shape functions.
\end{remark}

\subsection{Discretizing the \Ampere{} or Faraday equations in strong form}

In the article \cite{kraus2016gempic} the discretization ansatz was to consider fields in the spaces
\begin{equation}\label{fields_SF}
\t \phi_h \in V_h^0
\xrightarrow{ \mbox{$~ \grad ~$}}
\t \bE_h, \t \bA_h \in V_h^1
\xrightarrow{ \mbox{$~ \curl ~$}}
\t \bB_h \in V_h^2
\end{equation}
with $\t \phi_h$ and $\t \bA_h$ denoting discrete representations of the scalar and vector potentials,
and this has led to an approximation of \Ampere's and Faraday's laws in weak and strong form, respectively.
Although the analysis presented here readily applies to the ansatz \eqref{fields_SF}, it also covers the dual choice
\begin{equation}\label{fields_SA}
\bB_h \in V_h^1
\xrightarrow{ \mbox{$~ \curl ~$}}
\bE_h,\bA_h \in V_h^2
\xrightarrow{ \mbox{$~ \Div ~$}}
\phi_h \in V_h^3
\end{equation}
which leads to a new discrete model involving a strong \Ampere{} law and a weak Faraday law.
Throughout this article we will thus focus on this new ansatz \eqref{fields_SA},
and describe in Section~\ref{sec:SF} how our results apply to the `strong Faraday' ansatz \eqref{fields_SF}.

In both cases, the discrete equations in weak form will involve the discrete adjoints
to the strong differential operators, i.e.,
\begin{equation} \label{dw}
  \left\{
  \begin{aligned}
    & \grad_w : V_h^3 \to V_h^2, \qquad
        \int_\Omega (\grad_w \vp_h)\cdot \bF_h =  -\int_\Omega \vp_h \Div \bF_h
    \\
    & \curl_w : V_h^2 \to V_h^1, \qquad
        \int_\Omega (\curl_w \bF_h) \cdot \bC_h =  \int_\Omega \bF_h \cdot \curl \bC_h
    \\
    & \Div_w : V_h^1 \to V_h^0, \qquad
      \int_\Omega (\Div_w \bC_h)\psi_h = -\int_\Omega \bC_h \cdot \grad \psi_h
  \end{aligned}
  \right.
\end{equation}
for all $\vp_h \in V_h^3$, $\bF_h \in V_h^2$, $\bC_h \in V_h^1$, and $\psi_h \in V_h^0$.
These discrete operators may be seen as the discrete Riesz representants of the differential operators in distribution's sense.


\subsection{Discrete Action principle}
\label{sec:DAP}

We now derive a general geometric electromagnetic particle method where, following the
ansatz \eqref{fields_SA}, the \Ampere{} equation is discretized in a strong sense.
Here the coupling mechanism is essentially encoded in the abstract operators $\Pi^\ell$
that are only assumed to satisfy the commuting diagram properties, see Assumption~\ref{as:Pi_cd},
and in the shape function $S$ that must be admissible in the sense of Definition~\ref{def:adm_S}.

To do so we follow a discrete variational principle in the spirit of
\cite{Morrison:1980, Squire:2012, Kraus:2013:thesis, Evstatiev:2013},
based on Low's Lagrangian functional for the Vlasov--Maxwell equations~\cite{Low:1958.royal},
\begin{multline} \label{cL}
  \cL 
  = \sum_s \int f_s(t_0,\bx_0,\bv_0) \left( \big(m_s \bV + q_s \bA(t,\bX)\big)\cdot \bX' - \Big(\frac{m_s}{2}\bV^2 +q_s\phi(t,\bX)\Big)\right) \dd \bx_0 \dd \bv_0 \\
  +\frac {1}{2} \int_{\Omega}  |\grad\phi(t,\bx) +\bA'(t,\bx)|^2 \dd \bx  - \frac {1}{2}\int_{\Omega} |\curl \bA(t,\bx)|^2  \dd \bx.
\end{multline}
Here the curves $\bX=\bX(t;\bx_0,\bv_0)$, $\bX'=\bX'(t;\bx_0,\bv_0)$, $\bV=\bV(t;\bx_0,\bv_0)$ depend on time and
on the initial conditions, and we recall that in a variational derivation they represent {\em independent} variables of the
functional, in particular the prime symbol does not stand for a derivative.
We also note that a different set of characteristics is associated to each particle species, which has been left implicit here
for notational simplicity.

Formally, the Vlasov--Maxwell equations can be derived as the Euler-Lagrange equations associated
with this Lagrangian, as shown in \cite{Low:1958.royal}.
Here we will carefully apply this principle at the discrete level, starting from the discrete Lagrangian functional
\begin{multline} \label{cLh}
  \cL_h 
  = \sum_{p=1}^N w_p \left( \big(m_s \bV_p + q_s \bA^S(\bX_p)\big)\cdot \bX'_p - \Big(\frac{m_s}{2}\bV^2_p +q_s\phi^S(\bX_p)\Big)\right)  \\
  +\frac {1}{2} \int_{\Omega} |\grad_w\phi_h(\bx) +\bA'_h(\bx)|^2 \dd \bx  - \frac {1}{2}\int_{\Omega} |\curl_w \bA_h(\bx)|^2  \dd \bx.
\end{multline}
This Lagrangian is a function of discrete variables,
$ 
\cL_h = \cL_h(\arrX_N, \arrX'_N, \arrV_N, \bA_h, \bA'_h, \phi_h), 
$ 
where
$$
\arrX_N(t) = (\px_p(t))_{p = 1, \dots, N}, \quad
\arrX'_N(t) = (\px'_p(t))_{p = 1, \dots, N},
\quad
\arrV_N(t) = (\pv_p(t))_{p = 1, \dots, N} \quad \text{ in } (\RR^3)^N,
$$
are arbitrary collections of trajectories 
and $\bA_h(t), \bA'_h(t) \in V_h^2$, $\phi_h(t) \in V_h^3$ are arbitrary finite element potential fields.
In \eqref{cLh} the dependence on $t$ is implicit, and again we recall that the prime symbol does not mean a derivative,
as all these functions are independent in the variational derivation. 
Finally the coupling potentials are defined as
\begin{equation} \label{AphiS}
  \left\{
  \begin{aligned}
    & \bA^S(\px_p) := \sum_{\alpha = 1}^3 \uvec_\alpha \int_{\Omega} \Big( \bA_h \cdot \Pi^2(\uvec_\alpha S_{\px_p})\Big)  \dd \bx,
    \\
    &\phi^S(\px_p) := \int_{\Omega} \Big( \phi_h \Pi^3(S_{\px_p})\Big) \dd \bx
  \end{aligned}
  \right.
\end{equation}
where $S_{\px_p}(\bx) = S(\bx-\px_p)$ denotes the shape function centered on a particle. 
We note that $\cL_h$ is formally derived from the continuous functional \eqref{cL} by (i)
replacing the initial density $f_s$ by its Dirac approximation in \eqref{fN}, i.e.,
$
f^\delta_{s,N}(t^0, \bx^0, \bv^0) = \sum_{p=1}^N w_p \delta(\bx^0 - \px^0_p)\delta(\bv^0 - \pv^0_p),
$
(ii) using trajectories satisfying $(\bX,\bV)(t; \bX^0_p, \bV^0_p) = (\bX_p, \bV_p)(t)$,
(iii) potential fields in the discrete (finite element) spaces, (iv) weak discrete differentials \eqref{dw}
instead of the exact ones,
and finally (v) the coupling fields \eqref{AphiS} defined with admissible shape (smoothing) functions in \eqref{fN}.
In the case of several species, each density $f_s$ is approximated by a different set of discrete particles,
so that we actually have $s = s(p)$ in \eqref{cLh}. For this reason it will be convenient to denote
in the sequel particle masses and charges by
\begin{equation} \label{mp_qp}
  m_p := w_p m_{s(p)} \qquad \text{ and } \qquad q_p := w_p q_{s(p)}, \qquad \text{ for } ~~ p = 1, \dots, N.
\end{equation}
The discrete Action functional is then defined as
\begin{equation}\label{cS}
\cS_h(\arrX_N, \arrV_N, \phi_h, \bA_h) := \int_{0}^{T}
\cL_h\big((\arrX_N, \tfrac{\dd}{\dd t} {\arrX_N},\arrV_N,\phi_h, \bA_h, \dt \bA_h)(t)\big) \dd t
\end{equation}
and following a discrete action principle we look for generalized trajectories that form an extremum of $\cS_h$.
We already point out that the resulting equations will only involve the fields
\begin{equation} \label{EBpot}
\bE_h := - \dt \bA_h - \grad_w \phi_h \in V^2_h \qquad \text{ and } \qquad \bB_h := \curl_w \bA_h \in V^1_h,
\end{equation}
hence they will be gauge-independent.
Formally, extremality conditions for $\cS_h$ are associated to the Euler-Lagrange equations of
the discrete Lagrangian functional \eqref{cLh}.
Thus we look for $\arrX_N$, $\arrV_N$, $\phi_h$, and $\bA_h$ such that the following relations
hold for all $t \in [0,T]$, with functional Gateaux derivatives evaluated at
$(\arrX_N, \arrX'_N, \arrV_N, \phi_h, \bA_h, \bA'_h) = (\arrX_N, \tfrac{\dd}{\dd t}  \arrX_N, \arrV_N, \phi_h, \bA_h, \dt \bA_h)$:
\begin{align}
  \label{EL-V} 
&\Bigsprod{\diff{\cL_h}{\arrV_N}}{\test{\arrV}_N} = 0,
&&\forall ~\test{\arrV}_N \in \RR^{3N}
\\
\label{EL-X} 
&\Bigsprod{\diff{\cL_h}{\arrX_N}}{\test{\arrX}_N} = \Bigsprod{\fdt \diff{\cL_h}{\arrX_N'}}{\test{\arrX}_N},
&&\forall ~ \test{\arrX}_N \in \RR^{3N}
\\
\label{EL-phi} 
&\Bigsprod{\diff{\cL_h}{\phi_h}}{\test{\psi}_h} = 0,
&&\forall ~\test{\psi}_h \in V_h^3
\\
\label{EL-A} 
&\Bigsprod{\diff{\cL_h}{\bA_h}}{\test{\bA}_h} = \Bigsprod{ \fdt \diff{\cL_h}{\bA_h'}}{\test{\bA}_h},
&&\forall ~ \test{\bA}_h \in V_h^2.
\end{align}
For the variations with respect to $\arrV_N$, we compute
$$ 
\Bigsprod{\diff{\cL_h}{\arrV_N}}{\test{\arrV}_N} :=
  \left.\frac{\dd }{\dd \epsilon}\right|_{\epsilon=0}
  \Big(\cL_h(\arrX_N, \tfrac{\dd}{\dd t} \arrX_N,\arrV_N+\epsilon \test{\arrV}_N, \phi_h, \bA_h)\Big)
=
  \sum_{p}  m_p \Big(\frac{\dd\px_p}{\dd t} -\pv_p\Big) \cdot \test{\pv}_p
$$ 
for an arbitrary $\test{\arrV}_N = (\test{\pv}_p)_{p = 1, \dots, N}$,
so that \eqref{EL-V} gives
\begin{equation} \label{dtx}
\frac{\dd\px_p}{\dd t}   = \pv_p  \qquad \qquad \text{for } p = 1, \ldots, N.
\end{equation}
Using the coupling potentials \eqref{AphiS}, we compute for the variations with respect to $\arrX_N$
\begin{equation}
\label{diffLX}
\left\{\begin{aligned}
&\Bigsprod{\diff{\cL_h}{\arrX'_N}}{\test{\arrX}_N} =
  \sum_{p} \big(m_p \pv_p + q_p \bA^S(t,\bX_p)\big) \cdot \test{\px}_p
  = \sum_{p} m_p \pv_p \cdot \test{\px}_p + q_p \int_{\Omega} \bA_h \cdot \Pi^2(\test{\px}_p S_{\px_p})
\\
&\Bigsprod{\diff{\cL_h}{\arrX_N}}{\test{\arrX}_N} =
  -\sum_{p} q_p \int_{\Omega} \Big( \bA_h \cdot \Pi^2(\pv_p (\test{\px}_p \cdot \grad S_{\px_p}))
      -  \phi_h \Pi^3 (\test{\px}_p \cdot \grad S_{\px_p}) \Big)
\end{aligned}\right.
\end{equation}
for an arbitrary $\test{\arrX}_N = (\test{\px}_p)_{p = 1, \dots, N}$.
We then write Equation \eqref{EL-X} for a variation of a single particle $1 \le p \le N$
along the unit basis vector $\uvec _\alpha \in \RR^3$ for some dimension $1\le \alpha \le 3$.
Thus we take $\test{\px}_{p'} = \delta_{p',p} \uvec _\alpha$, which gives
$$
\begin{aligned} 
\frac{m_p}{q_p} \frac{\dd \pv_p }{\dd t} \cdot \uvec_\alpha &=
\int_{\Omega} \bA_h \cdot \Pi^2\big(\uvec_\alpha (\pv_p \cdot \! \grad S_{\px_p}) - \pv_p (\uvec_\alpha \cdot \! \grad S_{\px_p})\big)
\\
&\mspace{200mu}
  - \int_{\Omega} \dt \bA_h \cdot \Pi^2(\uvec_\alpha S_{\px_p})
  + \int_{\Omega} \phi_h \Pi^3 (\uvec_\alpha \cdot  \! \grad S_{\px_p})
\\
& =
\int_{\Omega} \bA_h \cdot \Pi^2 \curl (\uvec_\alpha \times \pv_p S_{\px_p})
- \int_{\Omega} \dt \bA_h \cdot \Pi^2(\uvec_\alpha S_{\px_p})
+  \int_{\Omega} \phi_h \Pi^3 \Div (\uvec_\alpha S_{\px_p})
\\
& = \int_{\Omega} \bA_h \cdot \curl \Pi^1 (\uvec_\alpha \times \pv_p S_{\px_p})
- \int_{\Omega} \dt \bA_h \cdot \Pi^2(\uvec_\alpha S_{\px_p})
+  \int_{\Omega} \phi_h \Div \Pi^3 (\uvec_\alpha S_{\px_p})
\\
& = \int_{\Omega} \curl_w \bA_h \cdot \Pi^1 (\uvec_\alpha \times \pv_p S_{\px_p})
- \int_{\Omega} (\dt \bA_h + \grad_w \phi_h) \cdot \Pi^2(\uvec_\alpha S_{\px_p})
\\
& = \int_{\Omega} \bB_h \cdot \Pi^1 (\uvec_\alpha \times \pv_p S_{\px_p})
+ \int_{\Omega} \bE_h \cdot \Pi^2(\uvec_\alpha S_{\px_p})
 \end{aligned}
$$
where we have used the commuting diagram property \eqref{CD-nabla}--\eqref{CD-Div}
of the operators $\Pi^\ell$, and the definition \eqref{EBpot} of the fields in the last equality.
Using the linearity of the projection operator we rewrite the magnetic rotation term
as
\begin{equation} \label{Brot}
\int_{\Omega} \bB_h \cdot \Pi^1 (\uvec_\alpha \times \pv_p S_{\px_p})
= \sum_{\beta = 1}^{3} (\uvec_\alpha \times \pv_p)_\beta \int_{\Omega} \bB_h \cdot \Pi^1 (\uvec_\beta S_{\px_p})
= \big(\pv_p \times \bB^S(\px_p)\big) \cdot \uvec_\alpha
\end{equation}
with a coupling magnetic field defined at the particle position as
\begin{equation}\label{BS}
    \bB^S(\px_p) := \sum_{\alpha = 1}^3 \uvec_\alpha \int_{\Omega} \bB_h(\bx) \cdot \Pi^1(\uvec_\alpha S_{\px_p})(\bx) \dd \bx.
\end{equation}
Defining similarly the coupling electric field by
\begin{equation}\label{ES}
    \bE^S (\px_p) := \sum_{\alpha = 1}^3 \uvec_\alpha \int_{\Omega} \bE_h(\bx) \cdot \Pi^2(\uvec_\alpha S_{\px_p})(\bx) \dd \bx
\end{equation}
we arrive at a velocity equation of the form
\begin{equation} \label{dtv}
  \frac{\dd \pv_p }{\dd t} = \frac{q_p}{m_p} \Big(\bE^S (\px_p) + \pv_p \times \bB^S(\px_p)\Big).
\end{equation}
Turning to the variations with respect to $\bA_h$, using again \eqref{AphiS} we compute
\begin{equation}
  \label{diffLA}
\left\{
  \begin{aligned}
  &
  \Bigsprod{\diff{\cL_h}{\bA_h}}{\test{\bA}_h} =
    \sum_{p} q_p \int_{\Omega} \Pi^2(\pv_p S_{\px_p}) \cdot \test{\bA}_h
  - \int_{\Omega} (\curl_w \bA_h)\cdot(\curl_w \test{\bA}_h)
  \\
  &
  \Bigsprod{\diff{\cL_h}{\bA_h'}}{\test{\bA}_h} =
  \int_{\Omega}(\dt \bA_h+\grad_w\phi_h)\cdot \test{\bA}_h
  \end{aligned}
\right.
\end{equation}
so that Equation \eqref{EL-A} gives
$$
\int_{\Omega} \dt (\dt \bA_h+\grad_w\phi_h)\cdot \test{\bA}_h
+ \int_{\Omega} (\curl_w \bA_h)\cdot(\curl_w \test{\bA}_h)
= \sum_{p} q_p \int_{\Omega} \Pi^2(\pv_p S_{\px_p}) \cdot \test{\bA}_h.
$$
The latter can be rewritten only in terms of the fields \eqref{EBpot} and the particle current defined as
$$
\bJ^S_N(t,\bx) := \sum_s q_s \int_{\RR^3} \bv f^S_{s,N}(t, \bx, \bv) \dd \bv = \sum_{p} q_p \pv_p S_{\px_p(t)}(\bx),
$$
see \eqref{fN} and \eqref{mp_qp}, as
\begin{equation}\label{dtEA}
  - \int_{\Omega} \dt \bE_h \cdot \test{\bA}_h
  + \int_{\Omega} \curl \bB_h \cdot\test{\bA}_h
  = \int_{\Omega} (\Pi^2 \bJ_N^S) \cdot \test{\bA}_h
\end{equation}
where we have used again the definition of the weak operators \eqref{dw}.
Since both $-\dt \bE_h + \curl \bB_h$ and $\Pi^2 \bJ_N^S$ belong to $V_h^2$,
and \eqref{dtEA} holds for all $\test{\bA}_h \in V_h^2$,
it leads to an \Ampere{} equation in strong form,
\begin{equation} \label{dtE}
  - \dt \bE_h + \curl \bB_h = \Pi^2 \bJ^S_N.
\end{equation}
In turn, a weak Faraday equation involving the discrete curl \eqref{dw},
\begin{equation} \label{dtB}
  \dt \bB_h + \curl_w \bE_h = 0
\end{equation}
follows from the definition of the fields \eqref{EBpot}:
Indeed, for all $\test{\bB}_h \in V_h^1$ we have
$$ 
\int_{\Omega} \dt \bB_h \cdot \test{\bB}_h = \int_{\Omega} \dt \bA_h \cdot \curl \test{\bB}_h =
-\int_{\Omega} (\bE_h + \grad_w \phi_h )\cdot \curl \test{\bB}_h
= - \int_{\Omega} \bE_h \cdot \curl \test{\bB}_h
$$ 
which amounts to \eqref{dtB}, by using the fact that
$\int_{\Omega} \grad_w \phi_h \cdot \curl \test{\bB}_h = \int_{\Omega} \phi_h \Div \curl \test{\bB}_h = 0$.
For the variations with respect to $\phi_h$ we use once more \eqref{AphiS} and compute
$$
\Bigsprod{\diff{\cL_h}{\phi_h}}{\test{\phi}_h} =
  - \sum_{p} q_p \int_{\Omega} (\Pi^3 S_{\px_p})\test{\phi}_h
+ \int_{\Omega} (\dt \bA_h+\grad_w\phi_h) \cdot \grad_w \test{\phi}_h
$$ 
for an arbitrary $\test{\phi}_h \in V_h^3$, so that \eqref{EL-phi} gives
\begin{equation} \label{var-Poisson}
  \int_{\Omega} ( \dt \bA_h + \grad_w \phi_h ) \cdot \grad_w \test{\phi}_h = \sum_{p} q_p \int_{\Omega} (\Pi^3 S_{\px_p})\test{\phi}_h.
\end{equation}
Using the field $\bE_h$ defined in \eqref{EBpot} and noting that \eqref{var-Poisson}
must hold for all $\test{\phi}_h \in V_h^3$,
we arrive at a Gauss law in strong form,
\begin{equation} \label{G-elec}
\Div \bE_h = \Pi^3 \rho^S_N
\qquad \text{ with } \qquad
  \rho^S_N(t,\bx) := \sum_s q_s \int_{\RR^3} f^S_{s,N}(t, \bx, \bv) \dd \bv = \sum_{p=1}^N q_p S_{\px_p(t)}(\bx),
\end{equation}
see again \eqref{fN}, \eqref{mp_qp}.
Finally a discrete magnetic Gauss law, this time in weak form,
follows again from the definition \eqref{EBpot} of $\bB_h = \curl_w \bA_h$, writing that
\begin{equation} \label{divwB}
\int_{\Omega} (\Div_w \bB_h ) \psi_h = - \int_{\Omega} \bB_h \cdot \grad \psi_h = -\int_{\Omega} \bA_h \cdot \curl \grad \psi_h = 0
\qquad \forall \psi_h \in V_h^0.
\end{equation}

\subsection{The variational equations}

Gathering the findings of the variational derivation just detailed, we obtain a system of semi-discrete equations where
the fields $\bE_h = \bE_h(t) \in V_h^2$ and $\bB_h = \bB_h(t) \in V_h^1$
are governed by the discrete \Ampere{} and Faraday equations
\begin{equation} \label{sawf-solve}
  \left\{
  \begin{aligned}
    - &\dt \bE_h + \curl \bB_h = \Pi^2 \bJ^S_N
    \\
    &\dt \bB_h + \curl_w \bE_h = 0
  \end{aligned}
  \right.
  \qquad  \text{ with } \qquad
    \Pi^2 \bJ^S_N = \sum_{p = 1 \cdots N} q_p \Pi^2( \pv_p S_{\px_p})
\end{equation}
with a weak $\curl_w : V_h^2 \to V_h^1$ defined by \eqref{dw},
and particles follow the trajectory equations
\begin{equation} \label{sawf-push}
\left\{
\begin{aligned}
  &\frac{\dd \px_p}{\dd t}  = \pv_p
  \\
  &
  \frac{\dd \pv_p}{\dd t}  = \frac{q_p}{m_p} \big(\bE^S(\px_p) + \pv_p \times \bB^S(\px_p)\big)
\end{aligned}
\right.
\qquad \text{ for } p = 1, \dots, N
\end{equation}
with coupling fields defined by 
\eqref{BS}--\eqref{ES}, namely
\begin{equation} \label{sawf-fields}
  \bE^S(\px_p) = \sum_{\alpha=1}^3 \uvec_\alpha\int_{\Omega} \bE_h \cdot \Pi^2(\uvec_\alpha S_{\px_p}),
  \qquad
  \bB^S(\px_p) = \sum_{\alpha=1}^3 \uvec_\alpha\int_{\Omega} \bB_h \cdot \Pi^1(\uvec_\alpha S_{\px_p})
\end{equation}
where $(\uvec_1, \uvec_2, \uvec_3)$ is an orthonormal basis of $\RR^3$.
These evolution equations are completed with two discrete Gauss laws,
\begin{equation} \label{sawf-Gauss}
  \left\{
  \begin{aligned}
    &\Div \bE_h = \Pi^3 \rho^S_N
    \\
    &\Div_w \bB_h = 0
  \end{aligned}
  \right.
\end{equation}
with $\rho^S_N = \sum_{p=1}^N q_p S_{\px_p}$ and the weak divergence operator $\Div_w : V_h^1 \to V_h^0$
defined by \eqref{dw}. We note that here the first Gauss law has been derived
from the variational principle (considering variations in the electric potential),
whereas the second one follows from the definition of the magnetic field.

\subsection{Derivation of a discrete Hamiltonian and an associated Poisson bracket}

In this section we describe how the above variational equations can be associated with a discrete Poisson bracket.

Following Hamilton's method \cite[Sec.~VI.1.2]{HairerLubichWanner:2006}, we observe that our discrete Lagrangian
has two nonzero conjugate momenta given by \eqref{diffLX} and \eqref{diffLA}, which we may identify with their
Riesz representant in the proper spaces.
Assuming that the discrete solution satisfies the variational equations \eqref{sawf-solve}--\eqref{sawf-Gauss}, we have
$$
P_{\arrX_N} := \diff{\cL_h}{\arrX'_N} \equiv \big(m_p \pv_p + q_p \bA^S(t,\bX_p)\big)_{p=1, \dots, N}
\qquad \text{ and } \qquad
P_{\bA_h} := \diff{\cL_h}{\bA'_h} \equiv - \bE_h
$$
which allows to define a discrete Hamiltonian $\cH_h = \cH_h(\arrX_N, \arrV_N, \bA_h, \phi_h)$ as
$$
\begin{aligned}
  \cH_h &:= \bigsprod{P_{\arrX_N}}{\arrV_N} + \bigsprod{P_{\bA_h}}{\dt \bA_h} - \cL_h
  \\
  &= \sum_{p=1}^N \big(m_p \pv_p + q_p \bA^S(t,\bX_p)\big) \cdot \pv_p - \int_{\Omega} \bE_h \cdot \dt \bA_h - \cL_h.
\end{aligned}
$$
Using the form of the coupling potential \eqref{AphiS} and the variational Gauss law \eqref{sawf-Gauss} we have
$$
\sum_p q_p \phi^S(\bX_p) = \int_\Omega \phi_h \sum_p q_p \Pi^3(S_{\bX_p}) = \int_\Omega \phi_h \rho^S_N
  = \int_\Omega \phi_h \Div \bE_h = -\int_\Omega \bE_h  \cdot \grad_w \phi_h,
$$
so that the resulting Hamiltonian can be reformulated as a function of the fields \eqref{EBpot}, namely
\begin{equation} \label{cH}
  \cH_h(\arrX_N, \arrV_N, \bE_h, \bB_h) = \sum_{p=1}^N \frac{m_p}{2}\bV^2_p + \frac {1}{2} \int_{\Omega} |\bE_h|^2  \dd \bx + \frac {1}{2} \int_{\Omega} |\bB_h|^2  \dd \bx.
\end{equation}
By construction this Hamiltonian is preserved by any solution satisfying the Euler-Lagrange equations \eqref{EL-V}--\eqref{EL-A}.
Following \cite[Sec.~40-A]{Arnold_1989}, a discrete Poisson bracket $\{\cF_h, \cG_h\}$ can then be associated to the
evolution equations \eqref{sawf-solve}--\eqref{sawf-fields}, such that
\begin{equation} \label{PoissonForm}
  \frac{\dd}{\dd t} \cF_h(\arrX_N, \arrV_N, \bE_h, \bB_h) = \{\cF_h,\cH_h\}
\end{equation}
holds for an arbitrary functional $\cF_h$ of the discrete solution. 
To identify this bracket we may simply consider linear functionals of the form
defined by
$$
\cF_h = \cF_{\test{\arrX}_N,\test{\arrV}_N,\test{\bE}_h,\test{\bB}_h} :
(\arrX_N,\arrV_N,\bE_h,\bB_h) \mapsto
\sum_{p=1}^N \bX_p \cdot \test{\bX}_p + \bV_p \cdot \test{\bV}_p + \int_\Omega \bE_h \cdot \test{\bE}_h + \int_\Omega \bB_h \cdot \test{\bB}_h,
$$
and $\cG_h = \cH_h$. Since the Poisson bracket should be a bilinear antisymmetric expression of the derivatives
of its respective functionals, which read (upon identification with their proper discrete Riesz representant)
$$
\diff{\cF_h}{\px_p} = \test{\px}_p,
\quad \diff{\cF_h}{\pv_p} = \test{\pv}_p,
\quad \diff{\cF_h}{\bE_h} = \test{\bE}_h,
\quad \diff{\cF_h}{\bB_h} = \test{\bB}_h
$$
and (for $\cG_h = \cH_h$),
$$
\diff{\cG_h}{\px_p} = 0,
\quad \diff{\cG_h}{\pv_p} = m_p \pv_p,
\quad \diff{\cG_h}{\bE_h} = \bE_h,
\quad \diff{\cG_h}{\bB_h} = \bB_h,
$$
and observing by linearity of $\cF$ that \eqref{PoissonForm} just amounts to the
evolution equations \eqref{sawf-solve}--\eqref{sawf-fields} written in weak forms,
with $\test{\arrX}_N,\test{\arrV}_N,\test{\bE}_h,\test{\bB}_h$ as test fields,
we verify that \eqref{PoissonForm} holds with the following discrete bracket
\begin{multline}\label{eq:disc_bracket}
  \cbra{\mathcal{F}_h}{\mathcal{G}_h} = \sum_{p=1}^N \Bigg[
   \frac{1}{m_p}\left( \diff{\cF_h}{\px_p}\cdot\diff{\cG_h}{\pv_p}- \diff{\cF_h}{\pv_p} \cdot \diff{\cG_h}{\px_p}\right)
  + \frac {q_p}{m_p^2}\bB^S(\px_p)\cdot \Big(\diff{\cF_h}{\pv_p}\times\diff{\cG_h}{\pv_p}\Big)
  \\ +
\frac {q_p}{m_p} 
 \int_\Omega \left(  \Pi^2 \left(S_{\px_p}\diff{\cF_h}{\pv_p}\right)\cdot\diff{\cG_h}{\bE_h}
  - \diff{\cF_h}{\bE_h}\cdot\Pi^2\left( S_{\px_p}\diff{\cG_h}{\pv_p}\right) \right) \dd\bx \Bigg] \\
  + \int_\Omega  
    \left( \diff{\cF_h}{\bE_h} \cdot  \curl\diff{\cG_h}{\bB_h}  - \curl\diff{\cF_h}{\bB_h} \cdot \diff{\cG_h}{\bE_h}\right)\dd\bx
  \end{multline}
where we remind that the coupling magnetic field $\bB^S(\px_p)$ is defined in \eqref{sawf-fields}
and involves the projection operator $\Pi^1$. We observe that this field plays the role of a parameter of the bracket,
as do the shape functions centered on the particle positions, $S_{\px_p}$.
A different role is played by the electric coupling terms, which enter the bracket through the product of $\bV$-$\bE$ derivatives.

Below we will verify that this bracket is a (non-canonical) Poisson bracket in the sense of \cite[Def.~VII.2.4]{HairerLubichWanner:2006},
in particular it satisfies the Jacobi identity. We note that other brackets involving different coupling fields $\bB^S$ would still be antisymmetric,
and hence also energy-preserving. As the different projection operators are connected by the commuting diagram properties which have been used
in several steps of the least action principle derivation, such brackets would probably not be variational,
but they could maybe still satisfy the Jacobi identity.

\subsection{Semi-discrete conservation properties of the variational system}

One major property of the above derivation is that the resulting semi-discrete system
has a Poisson structure, under the very general assumption that the diagram \eqref{CD} is commuting.
\begin{theorem}  \label{th:Poisson}
  If the operators $\Pi^\ell$ satisfy Assumption~\ref{as:Pi_cd} and if the shape function $S$ is admissible
  in the sense of Definition~\ref{def:adm_S},
  then the discrete bracket \eqref{eq:disc_bracket} is a (non-canonical) Poisson bracket and the
  semi-discrete equations 
  \eqref{sawf-solve}--\eqref{sawf-fields} are a Poisson system in the sense
  of \cite[Def.~VII.2.4]{HairerLubichWanner:2006}.
\end{theorem}
This result, whose proof will be given in Section~\ref{sec:proof}, implies in particular
that the evolution equations \eqref{sawf-solve}--\eqref{sawf-fields} preserve all
the functionals $\cF_h$ such that
$$
\{\cF_h, \cH_h\} = 0,
$$
which includes the Hamiltonian itself, $\cF_h = \cH_h$, but also all the Casimirs
of the bracket \eqref{eq:disc_bracket} which are the functionals $\cC_h$ such that
$\{\cC_h, \cG_h\} = 0$ for all $\cG_h$, and new Casimirs may be derived using
the Jacobi identity, see e.g. \cite{HairerLubichWanner:2006}.
An important example is provided by the functionals
\begin{equation}\label{GaussCas}
  \cC_h : (\arrX_N,\arrV_N,\bE_h,\bB_h)
  \mapsto
  \int_\Omega \test{\phi}_h \Big(\Div \bE_h - \Pi^3\big(\sum_{p=1}^N q_p S_{\px_p}\big)\Big)
\end{equation}
associated to an arbitrary $\test{\phi}_h \in V^3_h$. The fact that they
are Casimirs will be verified just below, and it implies that the discrete Gauss law
$\Div \bE_h = \Pi^3 \rho^S_N$ is preserved by our equations.

Theorem~\ref{th:Poisson} will be most conveniently proven on a matrix form of the equations,
which we will describe in Section~\ref{sec:hamiltonian_matrix}.
However a few basic conservation properties can be proven with a direct argument.

\begin{theorem}  \label{th:cons}
  Under the conditions of Theorem~\ref{th:Poisson},
  the evolution equations \eqref{sawf-solve}--\eqref{sawf-fields} preserve
  the discrete Hamiltonian \eqref{cH}
  %
  as well as the variational Gauss laws \eqref{sawf-Gauss}.
\end{theorem}

\begin{remark}[weak Gauss law]
  Similarly as for the GEMPIC method \cite{kraus2016gempic},
  the magnetic Gauss law plays the role of a pseudo-Casimir, in the sense that its
  conservation is actually needed to establish that the evolution system has
  a discrete Hamiltonian structure.
  With a strong-\Ampere{} ansatz \eqref{fields_SA}, we observe that 
  this divergence-free constraint is only preserved in a weak sense, see \eqref{divwB}.
  Although this may seem very weak, we will see below that it is the natural
  discrete invariant that provides a Poisson structure for the resulting Hamiltonian system.
\end{remark}

\bproof
The preservation of the magnetic Gauss law readily follows from the weak Faraday equation
in \eqref{sawf-solve}, indeed we have
$$
- \frac {\dd}{\dd t} \int_{\Omega} \bB_h \cdot \grad \vp_h
  = \int_{\Omega} \curl_w \bE_h \cdot \grad \vp_h = \int_{\Omega} \bE_h \cdot \curl \grad \vp_h = 0
$$
for all $\vp_h \in V_h^0$, using again the definition \eqref{dw} of the weak curl operator.
Turning to the electric Gauss law, we use $\frac{\dd}{\dd t} \px_p(t) = \pv_p$ to compute for
an arbitrary smooth function $\psi$
$$
\frac{\dd}{\dd t} \int_{\Omega} \rho^S_N(t,\bx) \psi(\bx) \dd \bx
  = \sum_{p = 1}^N q_p \int_{\Omega} S(\tilde \bx) \pv_p \cdot \grad \psi(\tilde \bx+\px_p) \dd \tilde \bx
  = \int_{\Omega} \bJ^S_N(t,\bx) \cdot \grad \psi(\bx) \dd \bx
$$
which shows that the continuity equation
\begin{equation} \label{cont}
  \dt \rho^S_N + \Div \bJ^S_N = 0
\end{equation}
always holds in distribution's sense, independently of the discrete particle trajectories.
Taking next the divergence of the discrete \Ampere{} equation in \eqref{sawf-solve},
the commuting diagram property \eqref{CD-Div} (which holds thanks to the admissibility of $S$)
allows us to write
$$
\dt \Div \bE_h = -\Div \Pi^2 \bJ^S_N = -\Pi^3 \Div \bJ^S_N = \dt \Pi^3 \rho^S_N
$$
where the last equality follows from \eqref{cont} and from the time-invariance of the
operator $\Pi^3$.
Integrating over time this shows that the electric Gauss law is indeed preserved.
Another argument consists of verifying that any functional of the form \eqref{GaussCas} is indeed a Casimir.
To do so we compute that the (Riesz representants of the) functional derivatives of $\cC_h$ read
$$
\diff{\cC_h}{\px_p} 
  = q_p \sum_{\alpha=1}^3 \Big(\int_\Omega \test{\phi}_h \Pi^3(\uvec_\alpha \cdot \grad S_{\px_p}) \Big) \uvec_\alpha
  = q_p \sum_{\alpha=1}^3 \Big(\int_\Omega \test{\phi}_h \Pi^3\Div (\uvec_\alpha S_{\px_p})\big) \Big) \uvec_\alpha
$$
and
$$
\diff{\cC_h}{\bE_h} = - \grad_w \test{\phi}_h.
$$
As for the derivatives $\diff{\cC_h}{\pv_p}$ and $\diff{\cC_h}{\bB_h}$, they vanish.
For the discrete bracket \eqref{eq:disc_bracket} we thus find
$$
\begin{aligned}
  \cbra{\mathcal{C}_h}{\mathcal{G}_h} &= \sum_{p=1}^N
   \frac{q_p}{m_p}\int_\Omega \left(  \test{\phi}_h \Pi^3\Div\Big(\diff{\cG_h}{\pv_p} S_{\px_p}\Big)
   +    \grad_w \test{\phi}_h \cdot\Pi^2\Big( S_{\px_p}\diff{\cG_h}{\pv_p}\Big) \right) \dd\bx
   \\
  & \mspace{60mu} - \int_\Omega  
    \grad_w \test{\phi}_h \cdot  \curl\diff{\cG_h}{\bB_h}  \dd\bx.
\end{aligned}
$$
Here the first term vanishes for arbitrary vectors $\diff{\cG_h}{\pv_p} \in \RR^3$,
by using the commuting diagram property and the definition of the weak gradient operator.
As for the second term, a discrete integration by parts yields
$\int_\Omega \grad_w \test{\phi}_h \cdot  \curl\diff{\cG_h}{\bB_h}  \dd\bx
= - \int_\Omega \test{\phi}_h \Div \curl\diff{\cG_h}{\bB_h}  \dd\bx = 0$,
which establishes that $\cbra{\mathcal{C}_h}{\mathcal{G}_h} = 0$ for any $\cG_h$. 
Equation~\eqref{PoissonForm} applied to $\cF_h = \cC_h$ then shows that
the quantity $\Div \bE_h - \Pi^3 \rho^S_N$ is an invariant of the evolution system.
Finally to verify the energy conservation, we may simply observe that the bracket \eqref{eq:disc_bracket}
is antisymmetric, so that $\cF_h = \cH_h$ is an obvious invariant of \eqref{PoissonForm}.
A more pedestrian argument is to first compute using \eqref{sawf-solve}
$$
\frac{\dd}{\dd t} \Big(\frac 12 \int_{\Omega} \abs{\bE_h}^2 + \abs{\bB_h}^2 \Big)
  =  \int_{\Omega} \bE_h \cdot (\curl \bB_h - \Pi^2 \bJ^S_N) - \bB_h \cdot \curl_w \bE_h
  =  - \int_{\Omega} \bE_h \cdot \Pi^2 \bJ^S_N
$$
where we have used the adjoint definition of $\curl_w$,
and then, using the trajectory equations \eqref{sawf-push}--\eqref{sawf-fields},
$$ 
\frac{\dd}{\dd t} \Big(\sum_{p = 1}^N \frac {m_p}{2} \abs{\pv_p}^2 \Big)
  = \sum_{p = 1}^N q_p \pv_p \cdot \big(\bE^S(\px_p) + \pv_p \times \bB^S(\px_p)\big)
  = \sum_{p = 1}^N q_p \int_{\Omega}  \bE_h\cdot \Pi^2(\pv_p S_{\px_p})
  =  \int_{\Omega} \bE_h \cdot \Pi^2 \bJ^S_N
$$ 
which shows that the discrete energy \eqref{cH} is indeed constant over time.
\eproof

\section{Generic Gauss and momentum preserving schemes} \label{sec:momentum}

Similarly as for the method in \cite{kraus2016gempic},
the semi-discrete scheme derived above is in general {\em not} momentum-preserving.
However it is possible to describe a general variant 
that preserves both the Gauss laws and a discrete momentum. This modified scheme 
comes at the price of losing the discrete Hamiltonian (Poisson) structure and the
conservation of energy, but it may be preferred for problems where momentum preservation is critical.

\subsection{Particle-field coupling with discrete interior products}

Our momentum-preserving schemes rely on discrete interior products of the form
\begin{equation} \label{I_a}
I^\ell_{\uvec_\alpha} = \cA_{h,\alpha} \imath^\ell_{\uvec_\alpha} : V_h^{\ell+1} \to V_h^\ell
\end{equation}
which involve the continuous interior products
$\imath^\ell_{\uvec_\alpha} : \Vsp^{\ell+1} \to \Vsp^{\ell}$
associated with a canonical unit vector $\uvec_\alpha$, $\alpha \in \range{1}{3}$, namely 
\begin{equation} \label{i_a}
  \imath^0_{\uvec_\alpha} \bC := \bC \cdot \uvec_\alpha,
\qquad
  \imath^1_{\uvec_\alpha} \bF := \bF \times \uvec_\alpha,
\qquad
  \imath^2_{\uvec_\alpha} g := g \uvec_\alpha,
\end{equation}
and where $\cA_{h,\alpha}$ is a linear approximation operator, 
such that the operators $I^\ell_{\uvec_\alpha}$ 
map every discrete space to its predecessor in the sequence,
as stated in \eqref{I_a}.

As a key property, denoting by $d^0 = \grad$, $d^1 = \curl$ and $d^2 = \Div$,
we require that the associated
discrete Lie derivatives, defined as
$$
L^\ell_{h,\uvec_\alpha} := d^{\ell-1} I^{\ell-1}_{\uvec_\alpha} +  I^{\ell}_{\uvec_\alpha} d^{\ell}  : \quad  V^\ell_h \to V^\ell_h
$$
are antisymmetric, in the sense that $\int_\Omega G \cdot L^\ell_{h,\uvec_\alpha} G = 0$ for all $G \in V^\ell_h$,
$\ell \in \{1,2\}$ and $1\le \alpha \le 3$.
Specifically, the momentum preserving properties will rely on the following relations
\begin{equation} \label{I01_rel}
\int_{\Omega} \bC_h \cdot \grad I^0_{\uvec_\alpha} \bC_h = -\int_{\Omega} \bC_h \cdot I^1_{\uvec_\alpha} \curl \bC_h
\qquad \forall \bC_h  \in V_h^1
\end{equation}
and
\begin{equation} \label{I12_rel}
\int_{\Omega} \bF_h \cdot \curl I^1_{\uvec_\alpha} \bF_h = - \int_{\Omega} \bF_h \cdot I^2_{\uvec_\alpha} \Div \bF_h
\qquad \forall \bF_h  \in V_h^2.
\end{equation}

\subsection{Gauss and momentum-preserving schemes}

Using the discrete interior products described above, we obtain the following result.

\begin{theorem} \label{th:mom}
  The scheme obtained by coupling the discrete Maxwell equations
  \eqref{sawf-solve}--\eqref{sawf-push} with the modified particle equations 
  \begin{equation} \label{sawf-push-mod}
  \left\{
  \begin{aligned}
    &\frac{\dd \px_p}{\dd t}  = \pv_p
    \\
    &
    \frac{\dd \pv_p}{\dd t}  = \frac{q_p}{m_p} \big(\bE^S(\px_p) + \bR^S(\bB_h,\px_p,\pv_p) \big)
  \end{aligned}
  \right.
  \qquad \text{ for } p = 1, \dots, N
  \end{equation}
  with coupling fields defined as
  \begin{equation} \label{sawf-fields-mod}
  \left\{
  \begin{aligned}
    &\bE^S (\px_p) = \sum_{\alpha = 1}^3 \uvec_\alpha \int_{\Omega} \bE_h(\bx) \cdot (I^2_{\uvec_\alpha} \Pi^3 S_{\px_p})(\bx) \dd \bx
  \\
    &\bR^S (\bB_h,\px_p,\pv_p) = -\sum_{\alpha = 1}^3 \uvec_\alpha \int_{\Omega} \bB_h(\bx) \cdot (I^1_{\uvec_\alpha} \Pi^2 (\pv_p S_{\px_p}))(\bx) \dd \bx,
  \end{aligned}
  \right.
  \qquad \text{ for } \alpha \in \range{1}{3},
  \end{equation}
  preserves the discrete Gauss laws \eqref{sawf-Gauss}, as well as the discrete momentum 
  \begin{equation} \label{Pa}
    {\bs\cP}_{h}(t) =  \sum_{p = 1}^N m_p \pv_p(t) - \sum_{\alpha=1}^3 \uvec_\alpha \int_{\Omega} (I^1_{\uvec_\alpha} \bE_h(t,\bx)) \cdot \bB_h(t,\bx) \dd \bx.
  \end{equation}
\end{theorem}

\begin{remark}
  Given the form \eqref{I_a}--\eqref{i_a} of $I^1_{\uvec_\alpha}$ and the linearity of $\cA_{h,\alpha}$, we have
  $$
  \int_{\Omega} (I^1_{\uvec_\alpha} \bE_h) \cdot \bB_h =
  \int_{\Omega} ((\cA_{h,\alpha} \bE_h) \times \uvec_\alpha ) \cdot \bB_h
  = \int_{\Omega} (\bB_h  \times  (\cA_{h,\alpha} \bE_h)) \cdot \uvec_\alpha
  $$
  which makes clear how \eqref{Pa} approximates the exact momentum along $\uvec_\alpha$.
  Similarly, we have
  \begin{equation} \label{R_tens}
    \bR^S (\bB_h,\px_p,\pv_p)  \cdot \uvec_\alpha
        =  \int_{\Omega} \Big(\big(\cA_{h,\alpha} \Pi^2(\pv_p S_{\px_p})\big)  \times \bB_h \Big) \cdot \uvec_\alpha
  \end{equation}
  which shows that the discrete magnetic force involved in \eqref{sawf-push-mod} is indeed an
  approximation of the ``natural'' term $\pv_p \times \bB_h(\px_p)$.
  However it is not possible in general to write
  $\bR^S (\bB_h,\px_p,\pv_p)$ as a product of the form $\pv_p \times \bB^S(\px_p)$
  for some field $\bB^S$, because the approximation operators $\cA_{h,\alpha}$ involved in the
  trajectory equation 
  depend a priori on the component $\alpha$ of the latter.
\end{remark}

\bproof
We first observe that the arguments used in the proof of Theorem~\ref{th:cons} for the conservation
of the discrete Gauss laws did not rely on the particle trajectory equation, hence they are still valid
for the modified scheme.
Turning to the discrete momentum, we compute using \eqref{sawf-push-mod}
$$
\begin{aligned}
\frac{\dd }{\dd t} \sum_{p = 1}^N m_p \pv_p \cdot \uvec_\alpha
  &= \sum_{p = 1}^N q_p \int_{\Omega} \big(\bE_h \cdot (I^2_{\uvec_\alpha} \Pi^3 S_{\px_p}) - \bB_h \cdot (I^1_{\uvec_\alpha}
  \Pi^2 (\pv_p S_{\px_p}))\big)
  \\
  &= \int_\Omega \bE_h\cdot (I^2_{\uvec_\alpha} \Pi^3\rho^S_N)
  - \int_\Omega  \bB_h\cdot (I^1_{\uvec_\alpha} \Pi^2 \bJ^S_N).
\end{aligned}
$$
Using next \eqref{sawf-solve} we write
$$
 \begin{aligned}
 \frac{\dd }{\dd t} \int_\Omega I^1_{\uvec_\alpha} \bE_h \cdot \bB_h
 &= - \int_\Omega I^1_{\uvec_\alpha} \bE_h \cdot \curl_w \bE_h
 + \int_\Omega (I^1_{\uvec_\alpha} (\curl \bB_h - \Pi^2 \bJ^S_N)) \cdot \bB_h
    \\
    &=       - \int_\Omega \curl I^1_{\uvec_\alpha} \bE_h \cdot \bE_h
    + \int_\Omega (I^1_{\uvec_\alpha} \curl \bB_h) \cdot \bB_h
     - \int_\Omega (I^1_{\uvec_\alpha} \Pi^2 \bJ^S_N) \cdot \bB_h
     \\
     &=
     \int_\Omega (I^2_{\uvec_\alpha} \Div \bE_h) \cdot \bE_h
     - \int_\Omega (\grad I^0_{\uvec_\alpha} \bB_h) \cdot \bB_h
      - \int_\Omega (I^1_{\uvec_\alpha} \Pi^2 \bJ^S_N) \cdot \bB_h
        \\
        &=
           \int_\Omega (I^2_{\uvec_\alpha} \Pi^3\rho^S_N) \cdot \bE_h
           - \int_\Omega (I^1_{\uvec_\alpha} \Pi^2 \bJ^S_N) \cdot \bB_h
        = \frac{\dd }{\dd t} \sum_{p = 1}^N m_p \pv_p \cdot \uvec_\alpha
\end{aligned}
$$
where we have used the definition of the weak curl operator in the second equality,
the relations \eqref{I01_rel}--\eqref{I12_rel} in the third one and
the preservation of the discrete (weak and strong) Gauss laws in the last one.
\eproof

\subsection{Interior products based on directional averaging on tensor-product spaces}


In this section we show that a simple construction based on directional averaging allows to
design momentum-preserving schemes when the compatible sequence
$$
V^0_h
\xrightarrow{ \mbox{$~ \grad ~$} }
V^1_h
\xrightarrow{ \mbox{$~ \curl ~$} }
V^2_h
\xrightarrow{ \mbox{$~ \Div ~$} } 
V^3_h
$$
involves tensor-product spaces of the form
\begin{equation} \label{V0_tens}
  V^0_h = \UU^{1}_{h} \otimes \UU^{2}_{h} \otimes \UU^{3}_{h}
  := \Span\Big(\Big\{ \bx \mapsto \Lambda^{0,1}_{k_1}(x_1) \Lambda^{0,2}_{k_2}(x_2) \Lambda^{0,3}_{k_3}(x_3) :
      (k_1, k_2, k_3) \in \prod_{\alpha =1}^3 \range{1}{N^\alpha_0} 
      \Big\} \Big)
\end{equation}
and
\begin{equation} \label{V123_tens}
V^1_h = \begin{pmatrix}
    \VV^{1}_{h} \otimes \UU^{2}_{h} \otimes \UU^{3}_{h}
    \\ \UU^{1}_{h} \otimes \VV^{2}_{h} \otimes \UU^{3}_{h}
      \\  \UU^{1}_{h} \otimes \UU^{2}_{h} \otimes \VV^{3}_{h}
  \end{pmatrix},
\qquad
V^2_h = \begin{pmatrix}
  \UU^{1}_{h} \otimes \VV^{2}_{h} \otimes \VV^{3}_{h}
  \\ \VV^{1}_{h} \otimes \UU^{2}_{h} \otimes \VV^{3}_{h}
    \\  \VV^{1}_{h} \otimes \VV^{2}_{h} \otimes \UU^{3}_{h}
  \end{pmatrix},
\qquad
V^3_h = \VV^{1}_{h} \otimes \VV^{2}_{h} \otimes \VV^{3}_{h},
\end{equation}
where the univariate spaces $\UU^\alpha_h$, $\VV^\alpha_h$,
form an exact sequence along each dimension 
$\alpha \in \range{1}{3}$, 
\begin{equation} \label{seq1d}
  \RR \rightarrow
  \UU^{\alpha}_{h} = \Span\Big(\{\Lambda^{0,\alpha}_k : k \in \range{1}{N^\alpha_0}\}\Big)
  \xrightarrow{ \mbox{$~ \partial_\alpha ~$} }
  \VV^{\alpha}_{h} = \Span\Big(\{\Lambda^{1,\alpha}_k : k \in \range{1}{N^\alpha_1}\}\Big)
  \rightarrow \{0\}.
\end{equation}

\begin{lemma}
  Assume that the univariate sequences \eqref{seq1d} are exact,
  with spaces $\UU^{\alpha}_{h}$ invariant over translations of $\pm h_\alpha$,
  $1 \le \alpha \le 3$.
  Then the discrete interior products 
  $ 
  I^\ell_{\uvec_\alpha} = \cA_{h,\alpha} \imath^\ell_{\uvec_\alpha} 
  $ 
  defined by composing the exact interior products \eqref{i_a} with the directional averaging operator,
  \begin{equation} \label{cA_av}
    (\cA_{h,1} G)(\bx) := \frac {1}{2h_1} \int_{x_1-h_1}^{x_1+h_1} G(y_1, x_2, x_3) \dd y_1
  \end{equation}
  and similarly for $\alpha = 2, 3$,
  map $V^{\ell+1}_h$ to $V^{\ell}_h$. Furthermore, they satisfy
  the relations \eqref{I01_rel}--\eqref{I12_rel}.
\end{lemma}

\bproof
Let us show that $I^0_{\uvec_\alpha}$ maps $V^1_h$ to $V^0_h$.
For a generic basis function in $V^1_h$, of the form
$$
\bLambda^1_{\alpha,\bk}(\bx)
  = \uvec_\alpha \Lambda^{1,\alpha}_{k_\alpha} (x_\alpha)\prod_{\beta \neq \alpha} \Lambda^{0,\beta}_{k_\beta}(x_\beta),
$$
we observe that
$
(\imath^0_{\uvec_\alpha} \bLambda^1_{\alpha',\bk})(\bx) = \delta_{\alpha, \alpha'} \Lambda^{1,\alpha}_{k_\alpha} (x_\alpha)\prod_{\beta \neq \alpha} \Lambda^{0,\beta}_{k_\beta}(x_\beta)
$ using \eqref{i_a} and the tensor-product structure \eqref{V0_tens}--\eqref{seq1d}.
The exact sequence property \eqref{seq1d} then allows us to write
$\Lambda^{1,\alpha}_{k_\alpha} = \partial_\alpha \Gamma^{0,\alpha}_{k_\alpha}$
for some $\Gamma^{0,\alpha}_{k_\alpha} \in \UU^\alpha_h$,
which yields
$$
(I^0_{\uvec_\alpha} \bLambda^1_{\alpha,\bk})(\bx)
= \big( \cA_{h,\alpha} \Lambda^{1,\alpha}_{k_\alpha} \big)(x_\alpha) \prod_{\beta \neq \alpha} \Lambda^{0,\beta}_{k_\beta}(x_\beta)
= \frac {1}{2h_\alpha} \left[ \Gamma^{0,\alpha}_{k_\alpha} \right]_{x_\alpha-h_\alpha}^{x_\alpha+h_\alpha}
 \prod_{\beta \neq \alpha} \Lambda^{0,\beta}_{k_\beta}(x_\beta)
$$
which belongs to $V^0_h$, according to \eqref{V0_tens} and the discrete translation invariance.
The argument for the other spaces is similar.
Turning to \eqref{I01_rel}--\eqref{I12_rel} we next observe that the directional averaging operators are of the form
$\cA_{h,\alpha} G = \mu_\alpha * G$ with a symmetric measure $\mu_\alpha (-\bx) = \mu_\alpha (\bx)$.
Thus,
$$
\int_\Omega G (\mu_\alpha * \partial_\beta G) 
  = \int_\Omega (\mu_\alpha * G) \partial_\beta G = - \int_\Omega \big(\partial_\beta (\mu_\alpha * G)\big) G
  = - \int_\Omega (\mu_\alpha * \partial_\beta G) G = 0
$$
for all $\alpha$, $\beta$, and any function $G$.
This allows to write a proof that is formally the same as for the continuous interior product \eqref{i_a}.
Thus, using that $(\curl \bC) \times \uvec_\alpha = \partial_\alpha \bC - \grad C_\alpha$
we have
$$
\int_{\Omega} \bC \cdot \grad I^0_{\uvec_\alpha} \bC
= \int_{\Omega} \bC \cdot (\mu_\alpha * \grad C_\alpha)
= \int_{\Omega} \bC \cdot \big(\mu_\alpha * (\partial_\alpha \bC - (\curl \bC)\times \uvec_\alpha)  \big)
= -\int_{\Omega} \bC \cdot I^1_{\uvec_\alpha} \curl \bC
$$
which proves \eqref{I01_rel}. The relation \eqref{I12_rel} follows by a similar argument.
\eproof

\section{The semi-discrete Hamiltonian system as a system of ordinary differential equations}
\label{sec:hamiltonian_matrix}

In this section, we express the variational particle method \eqref{sawf-solve}--\eqref{sawf-push}
as a system of ordinary differential equations.
This will allow us to introduce some useful notation for our general framework,
and to verify the Hamiltonian structure of the semi-discrete system.

\subsection{Commuting diagrams with degrees of freedom} 
\label{sec:pi_geo}

One practical approach to build commuting projection operators is to introduce one additional
layer in the diagram \eqref{CD}, consisting of coefficient spaces $\cC^\ell = \RR^{N_\ell}$
corresponding to the choice of specific bases for the finite-dimensional spaces
$V_h^\ell$ with dimension $N_\ell$.
This approach is somehow parallel to the geometric construction of \cite{kreeft2011mimetic} where
commuting de Rham complexes are described for differential forms. As we consider here a
a finite element setting, we will follow similar principles but our construction does not involve differential forms.
\begin{equation} \label{CD3}
\begin{tikzpicture}[ampersand replacement=\&, baseline] 
\matrix (m) [matrix of math nodes,row sep=3em,column sep=5em,minimum width=2em] {
   ~~ \Vsp^0 ~ \bbb
   \& ~~ \Vsp^1 ~ \bbb
    \& ~~ \Vsp^2 ~ \bbb
      \& ~~ \Vsp^3 ~ \bbb
\\
~~ \cC^0 ~ \bbb
  \& ~~ \cC^1 ~ \bbb
\& ~~ \cC^2 ~ \bbb
\& ~~ \cC^3 ~ \bbb
\\
~~ V_h^0 ~ \bbb
  \& ~~ V_h^1 ~ \bbb
\& ~~ V_h^2 ~ \bbb
\& ~~ V_h^3 ~ \bbb
\\
};
\path[-stealth]
(m-1-1) edge node [above] {$\grad$} (m-1-2)
        edge node [right] {$\arrsigma^0$} (m-2-1)
        edge [bend right=40] node [pos=0.2, left] {$\Pi^0$} (m-3-1)
(m-1-2) edge node [above] {$\curl$} (m-1-3)
        edge node [right] {$\arrsigma^1$} (m-2-2)
        edge [bend right=40] node [pos=0.2, left] {$\Pi^1$} (m-3-2)
(m-1-3) edge node [above] {$\Div$} (m-1-4)
        edge node [right] {$\arrsigma^2$} (m-2-3)
        edge [bend right=40] node [pos=0.2, left] {$\Pi^2$} (m-3-3)
(m-1-4) edge node [right] {$\arrsigma^3$} (m-2-4)
        edge [bend right=40] node [pos=0.2, left] {$\Pi^3$} (m-3-4)
(m-2-1) edge node[auto] {$ \matD^0 $} (m-2-2)
(m-2-2) edge node[auto] {$ \matD^1 $} (m-2-3)
(m-2-3) edge node[auto] {$ \matD^2 $} (m-2-4)
(m-2-1.285) edge node [right] {$\cI^0$} (m-3-1.75)
(m-2-2.285) edge node [right] {$\cI^1$} (m-3-2.75)
(m-2-3.285) edge node [right] {$\cI^2$} (m-3-3.75)
(m-2-4.285) edge node [right] {$\cI^3$} (m-3-4.75)
(m-3-1.105) edge node [pos=0.8, left] {$\arrsigma^0$} (m-2-1.255)
(m-3-2.105) edge node [pos=0.8, left] {$\arrsigma^1$} (m-2-2.255)
(m-3-3.105) edge node [pos=0.8, left] {$\arrsigma^2$} (m-2-3.255)
(m-3-4.105) edge node [pos=0.8, left] {$\arrsigma^3$} (m-2-4.255)
(m-3-1) edge node [above] {$\grad$} (m-3-2)
(m-3-2) edge node [above] {$\curl$} (m-3-3)
(m-3-3) edge node [above] {$\Div$} (m-3-4)
;
\end{tikzpicture}
\end{equation}
In this diagram the main novel ingredient is the degrees of freedom
${\arrsigma}^\ell = (\sigma^\ell_i)_{1 \le i \le N_\ell}$, which must be {\em unisolvent}
for the finite-dimensional spaces $V_h^\ell$ in the usual sense that they must be one-to-one
when restricted to these spaces. 
The spaces $\Vsp^\ell$ then denote the domains of these degrees of freedom, and as above we consider
a conforming discretization in the sense that $V_h^\ell \subset \Vsp^\ell$.
The other discrete entities can then be determined from the degrees of freedom.
\begin{itemize}
  \item The ``interpolation'' operators $\cI^\ell$ are characterized by the right-inverse property
  $\arrsigma^\ell \cI^\ell \arrg = \arrg$ for all $\arrg \in \cC^\ell$.
  In particular, the basis functions $\Lambda^\ell_i \in V_h^\ell$ defined by the usual duality relations
  \begin{equation}\label{dualbasis}
    \sigma^\ell_i(\Lambda^\ell_j) = \delta_{i,j} \qquad \text{ for } ~~ 1 \le i, j \le N_\ell
  \end{equation}
  correspond to $\Lambda^\ell_i = \cI^\ell \arre^\ell_i$ where $\arre^\ell_i = (\delta_{i,j})_{1 \le j \le N_\ell}$
  is a canonical basis vector of $\cC^\ell$.
  It is sometimes convenient to stack the basis functions into colum vectors
  $\arrLambda^\ell = (\Lambda^\ell_i)_{1 \le i \le N_\ell}$, and to use a matrix notation
  for stacked functionals evaluated on vectors of functions.
  With this convention, the duality relation \eqref{dualbasis} reads
  \begin{equation} \label{dual_mat}
    \arrsigma^\ell(\arrLambda^\ell) = \matI_{N_\ell}
    \quad \text{ with } \quad
      \arrsigma^\ell(\arrLambda^\ell)
      = \big(\sigma^\ell_i(\Lambda^\ell_j)\big)_{1 \le i, j \le N_\ell}.
  \end{equation}

  \item
  The matrices $\matD^\ell \in \RR^{N_{\ell+1} \times N_{\ell}}$ correspond to the differential operators
  $d^0 = \grad$, $d^1 = \curl$ and $d^2 = \Div$ in the respective bases, namely
  \begin{equation} \label{matD}
    \matD^\ell = \arrsigma^{\ell+1}(d^\ell \arrLambda^\ell)
        = \big(\sigma^{\ell+1}_i(d^\ell \Lambda^\ell_j)\big)_{1 \le i \le N_{\ell+1}, 1 \le j \le N_\ell}
  \end{equation}
  so that we have $\arrsigma^{\ell+1}(d^\ell G) = \arrsigma^{\ell+1}(\arrg^\top d^\ell \arrLambda^\ell) = \matD^\ell \arrg$
  for all $G = \arrg^\top \arrLambda^\ell \in V^\ell_h$ with $\arrg \in \cC^\ell$.
  \item
  The projection operators are defined as $\Pi^\ell = \cI^\ell \arrsigma^\ell : G \to \sum_i \sigma_i^\ell(G) \Lambda^\ell_i$,
  that is, 
  \begin{equation} \label{sigmaPi}
    \Pi^\ell G := (\arrsigma^\ell(G))^\top \arrLambda^\ell \qquad \text{ for } ~ G \in \Vsp^\ell,
  \end{equation}
  and they are characterized by the relations
  \begin{equation} \label{charPi}
  \sigma^\ell_i(\Pi^\ell G) = \sigma^\ell_i(G) \qquad \text{ for } \quad 1 \le i \le N_\ell,
  \end{equation}
  indeed we have
  $
  \arrsigma^\ell(\Pi^\ell G) = \arrsigma^\ell\big( (\arrsigma^\ell(G))^\top \arrLambda^\ell\big)
    = \arrsigma^\ell\big( \arrLambda^\ell\big) \arrsigma^\ell(G) = \arrsigma^\ell(G)
  $
  for all $G \in \Vsp^\ell$.
\end{itemize}
This setting proves particularly useful in practice,
as it allows to restate the commuting diagram properties
\eqref{CD} as a {\em linear relation} between degrees of freedom.
\begin{lemma} \label{lem:CDeq}
  The following properties are equivalent:
  \begin{itemize}
    \item[(i)]
    the projection operators \eqref{sigmaPi} satisfy the commuting diagram properties \eqref{CD},
    \begin{equation} \label{CDPi}
      \Pi^{\ell+1} d^\ell G = d^\ell \Pi^\ell G  \qquad \text{ for all } ~ G \in \Vsp^\ell,
    \end{equation}
    \item[(ii)] there exists a matrix $\matD^\ell \in \RR^{N_{\ell+1} \times N_{\ell}}$ such that
    \begin{equation} \label{CDsigma}
      \arrsigma^{\ell+1} (d^\ell G) = \matD^\ell \arrsigma^\ell(G) \qquad \text{ for all } ~ G \in \Vsp^\ell.
    \end{equation}
  \end{itemize}
  Moreover if \eqref{CDsigma} holds, then the matrix $\matD^\ell$ coincides with \eqref{matD}.
\end{lemma}

\bproof
The proof is a matter of elementary computations. 
For instance, \eqref{CDsigma} yields 
$$
\arrsigma^{\ell+1}(d^\ell \Pi^\ell G)
  = \matD^\ell \arrsigma^{\ell}(\Pi^\ell G)
  = \matD^\ell \arrsigma^{\ell}(G)
  = \arrsigma^{\ell+1}(d^\ell G)
  = \arrsigma^{\ell+1}(\Pi^{\ell+1} d^\ell G)
$$
where we have used twice the characterization \eqref{charPi}.

\eproof

\subsection{The semi-discrete Hamiltonian system in matrix form}

The introduction of a third layer in the commuting diagram offers the possibility
to rewrite the semi-discrete scheme \eqref{sawf-solve}--\eqref{sawf-push} as
a system of ordinary differential equations in matrix form.
To do so we collect all the dynamic variables in a global vector
\begin{equation*}
  \arrU = \begin{pmatrix}
  \arrX \\ \arrV \\ \arrE \\ \arrB
  \end{pmatrix}
\end{equation*}
where the (column) block-vectors $\arrX = \arrX_N = (\px_p)_{p = 1, \cdots, N}$
and $\arrV = \arrV_N = (\pv_p)_{p = 1, \cdots, N}\in (\RR^3)^N$ collect
all the particle positions and velocities as in Section~\ref{sec:DAP},
while the vectors $\arrE = \arrsigma^2(\bE_h) \in \RR^{N_2}$
and $\arrB = \arrsigma^1(\bB_h) \in \RR^{N_1}$ collect the coefficients
of the electric and magnetic fields in their respective bases. 
Using these degrees of freedom, we observe that the coupling fields \eqref{sawf-fields} read
$$
\left\{
\begin{aligned}
  &E^S_\alpha(\px_p) = \int_{\Omega} \bE_h(\bx) \cdot \Pi^2(\uvec_\alpha S_{\px_p})(\bx) \dd \bx
    = \sum_{i,j = 1}^{N_2} \sigma^2_i(\bE_h) \matM^2_{i,j} \sigma^2_j(\uvec_\alpha S_{\px_p})
  \qquad
  \\
  &B^S_\alpha(\px_p) = \int_{\Omega} \bB_h (\bx) \cdot \Pi^1 (\uvec_\alpha S_{\px_p})(\bx) \dd \bx
    = \sum_{i,j = 1}^{N_1} \sigma^1_i(\bB_h) \matM^1_{i,j} \sigma^1_j(\uvec_\alpha S_{\px_p})
\end{aligned}
\right.
$$
where $\matM^\ell$ is the standard finite-element mass matrix in the corresponding basis of $V_h^\ell$, $\ell =1, 2$,
\begin{equation} \label{matM}
  \matM^\ell_{i,j} = \int_\Omega \bLambda^\ell_i(\bx) \cdot \bLambda^\ell_j (\bx) \dd \bx, \qquad 1\le i,j\le N_\ell.
\end{equation}
The value of the coupling fields at the particle positions may then be expressed as block-vectors,
\begin{equation}\label{EBSX}
  \bE^S(\arrX) = \matS^2(\arrX) \matM^2 \arrE
  \qquad \text{ and } \qquad
  \bB^S(\arrX) = \matS^1(\arrX) \matM^1 \arrB \qquad \text{ in } (\RR^{3})^N,
\end{equation}
where $\matS^\ell(\arrX) \in (\RR^{3})^{N \times N_\ell}$ denotes the matrix with $(3\times 1)$ blocks
\begin{equation} \label{matS}
  \matS^\ell(\arrX)_{p,i} = \begin{pmatrix}
  \sigma^\ell_i(\uvec_1 S_{\px_p} )  &\sigma^\ell_i(\uvec_2  S_{\px_p} ) & \sigma^\ell_i(\uvec_3  S_{\px_p} )
  \end{pmatrix}^\top
  \qquad \text{ for } \quad 1 \le p \le N, \quad 1 \le i \le N_\ell.
\end{equation}
We finally let $\matr(\bb) = \big((\uvec_\alpha \times \uvec_\beta ) \cdot \bb \big)_{1 \le \alpha, \beta \le 3}
\in \RR^{3\times 3}$ be the rotation matrix
\begin{equation}\label{matr}
\matr(\bb) = \begin{pmatrix}
   0 & b_3 & - b_2\\
   -b_3 & 0 & b_1 \\
   b_2 & -b_1 & 0
  \end{pmatrix}
  \quad \text{ such that }
  \quad \bv \times \bb  = \matr(\bb) \bv \quad \text{ for all } \bv, \bb \in \RR^3,
\end{equation}
and we denote by
$\matR(\bb(\arrX)) \in (\RR^{3\times 3})^{N \times N}$
the block-diagonal rotation matrix with blocks
\begin{equation}\label{matR}
  \matR(\bb(\arrX))_{p,p} = \matr(\bb(\px_p)).
\end{equation}
Then the particle trajectory equations \eqref{sawf-push}--\eqref{sawf-fields},
\begin{equation} \label{sawf-push-fields}
\left\{
\begin{aligned}
  &\frac{\dd \px_p}{\dd t} = \pv_p
  \\
  &
  \frac{\dd \pv_p}{\dd t} = \frac{q_p}{m_p}
    \big(\bE^S(\px_p) + \pv_p \times \bB^S(\px_p) \big)
\end{aligned}
\right.
\qquad \text{ for }  p = 1, \dots, N
\end{equation}
can be written in the block-matrix form  
\begin{equation}\label{XVmat}
  \left\{
  \begin{aligned}
    \frac{\dd \arrX}{\dd t} &= \arrV,\\
    \frac{\dd \arrV}{\dd t} &= \matW_{\frac qm} \big( \matS^2(\arrX) \matM^2 \arrE + \matR^1(\arrX,\arrB)\arrV \big) ,
  \end{aligned}
  \right.
\end{equation}
where
$\matW_{\frac qm} = \diag(\frac{q_p}{m_p} : p \in \range{1}{N})$
is the diagonal weighting matrix carrying the particles charge to mass ratios,
and where we have denoted
\begin{equation}\label{matR1}
  \matR^1(\arrX,\arrB) = \matR(\bB^S(\arrX)) = \matR(\matS^1(\arrX)\matM^1 \arrB ) \qquad \text{ in } (\RR^{3\times 3})^{N\times N}
\end{equation}
the block-diagonal rotation matrix associated with the coupling magnetic field. Observe that its diagonal blocks
read $
\matR^1(\arrX,\arrB)_{p,p}
  = \big( \bsigma^1(\uvec_\alpha \times \uvec_\beta S_{\px_p})^\top \matM^1 \arrB \big)_{1\le \alpha, \beta \le 3}.
$

Turning to the field equations \eqref{sawf-solve},
we see that the strong \Ampere{} equation can be expressed directly on the degrees of freedom
$\arrsigma^2$. From the characterization of the projection operator \eqref{charPi}
we have
$
\sigma^2_i(\Pi^2 \bJ^S_N) = \sum_{p = 1 \cdots N} q_p \sigma^2_i(\pv_p S_{\px_p}),
$
hence our \Ampere{} equation takes the form
\begin{equation} \label{Emat}
  \frac{\dd \arrE}{\dd t} - \matC \arrB = - \matS^2(\arrX)^\top \matW_q \arrV
\end{equation}
with $\matC = \matD^1$ the matrix of the operator $\curl: V_h^1 \to V_h^2$, see \eqref{matD},
$\matS^2(\arrX)$ the matrix defined in \eqref{matS} and
$\matW_q$ the diagonal weighting matrix carrying the particles charges.
Finally the weak Faraday equation is tested against the basis functions $\bLambda^1_i$.
By definition of the weak curl operator \eqref{dw} this yields
\begin{equation} \label{dtBmat}
  \matM^1\frac{\dd \arrB}{\dd t} + \matC^\top \matM^2 \arrE = 0
\end{equation}
with $\matM^1$ and $\matM^2$ the mass matrices recalled in \eqref{matM}.

Finally, rewriting the discrete Hamiltonian $\cH_h(\arrX_N, \arrV_N, \bE_h, \bB_h)$ as a function of the array variables
\begin{equation}
\sfH(\arrU) = \tfrac 12 \arrV^\top \matW_m \arrV + \tfrac 12 \arrE^\top \matM^2 \arrE + \tfrac 12 \arrB^\top \matM^1 \arrB,
\end{equation}
with $\matW_m$ the diagonal weighting matrix carrying the particle masses, see~\eqref{cH},
we obtain for the corresponding derivatives
$$
\nabla_{\arrU} \sfH(\arrU) = \begin{pmatrix}
  \nabla_{\arrX} \sfH
  \\
  \nabla_{\arrV} \sfH
  \\
  \nabla_{\arrE} \sfH
  \\
  \nabla_{\arrB} \sfH
\end{pmatrix}(\arrU)
= \begin{pmatrix}
  \arr{0}
  \\
  \matW_m \arrV
  \\
  \matM^2 \arrE
  \\
  \matM^1 \arrB
  \end{pmatrix}
$$
which allows us to rewrite the equations \eqref{sawf-solve}--\eqref{sawf-fields} as a system of ODEs
\begin{equation} \label{ham}
  \frac{\dd \arrU}{\dd t} = \matJ(\arrU)\nabla_\arrU \sfH(\arrU)
\end{equation}
with a structure matrix given by
\begin{equation}\label{eq:poisson-mat}
 \matJ(\arrU) =
 \begin{pmatrix}
   0 &  \matW_{\frac 1m} & 0& 0 \\
  -\matW_{\frac 1m} & \matW_{\frac qm} \matR^1(\arrX, \arrB) \matW_{\frac 1m} & \matW_{\frac qm} \matS^2(\arrX) & 0\\
  0 & -\matS^2(\arrX)^\top\matW_{\frac qm} & 0 & \matC (\matM^1)^{-1} \\
  0& 0 & - (\matM^1)^{-1}  \matC^\top &0\\
  \end{pmatrix}.
\end{equation}
In particular, System \eqref{ham} may be rewritten in the form of a Poisson system
\begin{equation} \label{hamP}
\frac{\dd \arrU}{\dd t} = \{\arrU,\sfH\}(\arrU)
\end{equation}
with a discrete bracket defined as
$\{ \sfF, \sfG \}(\arrU) := (\nabla_\arrU \sfF)^\top \matJ(\arrU) \nabla_{\arrU} \sfG$, that is,
\begin{equation} \label{bracket}
  \begin{aligned}
  \{ \sfF, \sfG\}(\arrU)
  &= (\nabla_\arrX \sfF)^\top \matW_{\frac 1m} \nabla_\arrV \sfG - (\nabla_\arrV \sfF)^\top \matW_{\frac 1m} \nabla_\arrX \sfG
  \\
  &\qquad +
    (\nabla_\arrV \sfF)^\top\matW_{\frac qm}\matR^1(\arrX,\arrB)\matW_{\frac 1m}\nabla_\arrV \sfG
  \\
  &\qquad \qquad +
  ( \nabla_\arrV \sfF)^\top \matW_{\frac qm}\matS^2(\arrX)\nabla_\arrE \sfG
      - (\nabla_\arrE \sfF)^\top \matS^2(\arrX)^\top \matW_{\frac qm}\nabla_\arrV \sfG
  \\
  &\qquad \qquad \qquad +
  (\nabla_\arrE \sfF)^\top \matC(\matM^1)^{-1}  \nabla_\arrB \sfG - (\nabla_\arrB \sfF)^\top(\matM^1)^{-1}  \matC^\top\nabla_\arrE \sfG.
  \end{aligned}
\end{equation}
Note that this is just the matrix form of the discrete bracket $\{\cF_h,\cG_h\}$ given in \eqref{eq:disc_bracket},
where $\cF_h$ and $\cG_h$ are the same functionals as $\sfF$ and $\sfG$ but seen as functions
of the finite element fields $\bE_h, \bB_h$.

Compared with the Poisson matrix found in \cite[Eq.~(4.29)]{kraus2016gempic},
we observe that the main difference lies in the fact that the particle-field
coupling blocks now involve the degrees of freedom of the smoothed particles
through the matrices $\matS^2$ and $\matS^1$ which involve
the generic commuting diagram operators $\Pi^2$ and $\Pi^1$, see \eqref{matS} and \eqref{matR1}.
In particular, the similarity of both matrices allows us to easily verify the Poisson
structure of the semi-discrete system \eqref{sawf-solve}--\eqref{sawf-fields}.

\subsection{Proof of Theorem~\ref{th:Poisson}}
\label{sec:proof}
Since we have rewritten our equations in a matrix form, it suffices to show that
$\matJ$ is a Poisson matrix in the sense of \cite[Def.~VII.2.4]{HairerLubichWanner:2006},
i.e., that it is skew-symmetric and it satisfies the matrix Jacobi identity.
This will show that \eqref{bracket} is a (non-canonical) Poisson bracket and that
\eqref{hamP}, namely \eqref{ham}, is a Poisson system.

Using that weighting matrices like $\matW_{\frac qm}$ are diagonal,
and that $\matR^1(\arrX,\arrB)$ is skew-symmetric,
we easily verify that $\matJ = - \matJ^T$.
To verify the matrix Jacobi identity, we then observe that $\matJ$ has the same form
as the one involved in the original GEMPIC scheme, see \cite[Eq.~(4.29)]{kraus2016gempic}, with
$\matC\matM_1^{-1}$ and $\matLambda^1(\arrX)\matM_1^{-1}$ replaced by $(\matM^1)^{-1}\matC^\top$
and $\matS^2(\arrX)$, respectively (the mass and curl matrices being defined for different spaces,
due to the different ansatz in the fields).
We also note that $\matR^1$ plays the role of the magnetic rotation matrix $\matB$ in \cite{kraus2016gempic},
with smoothed coupling terms as already observed.
In particular, we may follow the same reasoning to verify that it satisfies the Jacobi identity,
which amounts to verifying that the analog of Eqs.~(4.34) and (4.38) hold in our case.
Using the block-diagonal matrix $\matR^1(\arrX,\arrB)$ defined by \eqref{matR1},
and taking $(p,\alpha)$, $(p,\beta)$, and $(p,\gamma)$ as multi-indices
corresponding to $b$, $c$, and $d$, Equation (4.34) becomes (for $q_p \neq 0$)
\begin{equation} \label{Jacobi1}
\frac{\partial \matR^1(\arrX,\arrB)_{(p,\alpha),(p,\beta)}}{\partial x_{p,\gamma}}
+ \frac{\partial \matR^1(\arrX,\arrB)_{(p,\beta),(p,\gamma)}}{\partial x_{p,\alpha}}
  + \frac{\partial \matR^1(\arrX,\arrB)_{(p,\gamma),(p,\alpha)}}{\partial x_{p,\beta}}
    = 0
    \qquad \forall p, \alpha, \beta, \gamma.
\end{equation}
Using the expression
$
\matR^1(\arrX,\arrB)_{(p,\alpha),(p,\beta)}
  = \bsigma^1(\uvec_\alpha \times \uvec_\beta S_{\px_p})^\top \matM^1 \arrB
$
seen above, this amounts to
$$
\arrB^\top \matM^1 \bsigma^1\Big(
  \uvec_\alpha \times \uvec_\beta (\partial_\gamma S_{\px_p})
    + \uvec_\beta \times \uvec_\gamma (\partial_\alpha S_{\px_p})
      + \uvec_\gamma \times \uvec_\alpha (\partial_\beta S_{\px_p})
  \Big)
    = 0.
$$
By antisymmetry, we see that the function in parentheses vanishes if two of the components coincide,
so that we may assume w.l.o.g. that $(\alpha, \beta, \gamma) = (1,2,3)$. Then this
function is just $\nabla S_{\px_p}$ and the above equation amounts to
$$
0 = \arrB^\top \matM^1 \bsigma^1\big(\nabla S_{\px_p})
    = \int_\Omega \bB_h \cdot \Pi^1(\nabla S_{\px_p}) = \int_\Omega \bB_h \cdot \nabla \Pi^0(S_{\px_p})
    = - \int_\Omega (\Div_w \bB_h) \Pi^0(S_{\px_p}),
$$
where we have used the commuting diagram property and the admissibility
of the shape function $S$. The desired equality \eqref{Jacobi1} then follows
from the discrete magnetic Gauss law, see \eqref{sawf-Gauss}.

The second equality to verify is the analog of Equation (4.37) from \cite{kraus2016gempic},
which reads here (given the above matrix correspondence and correcting a typo on the sign of the right-hand side)
\begin{equation} \label{Jacobi2}
\frac{\partial \matS^2(\arrX)_{(p,\alpha),i}}{\partial x_{p,\beta}}
- \frac{\partial \matS^2(\arrX)_{(p,\beta),i}}{\partial x_{p,\alpha}}
= - \sum_{j = 1}^{N_1} \frac{\partial \matR^1(\arrX,\arrB)_{(p,\alpha),(p,\beta)}}{\partial {\mathsf B}_j}
        \big((\matM^1)^{-1}\matC^\top\big)_{j,i}
        \qquad \forall p, \alpha, \beta, i.
\end{equation}
By antisymmetry of $\matR^1$, we see that both sides vanish for $\alpha = \beta$,
so let us assume w.l.o.g. that $(\alpha,\beta) = (1,2)$. Then
$\matR^1(\arrX,\arrB)_{(p,\alpha),(p,\beta)} = \bsigma^1(\uvec_3 S_{\px_p})^\top \matM^1 \arrB$
and by differentiating these entries and those of the matrix $\matS^2$, see \eqref{matS}, the equality becomes
$$ 
\sigma^2_i\big(\uvec_1 (\partial_2 S_{\px_p}) \big)
  - \sigma^2_i\big(\uvec_2 (\partial_1 S_{\px_p}) \big)
  = \big(\matC \bsigma^1(\uvec_3 S_{\px_p})\big)_{i} \qquad \text{ for } i = 1, \dots, N_1.
$$ 
In vector terms this writes
$\bsigma^2\big(\curl (\uvec_3 S_{\px_p}) \big)
  = \matC \bsigma^1(\uvec_3 S_{\px_p}),
$
which directly follows from the commuting diagram property as seen in Lemma~\ref{lem:CDeq}. Thus
\eqref{Jacobi2} holds, which shows that $\matJ$ satisfies the Jacobi identity and is indeed a Poisson matrix. 
\eproof

\subsection{Propagation in time}
\label{sec:time}

Based on its Poisson structure, geometric time propagation schemes can be derived for our variational system in the same way in \cite{kraus2016gempic}. More precisely, a variational integrator can be derived from a Hamiltonian splitting, that yields a scheme that is explicit in time. We refer to~\cite[Sec.~5.1]{kraus2016gempic} where the resulting equations are detailed for the weak \Ampere{} case and delta shape functions. This kind of splitting has originally been proposed for the Vlasov--Maxwell system in  \cite{He:2015,Xiao:2015} as a Hamiltonian splitting and later been constructed from a fully discrete action principle in \cite{Xiao:2018}. On the other hand, energy-conserving time propagators can be derived by an antisymmetric splitting of the Poisson matrix combined with a suitable discrete-gradient time propagation of the substeps as explained in \cite{kormann2020energy}. In our numerical experiments, we consider the energy- and Gauss-conserving discrete-gradient method from \cite{kormann2020energy}, which demonstrates the best the conservation properties of the phase-space discretization, and the Hamiltonian splitting from \cite{kraus2016gempic} due to its simplicity.

\section{Generalization and application to the strong Faraday model} \label{sec:SF}

Before turning to the description of particular discretizations of Maxwell's equations,
it may be useful to pause for a moment and make some comments on the above findings.
In our variational derivation we have explicitly required that \eqref{CD} was a commuting diagram,
and by doing so we have made two implicitly assumptions: first, we have considered that the discrete sequence
involved {\em strong} differential operators, which corresponds to a conforming discretization. Second, we have
referred to the operators $\Pi^\ell$ as {\em projection} operators. Although these are standard properties to assume,
they played no particular role in our analysis, be it in the variational derivation of Section~\ref{sec:variational},
or in the proof of its Hamiltonian structure. In particular, our results directly apply to a more general setting
of the form
\begin{equation}\label{tCD}
\begin{tikzpicture}[baseline=(current  bounding  box.center)]
  \matrix (m) [matrix of math nodes,row sep=3em,column sep=4em,minimum width=2em] {
    \bbb \Vsp^0 \bbb
      & \bbb \Vsp^1 \bbb
          & \bbb \Vsp^2 \bbb%
              & \bbb \Vsp^3 \bbb
    \\
           \wt V_h^0 & \wt V_h^1 & \wt V_h^2 & \wt V_h^3
    \\
    };
  \path[-stealth]
  (m-1-1) edge node [above] {$\grad$} (m-1-2)
  (m-1-2) edge node [above] {$\curl$} (m-1-3)
  (m-1-3) edge node [above] {$\Div$}  (m-1-4)
  (m-1-1) edge node [left]  {$\wt\Pi^0$} (m-2-1)
  (m-1-2) edge node [left]  {$\wt\Pi^1$} (m-2-2)
  (m-1-3) edge node [left]  {$\wt\Pi^2$} (m-2-3)
  (m-1-4) edge node [left]  {$\wt\Pi^3$} (m-2-4)
  (m-2-1.15) edge node [above] {$\wt\grad$} (m-2-2.165)
  (m-2-2.15) edge node [above] {$\wt\curl$} (m-2-3.165)
  (m-2-3.15) edge node [above] {$\wt\Div$}  (m-2-4.165)
  (m-2-4.195) edge [dashed] node [below] {$\wt\grad_*$} (m-2-3.345)
  (m-2-3.195) edge [dashed] node [below] {$\wt\curl_*$} (m-2-2.345)
  (m-2-2.195) edge [dashed] node [below] {$\wt\Div_*$}  (m-2-1.345)
  ;
\end{tikzpicture}
\end{equation}
where the discrete differential operators $\wt\grad$, $\wt\curl$, $\wt\Div$
no longer need to coincide with the exact ones
(in particular, the discrete spaces $\wt V^\ell_h$ need not be
conforming in $H^1$, $H(\curl)$ and $H(\Div)$), and the $\wt \Pi^\ell$
no longer need to be projection operators. In this generalized setting the only assumptions are that:
\begin{itemize}
  \item[(i)] the solid diagram in \eqref{tCD} commutes,
  \item[(ii)] the lower discrete differential are adjoint to the upper ones in the sense of \eqref{dw}, namely
  $
  \int_\Omega \t \vp_h \wt \Div_* \t \bF_h = -\int_\Omega (\wt \grad \t \vp_h)\cdot \t \bF_h
  $
  must hold for all $\t \vp_h \in \wt V^0_h$ and $\t \bF_h \in \wt V^1_h$, and so on.
\end{itemize}
Our variational derivation then applies verbatim, starting from the discrete Lagrangian
\begin{multline*} 
\wt\cL_h(\arrX_N, \arrX'_N, \arrV_N, \t\bA_h, \t\bA'_h, \t \phi_h)
= \sum_{p=1}^N \left( \big(m_p \bV_p + q_p \tilde \bA^S(\bX_p)\big)\cdot \bX'_p
  - \Big(\frac{m_p}{2}\bV^2_p +q_p\t \phi^S(\bX_p)\Big)\right)
  \\
  +\frac {1}{2} \int_{\Omega} |\wt \grad_* \t\phi_h(\bx) + \t\bA'_h(\bx)|^2 \dd \bx  - \frac {1}{2}\int_{\Omega} |\wt \curl_* \t\bA_h(\bx)|^2  \dd \bx
\end{multline*}
with particle arrays $\arrX_N, \arrX'_N, \arrV_N \in (\RR^3)^N$,
discrete fields $\t\bA_h, \t\bA'_h \in \wt V_h^2$, $\t\phi_h \in \wt V_h^3$
and coupling potentials defined as in \eqref{AphiS}.
The resulting variational equations, analog to \eqref{sawf-solve}--\eqref{sawf-Gauss},
read
\begin{equation} \label{sfwa-solve}
  \left\{
  \begin{aligned}
    - &\dt \t \bE_h + \wt \curl \t \bB_h = \wt\Pi^2 \bJ^S_N
    \\
    &\dt \t \bB_h + \wt \curl_* \t \bE_h = 0
  \end{aligned}
  \right.
  \qquad \text {and } \qquad
  \left\{
  \begin{aligned}
    &\frac{\dd \px_p}{\dd t}  = \pv_p
    \\
    &
    \frac{\dd \pv_p}{\dd t}  = \frac{q_p}{m_p} \big(\t \bE^S(\px_p) + \pv_p \times \t \bB^S(\px_p)\big)
  \end{aligned}
  \right.
\end{equation}
with coupling fields defined similarly as in \eqref{sawf-fields}.
Our analysis then shows that these general equations preserve both the corresponding
discrete Gauss laws and the Hamiltonian, and that they have a discrete Poisson structure.
This allows to extend our results to a wider range of discrete settings,
including 
the structure-preserving DG-type Conga discretizations developed in
\cite{Campos-Pinto.Sonnendrucker.2016.mcomp,campos_pinto_compatible_2017_I}
where both $\bE$ and $\bB$ are represented in broken finite element spaces.
Our results also apply to the discrete ansatz \eqref{fields_SF}
corresponding to a strong Faraday equation. For this case we may consider a conforming (strong)
discretization of the form \eqref{CD}, and set 
\begin{equation} \label{tVh}
  \wt V_h^0 := V_h^3, \qquad
  \wt V_h^1 := V_h^2, \qquad
  \wt V_h^2 := V_h^1, \qquad
  \wt V_h^3 := V_h^0
\end{equation}
so that the ansatz \eqref{fields_SF} takes a form similar to the one \eqref{fields_SA}
considered above, namely
\begin{equation}\label{tfields_SF2}
\t \bB_h \in \wt V_h^1
\xrightarrow{ \mbox{$~ \wt \grad ~$}}
\t \bE_h, \t \bA_h \in \wt V_h^2
\xrightarrow{ \mbox{$~ \wt \curl ~$}}
\t \phi_h \in \wt V_h^3.
\end{equation}
A commuting diagram \eqref{tCD} involving the spaces \eqref{tVh} can then be obtained as follows:
define the commuting (upper) discrete differential operators as the
weak operators \eqref{dw}, i.e.
$$
\wt \grad := \grad_w, \qquad
\wt \curl := \curl_w, \qquad
\wt \Div := \Div_w,
$$
use the strong ones for the adjoint (lower) operators,
$$
\wt \grad_* := \grad, \qquad
\wt \curl_* := \curl, \qquad
\wt \Div_* := \Div.
$$
and for the projection operators $\wt \Pi^\ell$
simply take the $L^2$ projections on the discrete spaces,
$$
\sprod{\wt \Pi^\ell G}{G_h} = \sprod{G}{G_h}
\qquad \text{ for } G \in V^\ell, ~~  G_h \in V^\ell_h.
$$
The commutation property is indeed easily verified: For the $\grad$ operator,
using the embedding $\Div : \wt V_h^1 = V_h^2 \to V_h^3 = \wt V_h^0$ and the characterization of
$L^2$ projections, we can write
$$
\sprod{\wt \Pi^1 \grad \psi}{\t \bC_h} =  \sprod{\grad \psi}{\t \bC_h} = -\sprod{\psi}{\Div \t \bC_h}
= - \sprod{\wt \Pi^0 \psi}{\Div \t \bC_h} = \sprod{\wt \grad \wt \Pi^0 \psi}{\t \bC_h}
$$
for all $\psi \in V^0$ and $\t \bC_h \in \wt V_h^1$,
which shows that
$\wt \Pi^1 \grad = \wt \grad \wt \Pi^0$ holds on $\Vsp^0$ (which may be taken here as $H^1(\Omega)$).
The same argument also applies for the operators $\curl$ and $\Div$.
With this construction one recovers the Hamiltonian particle method
of \cite{kraus2016gempic}, with general shape functions.
The discrete Poisson matrix thus takes the same form,
with particle-field coupling terms encoded in block matrices $\matLambda^\ell(\arrX) \in (\RR^3)^{N\times N_\ell}$,
$\ell = 1,2$, with generic $(3 \times 1)$ blocks
\begin{equation} \label{matL}
  \matLambda^\ell(\arrX)_{p,i} = \int_\Omega \bLambda^\ell_i(\bx) \cdot S_{\px_p}(\bx) \dd \bx
  \qquad \text{ for } \quad 1 \le p \le N, \quad 1 \le i \le N_\ell
\end{equation}
which extend the corresponding matrices in \cite{kraus2016gempic} to the case of
a general shape function $S$.

\section{Application to tensor-product spline and Fourier field solvers} 
\label{sec:interpolation_histopolation}

In this section, we apply the above method to the case of tensor-product finite element spaces
defined on cartesian domains.
Following the interpolation / histopolation approach of \cite{Gerritsma.2011.spec, kreeft2011mimetic},
we review a general method for designing commuting diagrams, which is based on geometric degrees of freedom
that can then be associated to finite element spaces of various types. In this article, we detail two
applications, one using splines and another one using truncated Fourier spaces.

\subsection{Geometric degrees of freedom with commuting properties}
\label{sec:geodofs}

Let us equip the cartesian domain $\Omega = [0,L]^3$ 
with a tensor-product grid using $M_\alpha$ nodes along each dimension $\alpha$,
\begin{equation} \label{gx}
\gx_\bsm = (x_{1,m_1}, x_{2,m_2}, x_{3,m_3})
  \qquad \text{ with } ~
  \bsm \in \range{1}{\bM} := \prod_{\alpha = 1}^3 \range{1}{M_\alpha}.
\end{equation}
On this mesh, we consider evaluation functionals defined on the various geometric elements:
\begin{itemize}
  \item point evaluations on the nodes
  \begin{equation} \label{scrP}
    \scrP_\bsm(G) := G(\bx_\bsm),
  \end{equation}

  \item edge integrals along some dimension $1\le \alpha\le 3$,
  \begin{equation} \label{scrE}
  \scrE_{\alpha,\bsm}(G) := \int_{\tte_{\alpha, \bsm}} G
  \quad \text{ with } \quad \tte_{\alpha,\bsm} = [\gx_{\bsm-\be_\alpha}, \gx_\bsm],
  \end{equation}

  \item face integrals normal to some dimension $1\le \alpha\le 3$,
  \begin{equation} \label{scrF}
  \scrF_{\alpha,\bsm}(G) := \int_{\ttf_{\alpha,\bsm}} G
  \quad \text{ with } \quad \ttf_{\alpha, \bsm} = [\tte_{\alpha+1, \bsm-\be_{\alpha-1}}, \tte_{\alpha+1,\bsm}],
  \end{equation}

  \item and cell integrals
  \begin{equation} \label{scrC}
    \scrC_{\bsm}(G) :=
    \int_{\ttc_{\bsm}} G
    \quad \text{ with } \quad \ttc_{\bsm} = [\ttf_{1,\bsm-\be_{1}}, \ttf_{1,\bsm}],
  \end{equation}
\end{itemize}
where we have denoted by $[a,b]$ the convex hull of $a \cup b$.
A set of ``geometric'' degrees of freedom can then be derived from these local functionals:
\begin{equation}\label{eq:sigma_interp_histop}
\left\{
\begin{aligned}
&\hat \sigma^0_{\bsm}(\vp) := \scrP_\bsm(\vp)
    \quad &&\text{ for } \vp \in \Vsp^0
\\
&\hat \sigma^1_{\alpha, \bsm}(\bC) := \scrE_{\alpha,\bsm}(\bC \cdot \uvec_\alpha)
    \quad &&\text{ for } \bC \in \Vsp^1 
\\
&\hat \sigma^2_{\alpha, \bsm}(\bF) := \scrF_{\alpha,\bsm}(\bF \cdot \uvec_\alpha)
    \quad &&\text{ for } \bF \in \Vsp^2
\\
&\hat \sigma^3_{\bsm}(g) := \scrC_{\bsm}(g)
    \quad &&\text{ for } g \in \Vsp^3
\end{aligned}\right. 
\quad \text{ and for } ~ \alpha \in \range{1}{3}, ~ \bsm \in \range{1}{\bM}.
\end{equation}
\extended{
  In Selalib a different scaling is used for the dofs, namely
  \begin{equation}\label{eq:sigma_interp_histop_selalib}
  \left\{
  \begin{aligned}
  &\check \sigma^0_{\bsm}(\vp) := \scrP_\bsm(\vp)
      \quad &&\text{ for } \vp \in \Vsp^0
  \\
  &\check \sigma^1_{\alpha, \bsm}(\bC) := \frac{1}{h_\alpha}\scrE_{\alpha,\bsm}(\bC \cdot \uvec_\alpha)
      \quad &&\text{ for } \bC \in \Vsp^1 
  \\
  &\check \sigma^2_{\alpha, \bsm}(\bF) := \frac{h_\alpha}{h_1h_2h_3}\scrF_{\alpha,\bsm}(\bF \cdot \uvec_\alpha)
      \quad &&\text{ for } \bF \in \Vsp^2
  \\
  &\check \sigma^3_{\bsm}(g) := \frac{1}{h_1h_2h_3} \scrC_{\bsm}(g)
      \quad &&\text{ for } g \in \Vsp^3
  \end{aligned}\right. 
  \quad \text{ and for } ~ \alpha \in \range{1}{3}, ~ \bsm \in \range{1}{\bM}.
\end{equation}
}
If these degrees of freedom are associated to spaces $V_h^\ell$, $0 \le \ell \le 3$, of respective dimensions
\begin{equation}
  \label{NM}
  N_0 = N_3 = M
\qquad \text{ and } \qquad
N_1 = N_2 = 3 M,
\qquad \text{ with } \qquad M := M_1M_2M_3,
\end{equation}
and for which they are unisolvent, then they
define a unique set of dual basis functions $\hat \Lambda^\ell_i$
according to \eqref{dualbasis}, which may also be called ``geometric'':
for the space $V_h^0$ for example these basis functions correspond to the interpolatory basis
associated with the nodes $\gx_\bsm$, for the space $V_h^3$ they correspond to histopolation
basis functions, and for the intermediate spaces they involve a combination of both.
A key property of this construction is the following.

\begin{lemma}\label{lemma:interpolation_histopolation_commuting_diagram}
The degrees of freedom defined by \eqref{eq:sigma_interp_histop}
are well-defined on the domains
$$
\Vsp^0 = W^{1}_{\rm per, 1,2,3},
\quad
\Vsp^1 = W^{1}_{\rm per,2,3} \times W^{1}_{\rm per,3,1} \times W^{1}_{\rm per,1,2},
\quad
\Vsp^2 = W^{1}_{\rm per,1} \times W^{1}_{\rm per,2} \times W^{1}_{\rm per,3},
\quad
\Vsp^3 = L^1_{\rm per},
$$
where we have denoted by $L^1_{\rm per}$ the space of $L$-periodic and locally $L^1$ functions, and
by
\begin{equation} \label{W1a}
\left\{\begin{aligned}
&W^{1}_{\rm per, 1,2,3} := \{ G \in L^1_{\rm per} : \partial_1 \partial_2 \partial_3 G \in L^1_{\rm per}\}
\\
&W^{1}_{\rm per, \alpha, \beta } := \{ G \in L^1_{\rm per} : \partial_\alpha \partial_\beta G \in L^1_{\rm per} \}
\\
&W^{1}_{\rm per, \alpha} := \{ G \in L^1_{\rm per} : \partial_\alpha G \in L^1_{\rm per}\}
\end{aligned}
\right.
\end{equation}
anisotropic Sobolev spaces of $W^{s,1}$ type.
Moreover if the $\hat \arrsigma^\ell$ are unisolvent on the spaces $V_h^\ell$, then the resulting projection operators
$\hat \Pi^\ell$ characterized by the relations \eqref{charPi}, namely
$$
\left\{
\begin{aligned}
&\hat \sigma^0_{\bsm}(\hat \Pi^0\vp) = \hat \sigma^0_{\bsm}(\vp)
\\
&\hat \sigma^1_{\alpha, \bsm}(\hat \Pi^1\bC) = \hat \sigma^1_{\alpha, \bsm}(\bC)
\\
&\hat \sigma^2_{\alpha, \bsm}(\hat \Pi^2\bF) = \hat \sigma^2_{\alpha, \bsm}(\bF)
\\
&\hat \sigma^3_{\bsm}(\hat \Pi^3 g) = \hat \sigma^3_{\bsm}(g)
\end{aligned}\right.
\qquad \text{ for all } \alpha \in \range{1}{3}, ~ \bsm \in \range{1}{\bM},
$$
satisfy the commuting diagram property
$$
  d^\ell \Pi^\ell G = \Pi^{\ell+1} d^\ell G
  \quad \text{ for all } G \in \Vsp^\ell.
$$
\end{lemma}

\bproof
The fact that these degrees of freedom are well-defined on the above domains follows from standard
Sobolev inequalities, see e.g. \cite[Rem.~13]{Brezis.2010}. The commuting diagram properties
are then easy to verify by applying the Stokes formula and Lemma~\ref{lem:CDeq}. For the gradient for instance,
we consider some $\vp \in \Vsp^0$ and compute
$$
\hat \sigma^1_{\bsm,\alpha}(\grad \vp) = \int_{\tte_{\bsm,\alpha}} \be_\alpha \cdot \grad \vp
= \vp(\gx_\bsm) - \vp(\gx_{\bsm-\be_\alpha}) = \hat \sigma^0_{\bsm}(\vp) - \hat \sigma^0_{\bsm-\be_\alpha}(\vp).
$$
According to Lemma~\ref{lem:CDeq}, this specifies the gradient matrix $\hat \matD^0 \in \RR^{N_1 \times N_0}$ such that
$$
\hat \arrsigma^1(\grad \vp) = \hat \matD^0 \hat \arrsigma^0(\vp)
$$
and also implies $\grad \hat \Pi^0 \vp = \hat \Pi^1 \grad \vp$.
The same argument works for the other operators.
\eproof

In the construction above, we see that the commuting properties rely only on the geometric nature
of the degrees of freedom, and not on the tensor-product structure of the grid.
However, this tensor-product structure allows us to specify the form of the differential matrices.
Setting $\varphi= \hat \Lambda^0_{\bk}$ in the proof of
Lemma \ref{lemma:interpolation_histopolation_commuting_diagram},
we find indeed the following representation of $\hat \matD^0$
\begin{equation} \label{hD0}
\hat \matD^0 = \begin{pmatrix}
\matI_{M_3} \otimes \matI_{M_2} \otimes \matd_1 \\
\matI_{M_3} \otimes \matd_2 \otimes \matI_{M_1} \\
\matd_3 \otimes \matI_{M_2} \otimes \matI_{M_1}
\end{pmatrix},
\end{equation}
where $\matI_{M_\alpha}$ is the identity matrix of size $M_\alpha \times M_\alpha$,
$\matd_{\alpha}$ is a univariate differential matrix
\begin{equation} \label{hda}
\matd_{\alpha} =  \begin{pmatrix}
1 & 0 & \ldots && 0 & -1 \\
-1 & 1 & 0  & & & \\
0 & -1 & 1 & 0 & &\\
&  \ddots & \ddots & \ddots & \ddots &\\
&& 0 & -1 & 1 & 0 \\
&&& 0& -1 & 1 \\
\end{pmatrix} \in \RR^{M_{\alpha} \times M_{\alpha}}, \quad \alpha \in \range{1}{3}
\end{equation}
and the Kronecker matrix product is defined as
$
(\mat{c} \otimes \mat{b} \otimes \mat{a})_{\bsm,\bn} = \mat{c}_{m_3,n_3}\mat{b}_{m_2,n_2}\mat{a}_{m_1,n_1}.
$
In the same way, we find
\begin{equation}\label{hD12}
  \hat \matD^1 = \begin{pmatrix}
  \matO_{M} & - \matd_{3} \otimes \matI_{M_2} \otimes \matI_{M_1}  & \matI_{M_3} \otimes \matd_{2} \otimes \matI_{M_3}\\
  \matd_{3} \otimes \matI_{M_2} \otimes \matI_{M_1}  & \matO_{M} &  -\matI_{M_3} \otimes \matI_{M_2} \otimes \matd_{1}\\
  -\matI_{M_3} \otimes \matd_{2} \otimes \matI_{M_1}  & \matI_{M_3} \otimes \matI_{M_2} \otimes \matd_{1} & \matO_{M} \\
  \end{pmatrix}
  \quad  \text{ and } \quad
  \hat \matD^2 = \left(\hat \matD^0\right)^\top
\end{equation}
where $\matO_{M}$ denotes the zero square matrix of size $M=M_1M_2M_3$.
\extended{With the normalization in SeLaLib, we have a scaling of $\frac{1}{h_\alpha}$ in $\matd_{\alpha}$.}
In practice, the basis functions $\hat \Lambda^\ell_i$ defined by the geometric degrees of freedom according to \eqref{dualbasis}
may not be the most convenient to use, either because they have no simple expression, or because
some other basis $\Lambda^\ell_i$ has better locality properties, or leads to simpler discrete Maxwell equations.
One then needs to determine the coefficients of the geometric projections in this new practical basis,
which amounts to finding degrees of freedom $\sigma^\ell_i$ that are dual to the practical basis functions
and lead to the same projection operator $\Pi^\ell = \hat \Pi^\ell$ as the geometric ones.
Using the stacked vector notation introduced in Section~\ref{sec:pi_geo}
for the geometric basis $\hat \arrLambda^\ell$ and the practical basis $\arrLambda^\ell$,
these new degrees of freedom $\arrsigma^\ell$ are characterized by the relations
$$
(\arrsigma^\ell(G))^\top \arrLambda^\ell = \Pi^\ell G = \hat \Pi^\ell G = (\hat \arrsigma^\ell(G))^\top \hat \arrLambda^\ell
  \qquad \text{ for all} ~ G \in \Vsp^\ell.
$$
Introducing the matrix
$
\matK^\ell = \hat \arrsigma^\ell(\arrLambda^\ell) = \big(\hat \sigma^\ell_m(\Lambda^\ell_k)\big)_{1 \le m, k \le N_\ell}
$
such that $\arrLambda^\ell = (\matK^\ell)^\top \hat \arrLambda^\ell$, this yields
$$
\arrsigma^\ell(G) = (\matK^\ell)^{-1} \hat \arrsigma^\ell(G),
$$
which gives a practical formula for computing the coefficients of the geometric projections in the practical basis.
Accordingly, the differential matrices in this new basis read
$$
\matD^\ell = \arrsigma^{\ell+1}(d^\ell \arrLambda^\ell)
  = \Big( \sum_{n,m} \big(\matK^{\ell+1}\big)^{-1}_{i,n} \hat \sigma^{\ell+1}_n(d^\ell \matK^\ell_{m,j} \hat \Lambda^\ell_m) \Big)_{i,j}
  = (\matK^{\ell+1})^{-1} \hat \matD^\ell \matK^\ell.
$$
Note that $\matK^0$ is a Vandermonde matrix when $\arrLambda^0$ is a monomial basis. For this reason
the matrices $\matK^\ell$ are sometimes referred to as a generalized Vandermonde matrices.

\subsection{Compatible finite elements based on B-splines}
\label{sec:splines}

Compatible finite elements based on splines on a Cartesian grid have been studied
by Buffa, Sangalli, V\'azquez and co-authors, see e.g. \cite{buffa2010isogeometric, Buffa:2011},
and in \cite{kraus2016gempic} they have been used to implement the strong Faraday GEMPIC formulation.
Here we describe how spline spaces can be used in conjunction with the geometric degrees of freedom
described in Section~\ref{sec:geodofs}.

For simplicity, we consider periodic boundaries and regular knot sequences with $M_\alpha$ knots per dimension.
Denoting by $N^{p}_{\alpha,k}$ the univariate B-spline of degree $p$ along $x_\alpha$, associated with the knots
$(k h_\alpha, \dots, (k+p+1) h_\alpha)$ where $h_\alpha = \frac{L}{M_\alpha}$, see e.g. \cite{Schumaker_2007},
the first space in the sequence consists of tensor-product splines of multi-variate degree $(p_1, p_2, p_3)$, namely
$$
V^0_h = \SS_{p_1,p_2,p_3} := \Span\Big( \big\{\Lambda^0_{\bk} : \bk \in \range{1}{\bM} \big\} \Big)
\qquad \text{ with } \qquad
\Lambda^0_{\bk}(\bx) := \prod_{\alpha = 1}^3 N^{p_\alpha}_{\alpha,k_\alpha}(x_\alpha)
$$
and the full sequence reads
$$
V^0_h \xrightarrow{ \mbox{$~ \grad ~$} }
V^1_h = \begin{pmatrix} \SS_{p_1-1,p_2,p_3} \\ \SS_{p_1,p_2-1,p_3} \\ \SS_{p_1,p_2,p_3-1} \end{pmatrix}
\xrightarrow{ \mbox{$~ \curl ~$} }
V^2_h = \begin{pmatrix} \SS_{p_1,p_2-1,p_3-1} \\ \SS_{p_1-1,p_2,p_3-1} \\ \SS_{p_1-1,p_2-1,p_3}\end{pmatrix}
\xrightarrow{ \mbox{$~ \Div ~$} } 
V^3_h = \SS_{p_1-1,p_2-1,p_3-1}.
$$
The fact that this is indeed a sequence follows from the well-known relation 
\begin{equation}\label{eq:spline_recursion}
  \frac{\mathrm{d}}{\mathrm{d} x_\alpha} N^p_{\alpha, k} = \frac{1}{h_\alpha}\big(N^{p-1}_{\alpha,k} - N^{p-1}_{\alpha, k+1}\big).
\end{equation}
Introducing for convenience the scaled $B$-splines along $x_\alpha$,
\begin{equation} \label{Dspline}
  D^p_{\alpha,k} = \frac{1}{h_\alpha} N^{p-1}_{\alpha,k}
\end{equation}
yields a particularly simple formula for the derivative operator in the corresponding basis.
In particular, it makes it convenient to equip the vector-valued spaces $V^1_h$, $V^2_h$ with the basis functions
$$
\bLambda^1_{\alpha, \bk}(\bx) := \uvec_\alpha D^{p_\alpha}_{\alpha,k_\alpha}(x_\alpha)
    \prod_{\beta \neq \alpha} N^{p_\beta}_{\beta,k_\beta}(x_\beta)
\qquad \text{ for } ~ \alpha \in \range{1}{3}, ~ \bk \in \range{1}{\bM},
$$
$$
\bLambda^2_{\alpha, \bk}(\bx) := \uvec_\alpha N^{p_\alpha}_{\alpha,k_\alpha}(x_\alpha)
    \prod_{\beta \neq \alpha} D^{p_\beta}_{\beta,k_\beta}(x_\beta)
\qquad \text{ for } ~ \alpha \in \range{1}{3}, ~ \bk \in \range{1}{\bM},
$$
and the last, scalar-valued space $V^3_h$ with
$$
\Lambda^3_{\bk}(\bx) := \prod_{\alpha =1}^{3} D^{p_\alpha}_{\alpha,k_\alpha}(x_\alpha)
\qquad \text{ for } \bk \in \range{1}{\bM}.
$$
In practice, B-splines are appealing because of their minimal support property,
however they are not dual to the geometric degrees of freedom defined in \eqref{eq:sigma_interp_histop}
so that new degrees of freedom must be computed as described at the end of Section~\ref{sec:geodofs}.
For the nodal degrees of freedom, the change of basis matrix reads
\begin{equation}
 \matK^0_{\bsm, \bk} = \hat \sigma^0_\bsm(\Lambda_\bk^0) = \Lambda_\bk^0(\bx_\bsm)
\end{equation}
and a common choice of interpolation nodes $\bx_\bsm$ consists of Greville points,
which coincide with the knot sequence for regular splines of odd degrees, and with their midpoints for even degrees.
More generally, we observe that $\matK^0$ is invertible as long as the degrees of freedom
$\arrsigma^0$ are unisolvent, which holds iff the grid satisfies the spline interpolation condition,
see e.g. \cite[Th.~4.61]{Schumaker_2007}.
Using the tensor-product structure and the locality of the B-splines, we see that $\matK^0$ is the Kronecker product
of three banded matrices, which are also circulant for regular Greville points. 
Moreover, as B-splines satisfy by construction
\begin{equation}\label{eq:spline_integral}
N^p_{\alpha,k}(x) = \int_{x-h_\alpha}^x D_{\alpha,k}^{p}(y) \dd y,
\end{equation}
see \eqref{eq:spline_recursion}, \eqref{Dspline}, we have
$$
  \hat \sigma^1_{\alpha, \bsm}(\bLambda^1_{\alpha,\bk})
    = \scrE_{\alpha,\bsm}( \bLambda^1_{\alpha,\bk}\cdot \uvec_\alpha )
    = \int_{(m_\alpha-1)h_\alpha}^{m_\alpha h_\alpha} \!\! D^{p_\alpha}_{\alpha,k_\alpha}(x) \dd x
        \prod_{\beta \neq \alpha} N^{p_\beta}_{\beta,k_\beta}(m_\beta h_\beta)
    = \Lambda^0_\bk(\bx_\bsm) = \hat \sigma^0_\bsm(\Lambda_\bk^0)
$$
hence the matrix block $\matK^{1,\alpha} = \big(\matK^1_{(\alpha,\bsm),(\alpha,\bk)}\big)_{\bsm,\bk}$ coincides with $\matK^0$, the other blocks
of $\matK^1$ being clearly zero.
Similarly we find that $\matK^3$ and $\matK^{2,\alpha}$ also coincide with $\matK^0$, so that (with obvious notation)
 %
$$
\left\{
\begin{aligned}
&\arrsigma^0(\vp) =  \left(\matK^0\right)^{-1} \hat \arrsigma(\vp),
\\
&\arrsigma^{1,\alpha}(\bC) = \left(\matK^0\right)^{-1} \hat \arrsigma^{1,\alpha}(\bC),
\\
&\arrsigma^{2,\alpha}(\bF) =\left(\matK^0\right)^{-1} \hat \arrsigma^{2,\alpha}(\bF),
\\
&\arrsigma^{3}(g) = \left(\matK^0 \right)^{-1} \hat \arrsigma^3(g).
\end{aligned}\right. 
$$
\extended{Since we use $N^{p-1}$ instead of $D^p$ in our implementation for the basis functions, we have the same scaling of the coefficients.}
From relation \eqref{eq:spline_recursion} we also see that the one-dimensional derivative matrices -- and hence, every $\matD^\ell$ --
are the same as for the geometric basis.
As for the three-dimensional mass matrices, they are the Kronecker product of the one-dimensional mass matrices
which are a circulant matrices with $2p_{\alpha}+1$ non-zero entries per row in each dimension.
Finally we note that the discrete interior products \eqref{I_a} based on the directional averaging operator \eqref{cA_av}
may be evaluated using the relation \eqref{eq:spline_integral}, writing e.g.
$$
I^0_{\uvec_\alpha} \bLambda^1_{\beta,\bk} (\bx)
= \cA_{h,\alpha} (\bLambda^1_{\beta,\bk}\cdot \uvec_\alpha)(\bx)
= \frac {\delta_{\alpha, \beta} }{2h_\alpha} \int_{x_\alpha-h_\alpha}^{x_\alpha+h_\alpha} \!\! D^{p_\alpha}_{\alpha,k_\alpha}
    \prod_{\gamma \neq \alpha} N^{p_\gamma}_{\gamma,k_\gamma}(x_\gamma)
  = \frac {\delta_{\alpha, \beta} }{2h_\alpha}\left(\Lambda^0_{\bk-\uvec_\alpha} + \Lambda^0_{\bk}\right)(\bx).
$$

\subsection{Compatible finite elements based on Fourier spaces}
\label{sec:fourier}

With periodic boundary conditions, another option is to consider a
sequence of compatible finite elements made of discrete Fourier spaces.
Such spectral elements are very common in particle solvers,
with particle-field interaction usually based on discrete
Fourier transforms and FFT algorithms.
Here we describe a coupling based on the geometric degrees of freedom described
in Section~\ref{sec:geodofs}.
To match the dimensions of the grid, we consider spaces with
$M_\alpha = 2K_\alpha + 1$ modes per dimension, of the form
$$
V^0_h = V^1_h = \Span\Big( \big\{\Lambda^0_{\bk} : \bk \in \range{-\bK}{\bK} \big\} \Big)
\qquad \text{ with } \qquad
\Lambda^0_{\bk}(\bx) := \ee^{\frac{2 \ii \pi \bk\cdot \bx}{L}}
  = \prod_{\alpha = 1}^3 \ee^{\frac{2 \ii \pi k_\alpha x_\alpha}{L}},
$$
where we have denoted $\range{-\bK}{\bK} = \prod_{\alpha=1}^3 \range{-K_\alpha}{K_\alpha}$, and
$$
V^2_h = V^3_h = \Span\Big( \big\{ \bLambda^2_{\alpha \bk} : \alpha \in \range{1}{3}, ~ \bk \in \range{1}{\bM} \big\} \Big)
\qquad \text{ with } \qquad
\bLambda^2_{\alpha, \bk}(\bx) := \uvec_\alpha \Lambda^0_{\bk}(\bx).
$$
These discrete spaces clearly form a de Rham sequence, as the derivative of a Fourier mode
is the same mode up to a complex scaling factor.

One interesting feature of the canonical modal basis is that it leads to diagonal Maxwell equations.
Indeed the differential matrices $\matD^\ell$ have the same simple block and Kronecker-product structure
as \eqref{hD0}--\eqref{hD12}, here with diagonal one-dimensional derivative matrices
\begin{equation*}
\matd^{\alpha} = \frac{2\ii \pi}{L_{\alpha}} \diag\big( - K_\alpha, \cdots, 0, \cdots, K_{\alpha}\big),
\end{equation*}
and the mass matrices are all diagonal due to the orthogonality of the basis functions,
with $\matM^{\ell} = L^3 \matI_{N_\ell}$ for the chosen normalization.

However, as the modal basis is not dual to the geometric degrees of freedom from Section~\ref{sec:geodofs},
we need to determine the proper change of basis formulas in order to apply the geometric
interpolation-histopolation projections $\hat \Pi^\ell$, as we did for the B-splines in the previous section.
To do so, it is convenient to consider regular interpolation nodes $\bx_\bsm = (m_1 h_1, m_2 h_2, m_3 h_3)$,
with $h_\alpha = \frac {L}{M_\alpha}$. The nodal change of basis matrix reads then
$$
\matK^0_{\bsm,\bk} = \hat \sigma^0_\bsm(\Lambda_\bk)
= \Lambda_\bk(\bx_\bsm)
  = \prod_{\alpha =1}^3 \ee^{\frac{2 \ii \pi k_\alpha m_\alpha}{M_\alpha}}
$$
which is a standard DFT matrix as well as its inverse,
$$
(\matK^0)^{-1} = \Big(\frac{1}{M} \prod_{\alpha = 1}^3\ee^{-\frac{2 \ii \pi k_\alpha m_\alpha}{M_\alpha}}\Big)_{\bk,\bsm}
= \frac 1M \big(\matK^0\big)^* =: \matF
$$
where we remind that $M = M_1M_2M_3$, see \eqref{NM}.
The interpolation operator in the modal basis then takes the well-known form
$$
\hat \Pi^0 (\vp) = \arrsigma^0(\vp)^\top \arrLambda^0
\qquad \text{ with } \qquad
\arrsigma^0(\vp) = \matF \hat \arrsigma^0(\vp)
  = \frac{1}{M} \sum_{\bsm \in \range{1}{\bM}} \vp(\bx_\bsm)
      \prod_{\alpha = 1}^3 \ee^{-\frac{2 \ii \pi k_\alpha m_\alpha}{M_\alpha}}.
$$
For the other projections in the sequence we proceed similarly as in Section~\ref{sec:splines}, noting that
\begin{equation*}
\int_{(m_{\alpha}-1)h_{\alpha}}^{m_\alpha h_\alpha} \ee^{\frac{2\ii \pi k_{\alpha} x_{\alpha}}{L}} \dd x_{\alpha}
  = T^\alpha_{k_\alpha} \ee^{\frac{2\ii \pi k_{\alpha} m_{\alpha}}{M_\alpha}}
  \qquad \text{ with } \qquad
  T^\alpha_{k_\alpha} := \begin{cases}
  h_\alpha & \text{if } k_{\alpha}=0,
  \\
  \frac{L}{2 \ii \pi k_\alpha} 
    \left(1-\ee^{-\frac{2 \ii \pi k_\alpha}{M_\alpha}}\right) &\text{else}.
\end{cases}
\end{equation*}
\extended{Since we scale the integrals by the inverse of the interval, the transformation matrix $T^\alpha$ is scaled by a factor $\frac{1}{h_\alpha}$.}
In particular, writing 
$
\matT^\alpha := \diag\big(\matT^\alpha_{\bk, \bk} = T^\alpha_{k_\alpha} : \bk \in \range{-\bK}{\bK}\big)
$
we find $\matK^{1,\alpha} = \matK^0 \matT^\alpha$
for the matrix block
$
\matK^{1,\alpha} = \big(\matK^1_{(\alpha,\bsm),(\alpha,\bk)}\big)_{\bsm,\bk}
$, and similarly 
$\matK^{2,\alpha} = \matK^0 \matT^{\alpha-1}\matT^{\alpha+1}$ and $\matK^{3} = \matK^0 \matT^1\matT^2\matT^3$.
The expression of the different degrees of freedom in the modal Fourier basis reads then
$$
\left\{
\begin{aligned}
&\arrsigma^0(\vp) =  \matF \hat \arrsigma^0(\vp),
\\
&\arrsigma^{1,\alpha}(\bC) = \left(\matT^{\alpha}\right)^{-1} \matF \hat \arrsigma^{1,\alpha}(\bC),
\\
&\arrsigma^{2,\alpha}(\bF) =  \left(\matT^{\alpha-1}\matT^{\alpha+1}\right)^{-1} \matF \hat \arrsigma^{2,\alpha}(\bF),
\\
&\arrsigma^3(g) = \left(\matT^{1}\matT^{2}\matT^{3}\right)^{-1} \matF \hat \arrsigma^{3}(g),
\end{aligned}\right. 
$$
where we note that all the $\matT^\alpha$ matrices are clearly diagonal and invertible.
To apply the discrete interior products \eqref{I_a} based on directional averaging \eqref{cA_av},
we finally need to evaluate
$$
(\cA_{h,\alpha} \Lambda^0_{\bk})(\bx)
= \frac{1}{2h_\alpha} \int_{x_\alpha-h_\alpha}^{x_\alpha+h_\alpha} \ee^{\frac{2\ii \pi \bk\cdot\bx}{L}} \dd \bx
= \sinc\Big(\frac{2\pi k_{\alpha}}{M_\alpha}\Big) \Lambda^0_{\bk}(\bx)
$$
for all $\alpha \in \range{1}{3}$ and $\bk \in \range{-\bK}{\bK}$.

\section{Numerical illustration in reduced phase space} \label{sec:numerical_experiments}

In this section, we will show some numerical results obtained with the proposed schemes in a reduced phase space. All results are obtained with an implementation of the strong \Ampere{} scheme 
within the SeLaLib library.
We study the variational semi-discretization as derived in Section~\ref{sec:variational}---which is energy conserving---as well
as the momentum-preserving semi-discretization as derived in Section~\ref{sec:momentum}.
For the basis of the finite element field solver, both splines and Fourier modes are considered.
The shape function is chosen to be a B-spline of varying degree.

As for the time discretization, we compare a Hamiltonian splitting scheme for both space discretization methods, see Section~\ref{sec:time}.
Only when considering the conservation properties we also provide results for the variational scheme with an energy-conserving discrete gradient time discretization. We use a time step of $\Delta t = 0.05$, the linear solvers use a tolerance of $10^{-15}$ and the nonlinear iterations in the discrete gradient method have a tolerance of $10^{-12}$.

\subsection{Physical model}

For the numerical study we consider a reduced phase space with one periodic spatial and one or two velocity dimensions, namely
$\bx=x_1 \in [0,L_1)$, $\bv=(v_1,v_2) \in \RR^2$, with unknowns of the form
$$
f = f(t,x_1,v_1,v_2),
\qquad \bE = ( E_1(t, x_1), E_2(t, x_1) ),
\qquad \bB = B_3(t, x_1).
$$
Moreover, we simulate an electron distribution in a neutralizing ion background, which differs
from the multi-species Vlasov--Maxwell system in that the average current is substracted from the total one
in order for the model to be momentum preserving.
In particular, the reduced Maxwell system then reads
\begin{equation}
  \left\{
  \begin{aligned}
    \frac{\partial E_1(t,x_1)}{\partial t} &= -J_1(t,x_1) + \frac{1}{L_1} \int_0^{L_1} J_1(t,y_1) \, \text{d} y_1
    \\
    \frac{\partial E_2(t,x_1)}{\partial t} + \frac{\partial B_3(t,x_1)}{\partial x_1} &= -J_2(t,x_1) + \frac{1}{L_1} \int_0^{L_1} J_2(t,y_1) \, \text{d} y_1
    \\
    \frac{\partial B_3(t,x_1)}{\partial t} + \frac{\partial E_2(t,x_1)}{\partial x_1} &= 0.
  \end{aligned}
  \right.
\end{equation}
 In some cases this model will be further reduced to 1d1v phase space by skipping $v_2$, $E_2$ and $B_3$,
 so that the equation for $E_1$ above remains as the only field equation.

As a first test case, we consider the Weibel instability in 1d2v phase-space as studied in \cite{kraus2016gempic} with an initial value of
\begin{align*}
f(t=0,x_1,v_1,v_2) &= \frac{1}{2\pi v_{th,1} v_{th,2}} \exp \left(- \frac{1}{2}\left( \frac{v_1^2}{v_{th,1}^2} + \frac{v_2^2}{v_{th,2}^2} \right) \right), \quad x_1 \in [0,2\pi/\sfk),\\
B_3(t=0,x_1) &= \beta \cos(\sfk x_1),\\
E_2(t=0,x_1) &= 0,
\end{align*}
and $E_1(t=0,x_1)$ is computed from Poisson's equation. The parameters are set to $v_{th,1} = \frac{0.02}{\sqrt{2}}$, $v_{th,2} = \sqrt{12}v_{th,1}$, $\sfk=1.25$, $\beta = 10^{-4}$.
As a reference solution, we use a simulation with a Fourier solver with $K = 30$ modes corresponding to $M = 61$ cells (i.e., grid points),
and $N = 10^5$ particles with a piecewise affine spline shape function $S$.

As a second test case, we consider the two-stream instability in 1d1v phase-space with initial value
$$
f(t=0,x_1,v_1) = (1+ \epsilon \cos(\sfk x_1)) \frac {1}{2\sqrt{2\pi}} \left(\exp\left(-\frac{(v_1+2.4)^2}{2}\right)-\exp\left(-\frac{(v_1-2.4)^2}{2}\right)\right)
$$
with parameters $\epsilon = 0.001$ and $\sfk = 0.2$. The initial field $E_1$ is again determined from Gauss' law.
For this test case, the reference solution is also produced with a Fourier solver and a piecewise affine spline as shape function,
but the grid resolution is reduced to 31 cells (and 15 modes) while the particle number is increased to $5 \cdot 10^{6}$.
Note that this test case requires a lot more particles to produce qualitative results compared to the Weibel test case.

In Sections~\ref{sec:shape} to \ref{sec:fem} below we study the influence of different numerical parameters
using the relevant energy curves for these two test cases, namely the magnetic and electric energy,
plotted in Figures~\ref{fig:weibel_energies} and \ref{fig:tsi_energies} respectively.
In Section~\ref{sec:cons} we finally compare the long-time conservation properties of the schemes,
looking at different error curves shown in Figure~\ref{fig:conservation}.


\subsection{Influence of the shape function}
\label{sec:shape}

We first study the influence of the shape function. Here, we expect two counteracting effects: On the one hand, a higher degree of the shape function yields smoother data for the field solver which can yield better results. On the other hand, higher order smoothing kernels smear out the influence of particles which yields a damping. This latter effect is clearly seen in the simulations with $M=7$ cells (grid points) of Figures~\ref{fig:weibel_shape_n7} and \ref{fig:tsi_shape_n7}. For this coarse resolution, low order splines give rather good results whereas higher order shapes lead to a visible damping in the instability growth rate for both test cases. Increasing the number of cells to $M=15$ while keeping the number of particles constant as in Figures~\ref{fig:weibel_shape_n15} and \ref{fig:tsi_shape_n15}, we observe both effects: In this case, the degree one spline yields too noisy data for the field solver, while a degree of e.g. seven yields too high damping and an intermediate degree of four yields rather accurate results. Our results also show that when increasing also the number of particles, the choice of the shape function is of lesser importance (cf.~Figures~\ref{fig:weibel_scheme_n15} and \ref{fig:tsi_scheme_n15}).


\subsection{Influence of the space semi-discretization}
\label{sec:space_dis}

\begin{figure}[ht]
\begin{subfigure}{.5\textwidth}
  \centering
  \includegraphics[width=.8\linewidth]{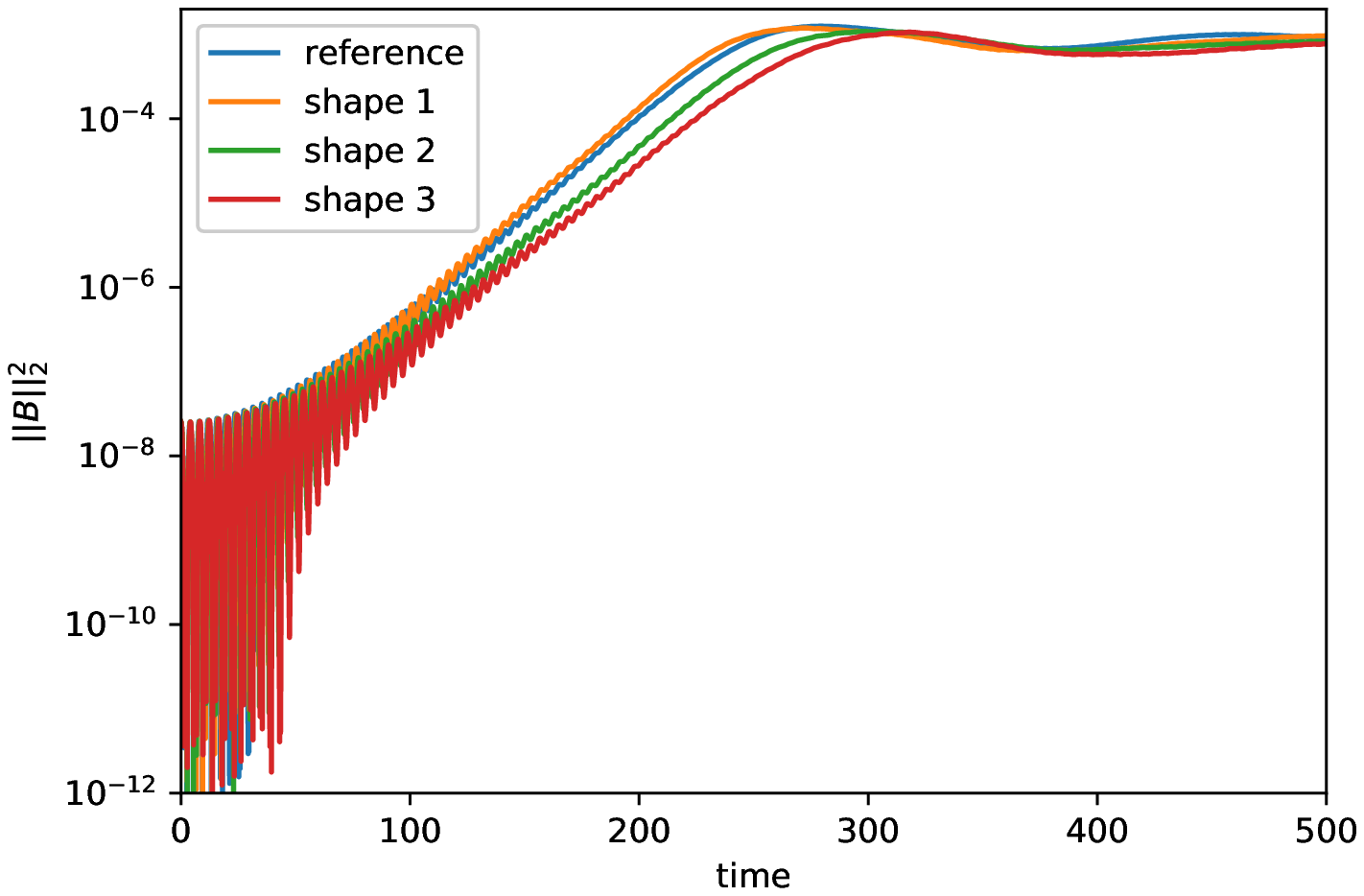}
  \caption{Varying shape function $S$ for 7 cells.}
  \label{fig:weibel_shape_n7}
\end{subfigure}
\begin{subfigure}{.5\textwidth}
  \centering
 \includegraphics[width=.8\linewidth]{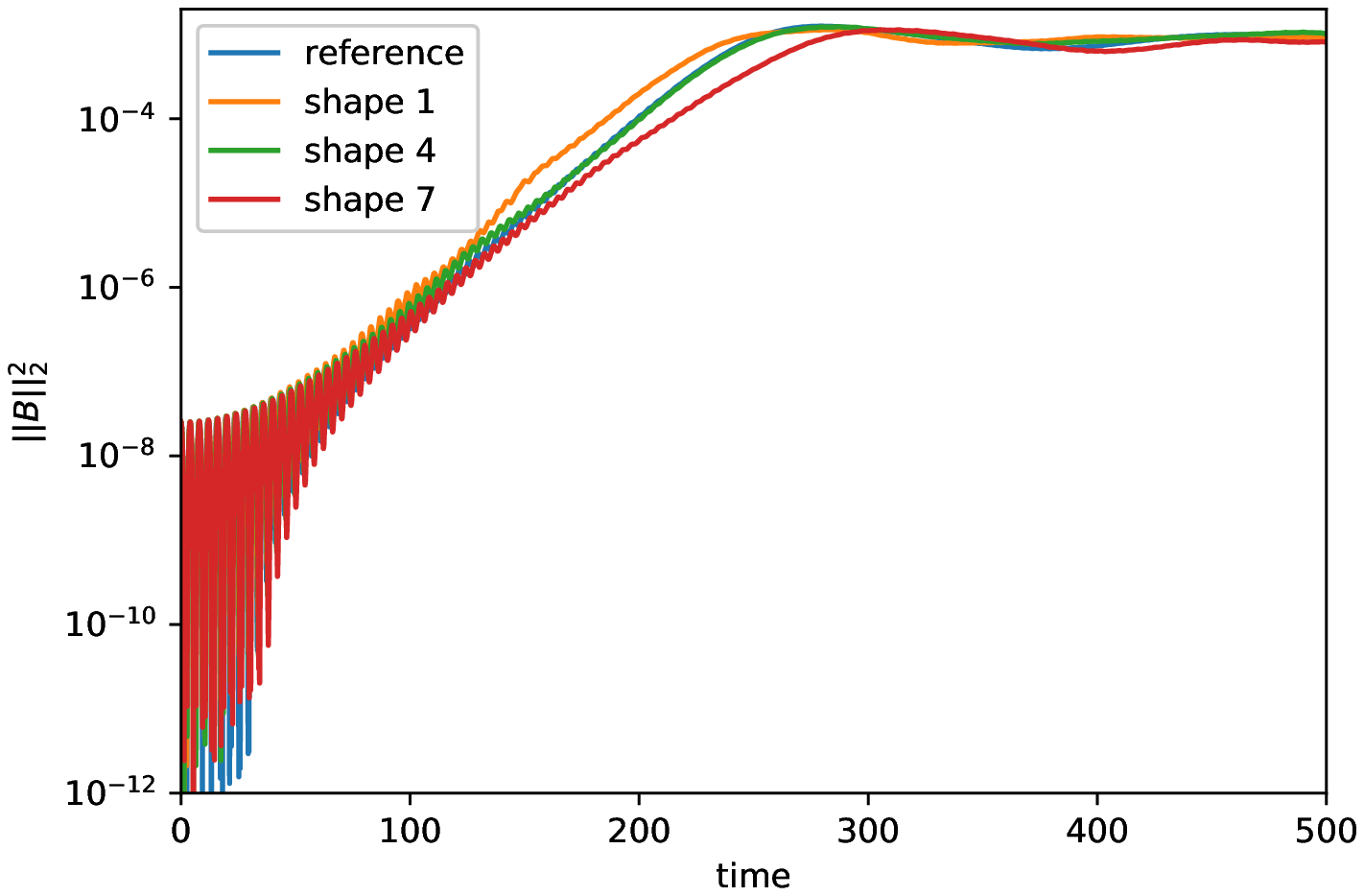}
  \caption{Varying shape function $S$ for 15 cells.}
  \label{fig:weibel_shape_n15}
\end{subfigure}
\begin{subfigure}{.5\textwidth}
  \centering
 \includegraphics[width=.8\linewidth]{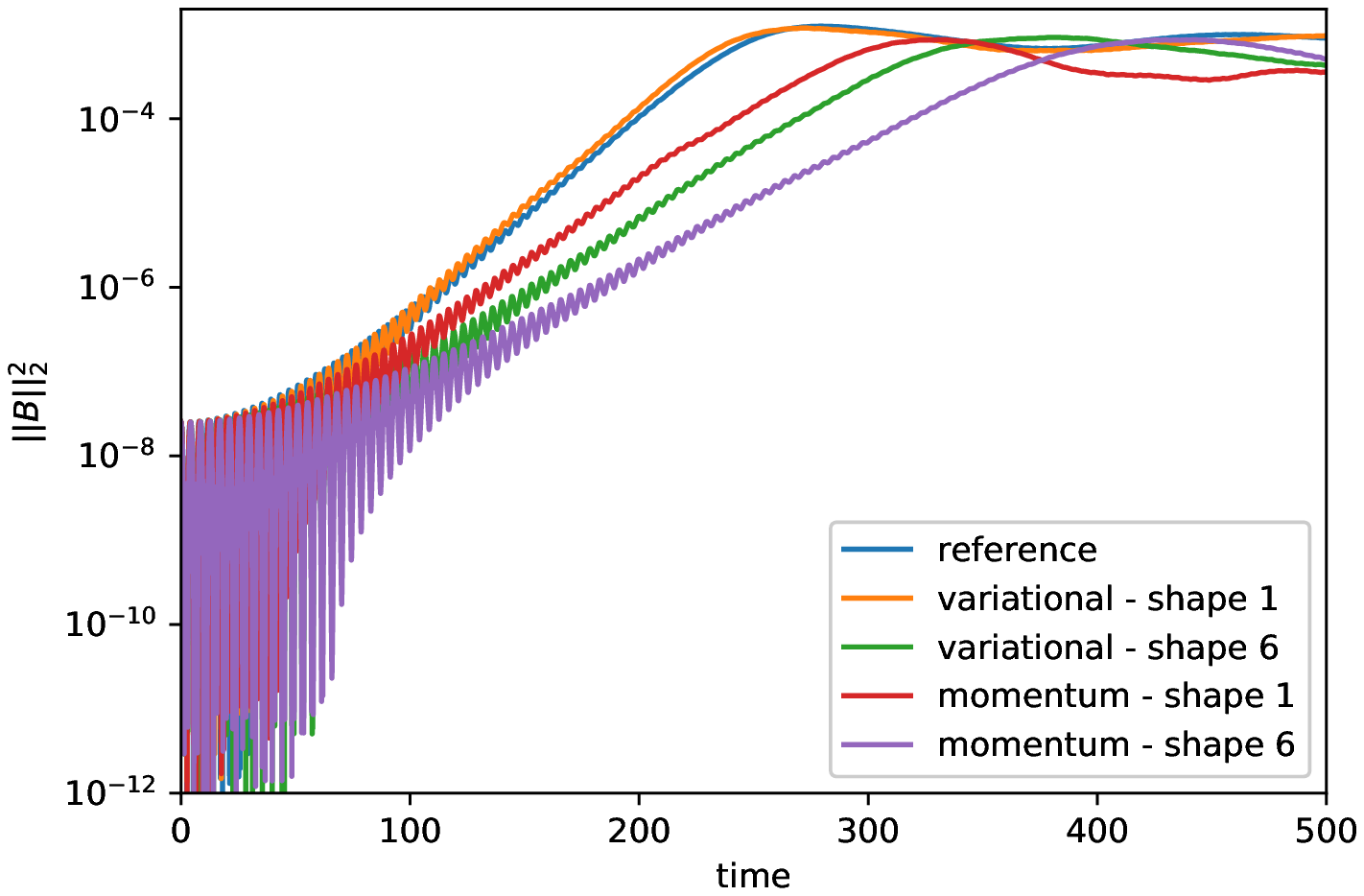}
  \caption{Varying scheme and shape $S$ for 7 cells.}
  \label{fig:weibel_scheme_n7}
\end{subfigure}
\begin{subfigure}{.5\textwidth}
  \centering
 \includegraphics[width=.8\linewidth]{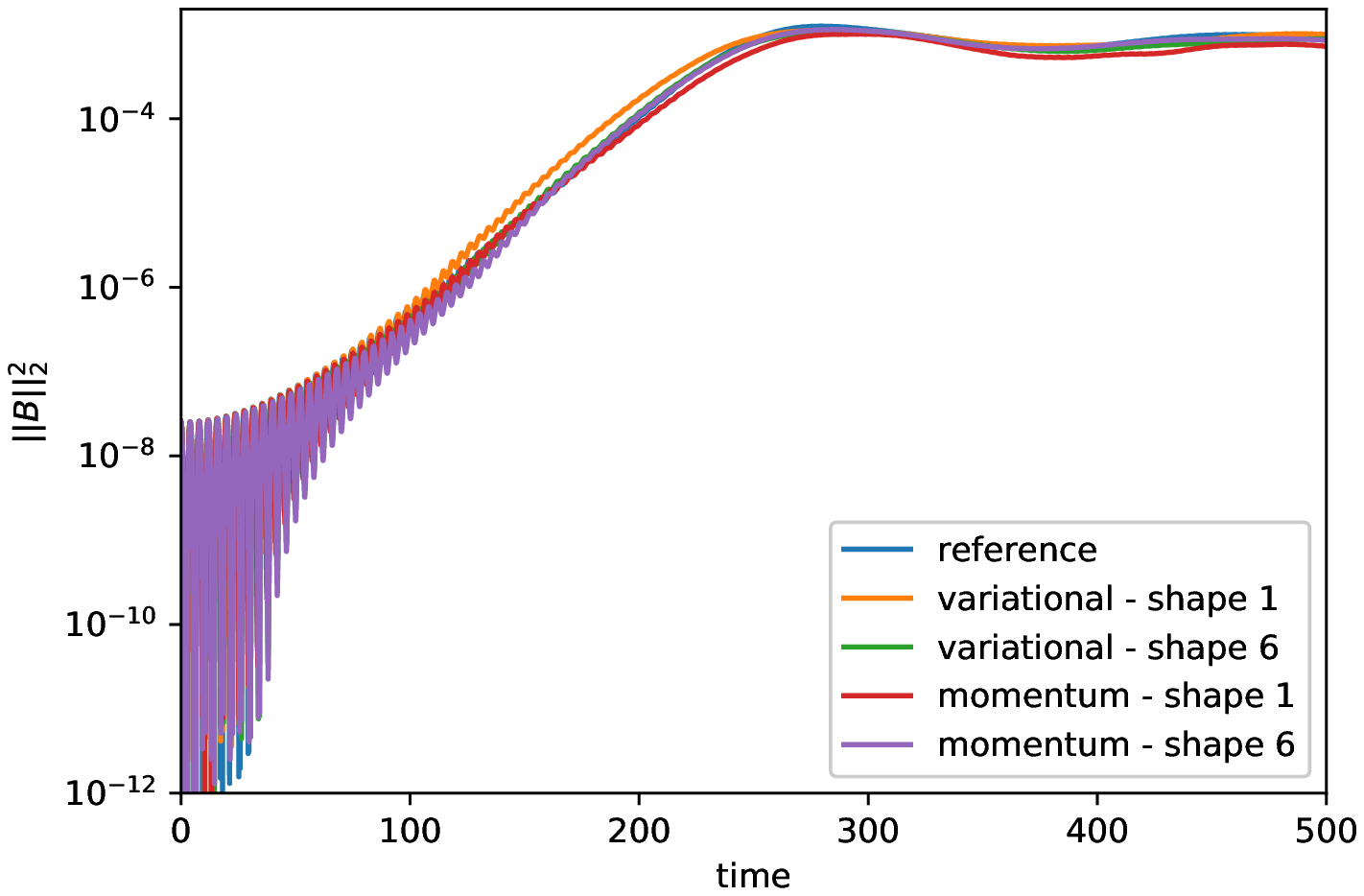}
  \caption{Varying scheme and shape $S$ for 15 cells, using more particles.}
  \label{fig:weibel_scheme_n15}
\end{subfigure}
\begin{subfigure}{.5\textwidth}
  \centering
 \includegraphics[width=.8\linewidth]{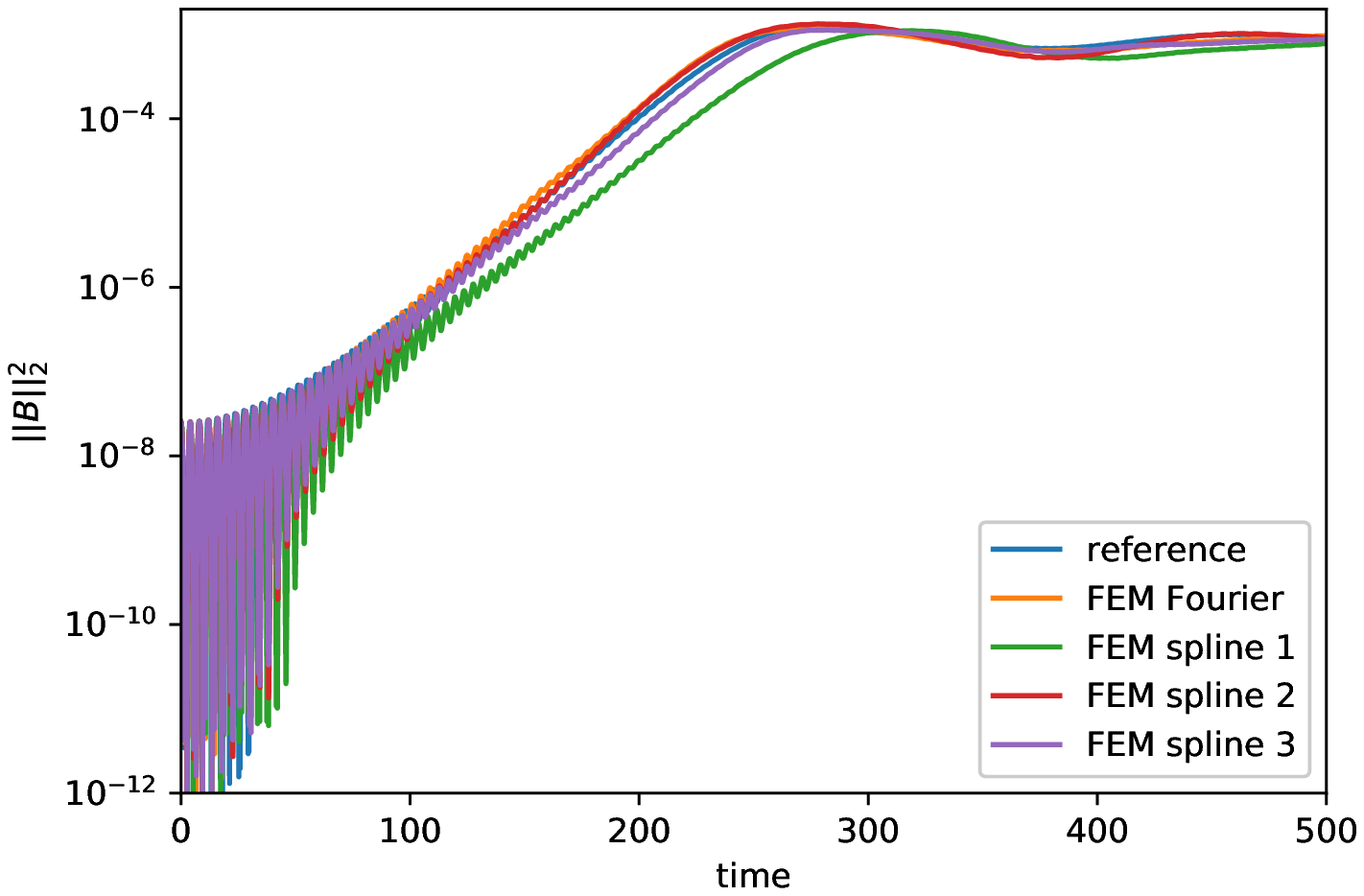}
  \caption{Varying finite element solver for 7 cells and first degree spline as shape function.}
  \label{fig:weibel_fem}
\end{subfigure}
\begin{subfigure}{.5\textwidth}
  \centering
 \includegraphics[width=.8\linewidth]{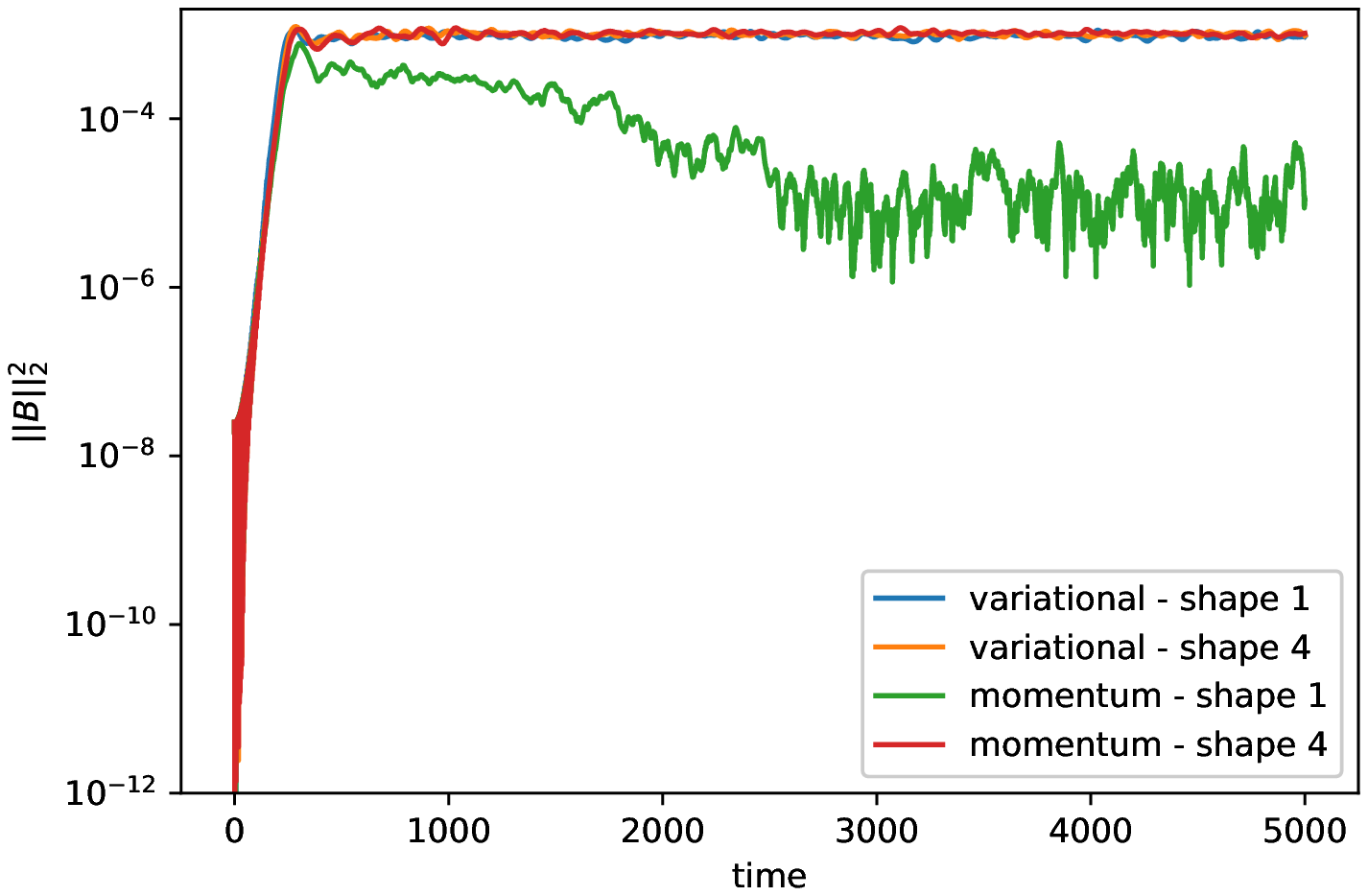}
  \caption{Varying scheme and shape $S$ for 15 cells.}
  \label{fig:weibel_long}
  \end{subfigure}
\caption{Weibel instability: Time evolution of the magnetic energy for various configurations. In all figures except (e), a spectral finite element solver is used and the degree of the shape function is given in the legend. In (e), the shape function is a spline of degree 1 and the legend indicates the degree of the finite element solver. The number of particles is 1000 in all figures except (d) where it is 4000. Figures (a), (b), and (e) show results with the variational scheme and figures (c),(d), and (f) compare the variational and the momentum-preserving schemes (see legend). All simulations use the Hamiltonian splitting time propagator. } 
\label{fig:weibel_energies}
\end{figure}

\begin{figure}[ht]
\begin{subfigure}{.5\textwidth}
  \centering
  \includegraphics[width=.8\linewidth]{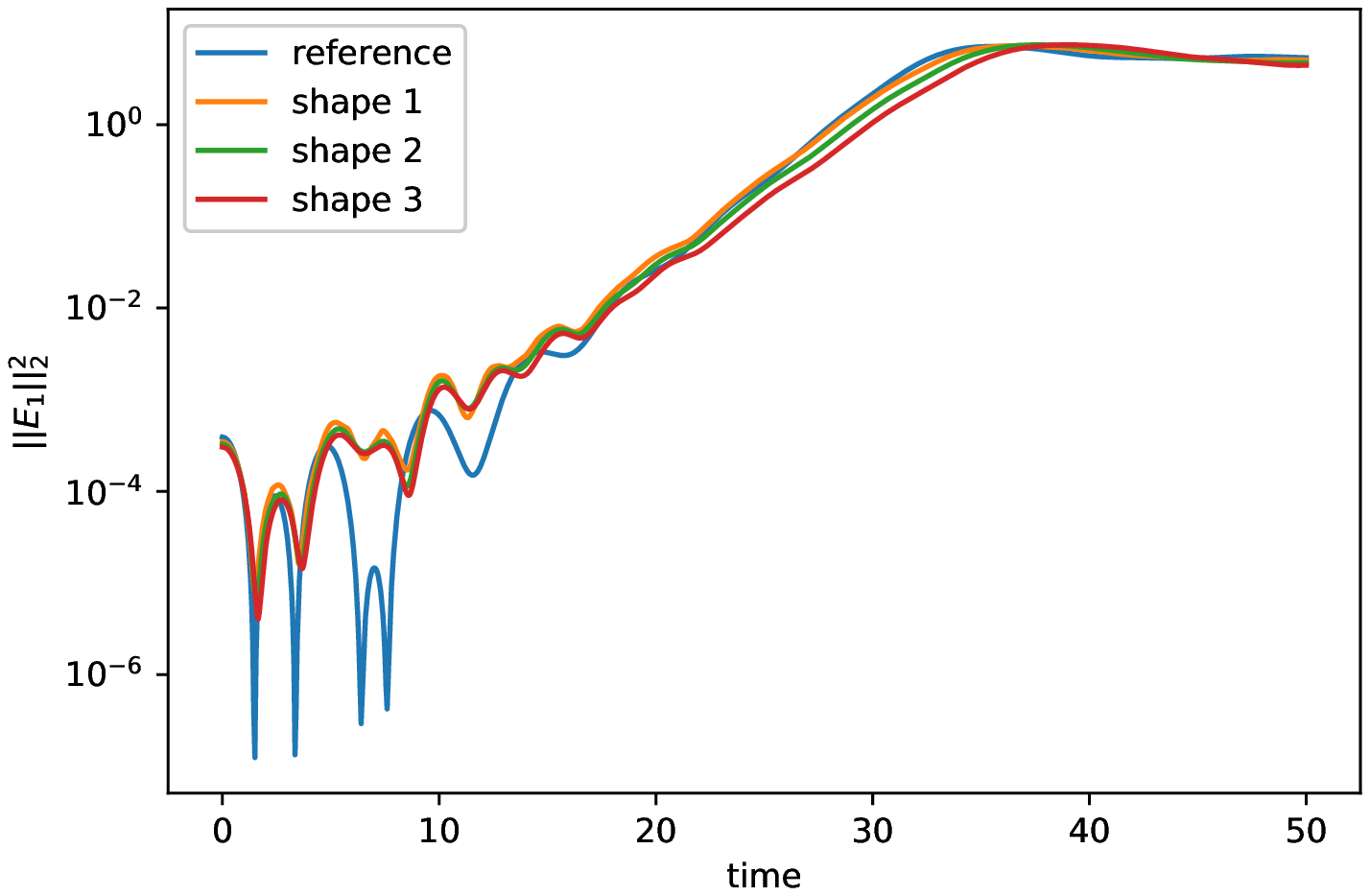}
  \caption{Varying shape function $S$ for 7 cells.}
  \label{fig:tsi_shape_n7}
\end{subfigure}
\begin{subfigure}{.5\textwidth}
  \centering
 \includegraphics[width=.8\linewidth]{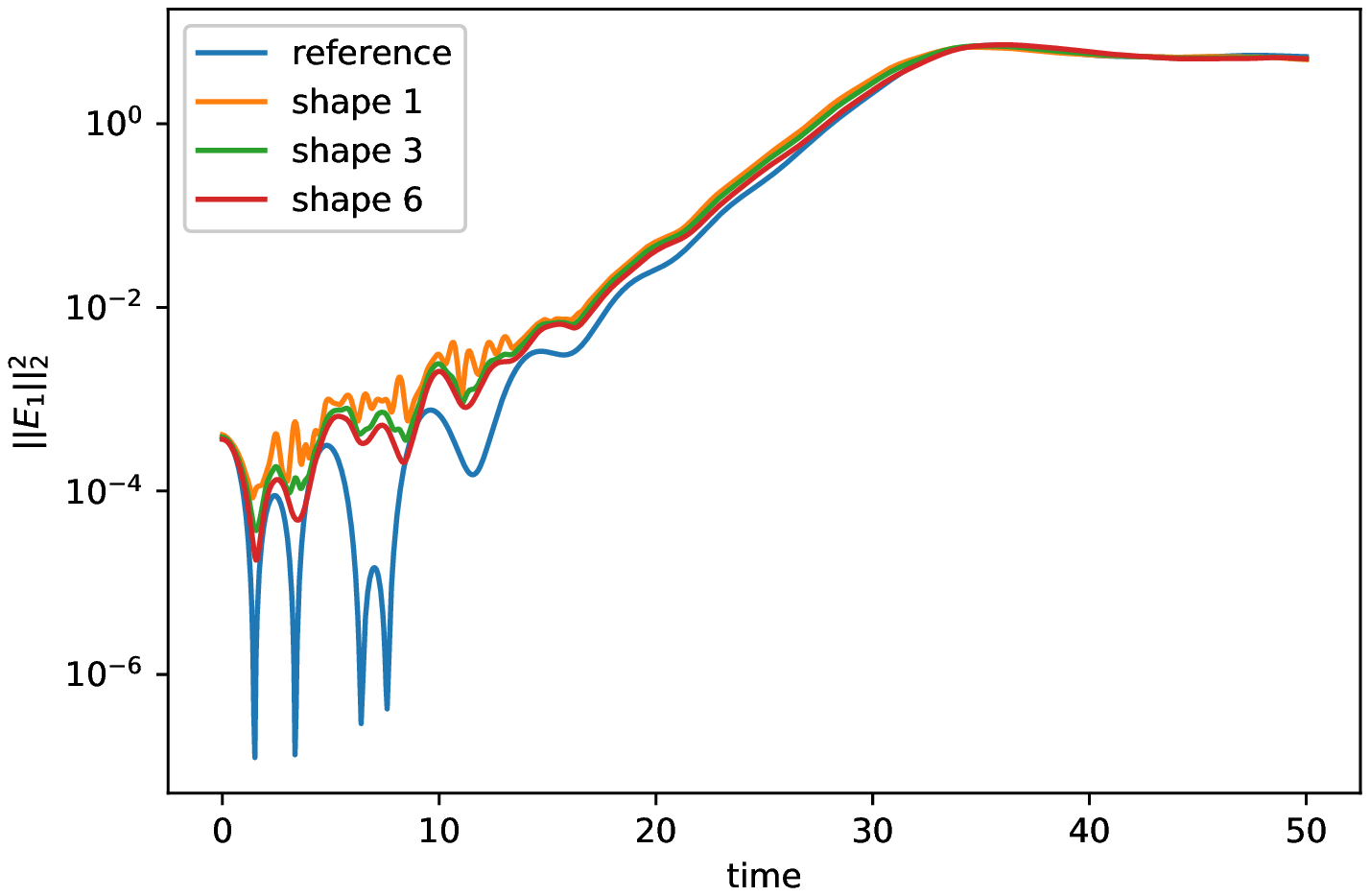}
  \caption{Varying shape function $S$ for 15 cells.}
  \label{fig:tsi_shape_n15}
\end{subfigure}
\begin{subfigure}{.5\textwidth}
  \centering
 \includegraphics[width=.8\linewidth]{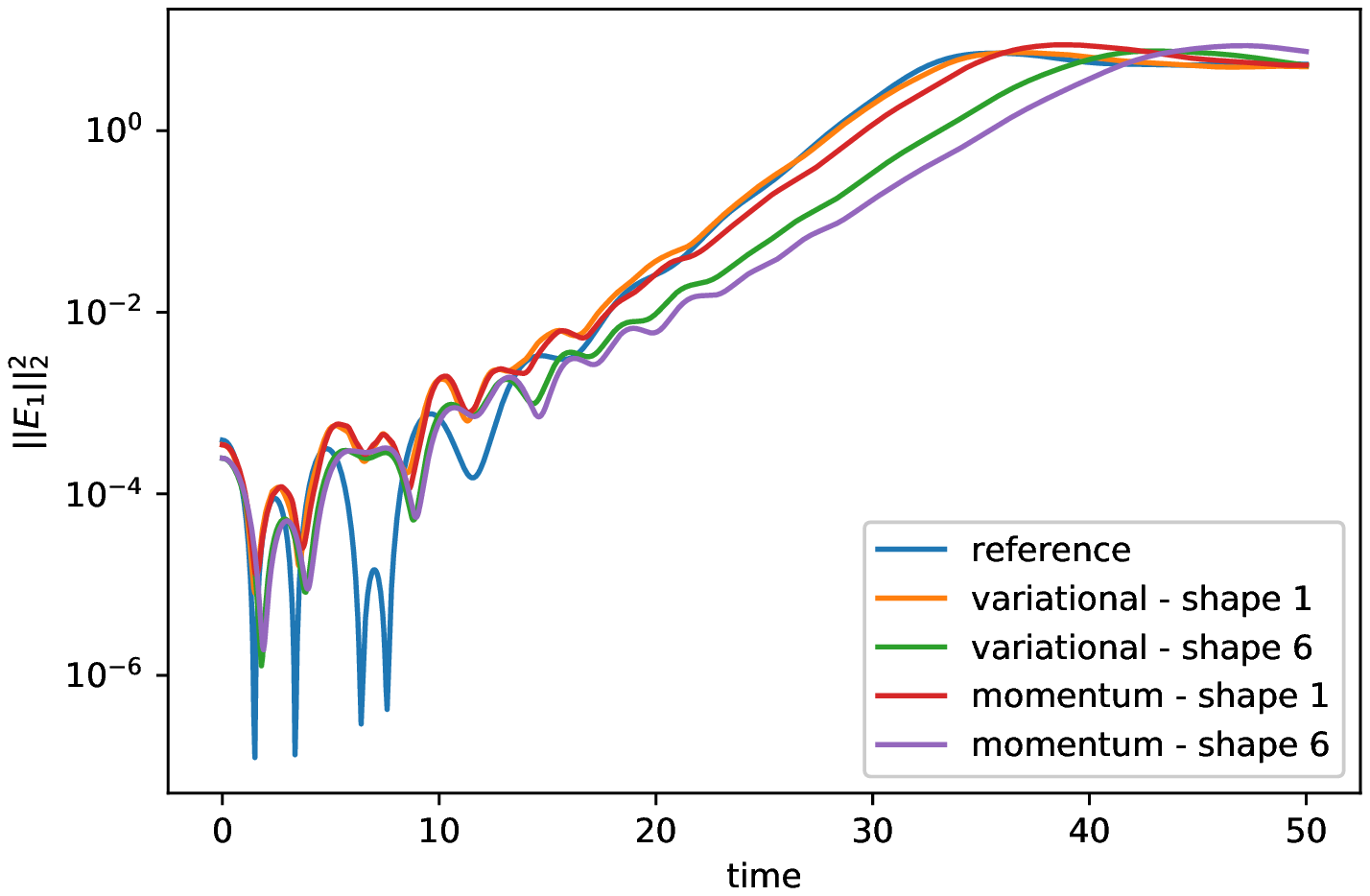}
  \caption{Varying scheme and shape  $S$ for 7 cells.}
  \label{fig:tsi_scheme_n7}
\end{subfigure}
\begin{subfigure}{.5\textwidth}
  \centering
 \includegraphics[width=.8\linewidth]{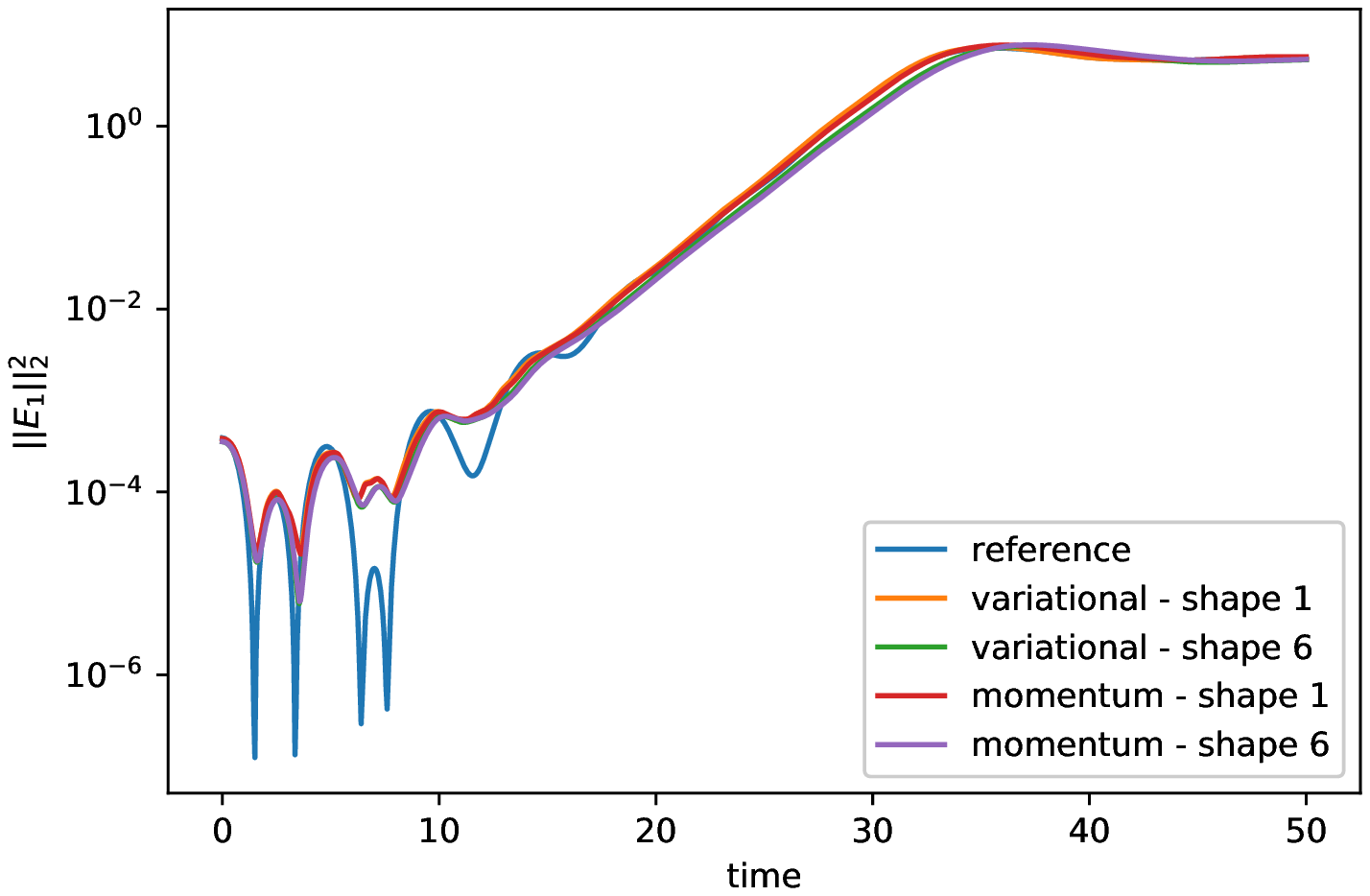}
  \caption{Varying scheme and shape $S$ for 15 cells, using more particles.}
  \label{fig:tsi_scheme_n15}
\end{subfigure}
\begin{subfigure}{.5\textwidth}
  \centering
 \includegraphics[width=.8\linewidth]{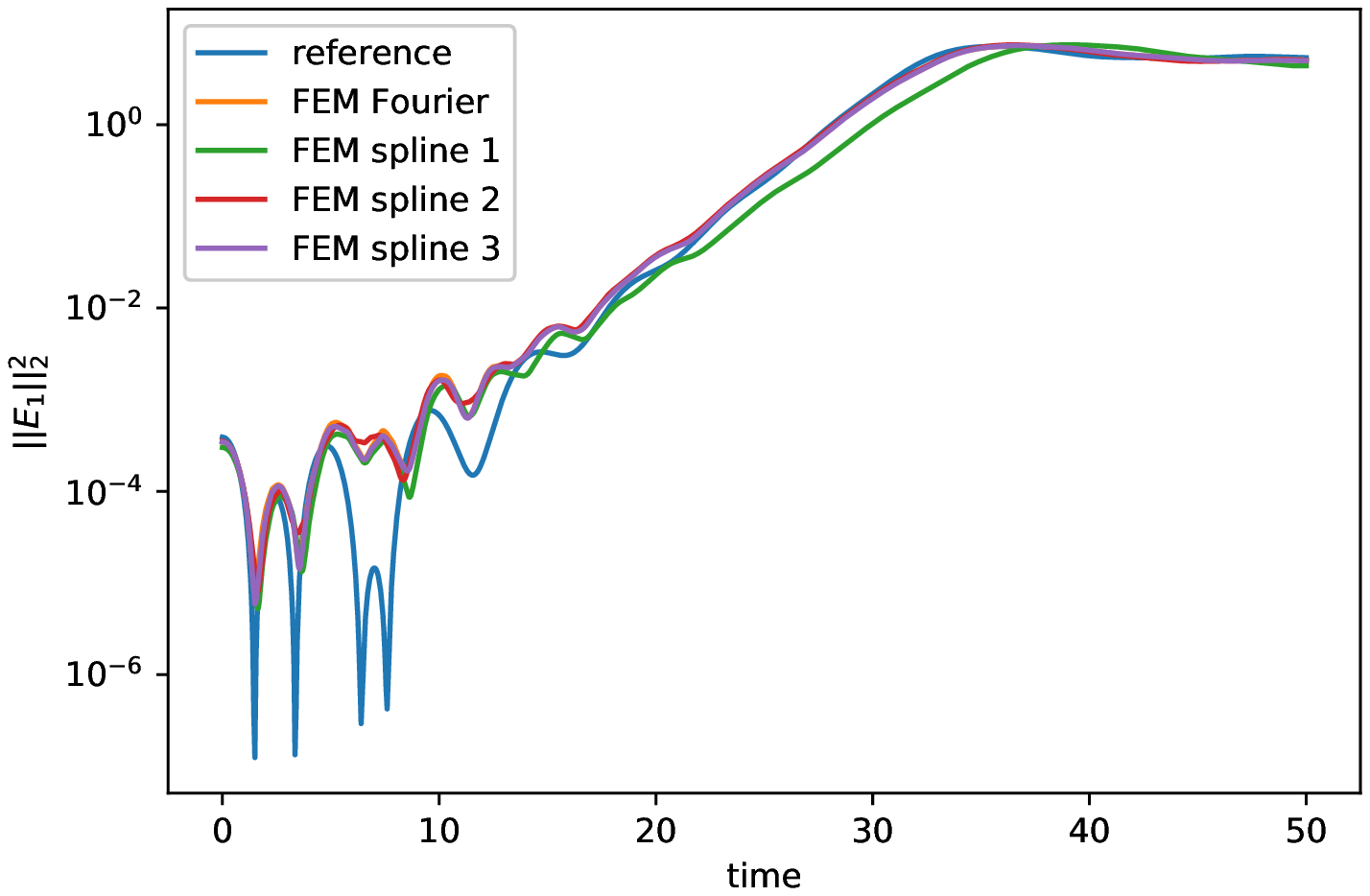}
  \caption{Varying finite element solver for 7 cells and first order spline as shape function.}
  \label{fig:tsi_fem}
\end{subfigure}
\begin{subfigure}{.5\textwidth}
  \centering
 \includegraphics[width=.8\linewidth]{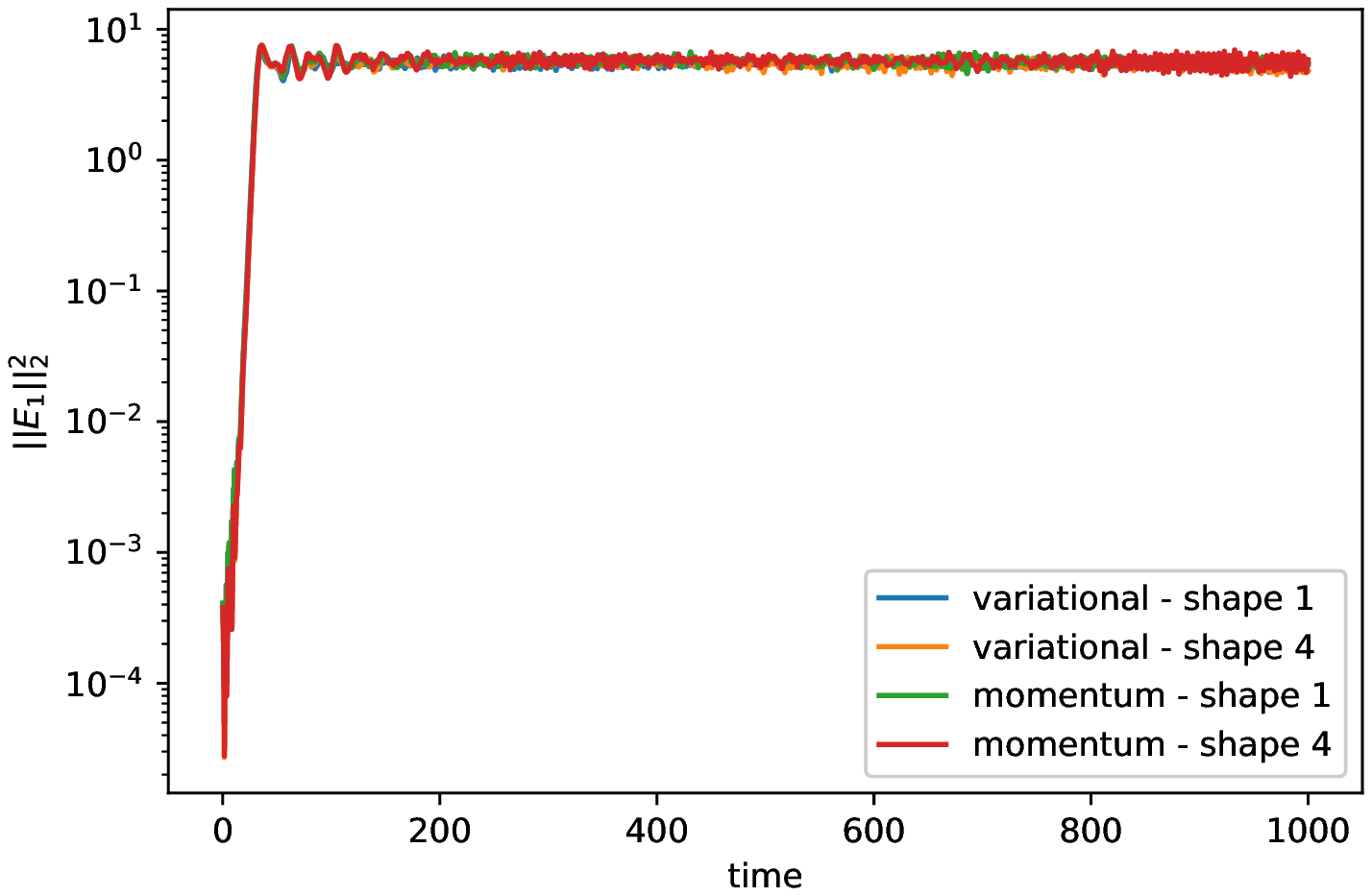}
  \caption{Varying scheme and shape $S$ for 7 cells.}
  \label{fig:tsi_long}
\end{subfigure}

\caption{Two-stream instability: Time evolution of the first component of the electric energy for various configurations. In all figures except (e), a spectral finite element solver is used and the degree of the spline shape is given in the legend. In (e), the shape function is a spline of degree 1 and the legend indicates the degree of the finite element solver. The number of particles is 48000 in all figures except (d) where it is 192000. Figures (a), (b), and (e) show results with the variational scheme and figures (c),  (d), and (f) compare the variational and the momentum-preserving schemes (see legend). All simulations use the Hamiltonian splitting time propagator. } 
\label{fig:tsi_energies}
\end{figure}

In Figures \ref{fig:weibel_scheme_n7}--\ref{fig:weibel_scheme_n15} and \ref{fig:tsi_scheme_n7}--\ref{fig:tsi_scheme_n15}
we next compare the variational scheme presented in Section~\ref{sec:variational} with the momentum-preserving variant from Section~\ref{sec:momentum}.
Here we use the spectral finite element solver and a Hamiltonian splitting time discretization. With this configuration, the momentum-preserving scheme yields clearly worse results for the coarse resolution runs (in Figures~\ref{fig:weibel_scheme_n7} and \ref{fig:tsi_scheme_n7}), as the instability growth rate is damped similarly as with higher order shape functions. With increased resolution (i.e., using twice as many cells and four times as many particles for
both test cases), we find that both schemes yield rather good results for various orders of the shape function (in Figures~\ref{fig:weibel_scheme_n15} and \ref{fig:tsi_scheme_n15}). Finally, we see in Figure~\ref{fig:weibel_long} that the long-time accuracy of the variational semi-discretization can be significantly better than that of the momentum-preserving one: here the Weibel instability is run with a small number of particles and we find a qualitatively wrong behavior for the momentum-preserving scheme using a piecewise affine shape function, where other schemes perform correctly. In Figure~\ref{fig:tsi_long} a similar comparison is done with the two-stream instability, using a higher particle resolution as required for this test case to produce qualitatively correct results. The long-time behavior is then found to be qualitatively good for the different schemes and shapes.

\subsection{Influence of the finite element solver}
\label{sec:fem}

In Figures \ref{fig:weibel_fem} and \ref{fig:tsi_fem} we then compare the different field solvers, namely the spectral solver and finite element solvers based on splines of degree one to three. Using a piecewise affine spline for the shape function and low resolution runs
 we find that the accuracy of the low order fem solver is of bad quality and it improves for higher orders and for the spectral solver.
 This observation holds for the two test cases.

\subsection{Conservation properties}
\label{sec:cons}

\begin{figure}
  \begin{subfigure}{.5\textwidth}
  \centering
 \includegraphics[width=.8\linewidth]{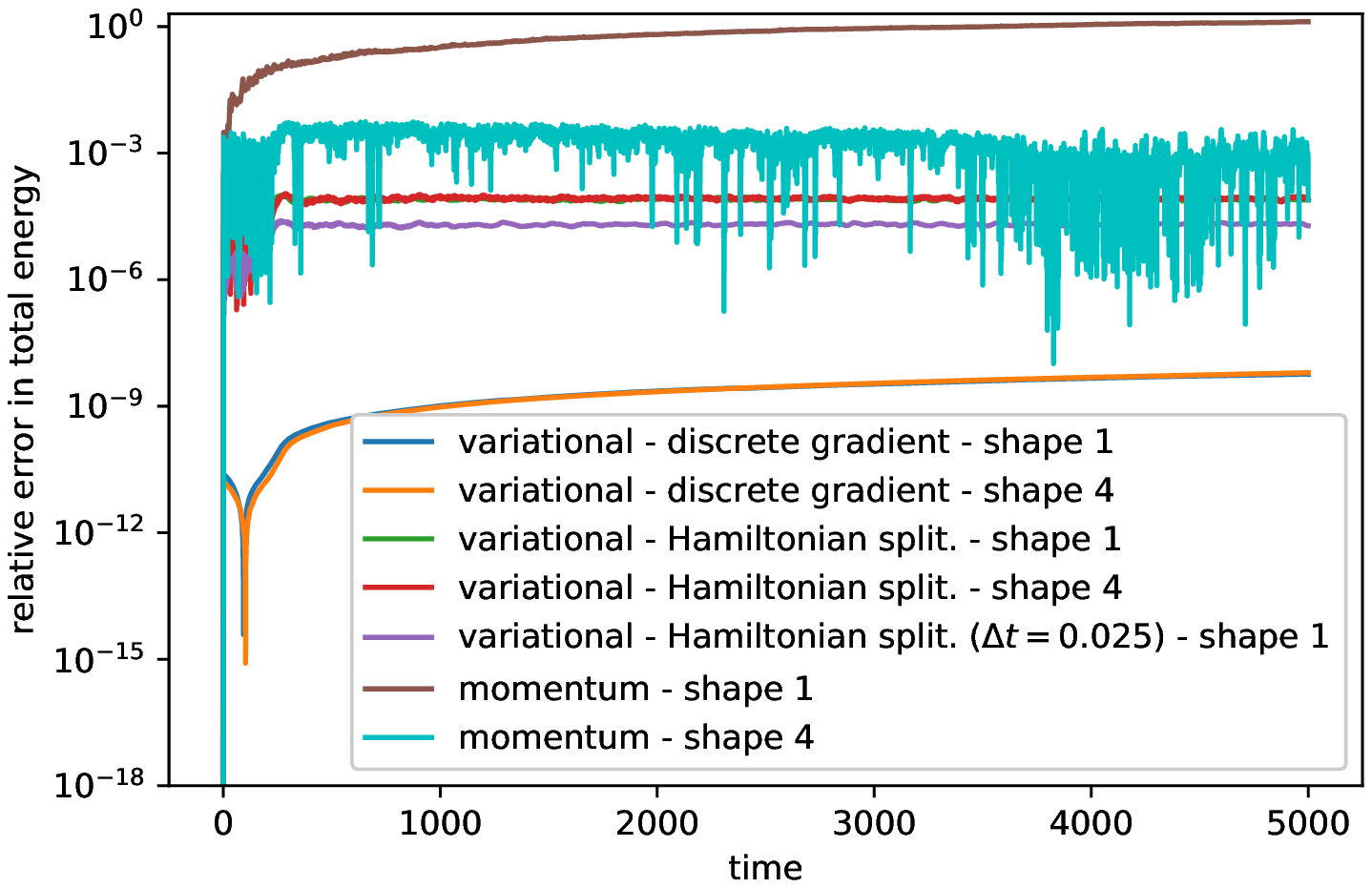}
  \caption{Weibel instability: energy errors.}
  \label{fig:weibel_toten}
  \end{subfigure}
\begin{subfigure}{.5\textwidth}
  \centering
 \includegraphics[width=.8\linewidth]{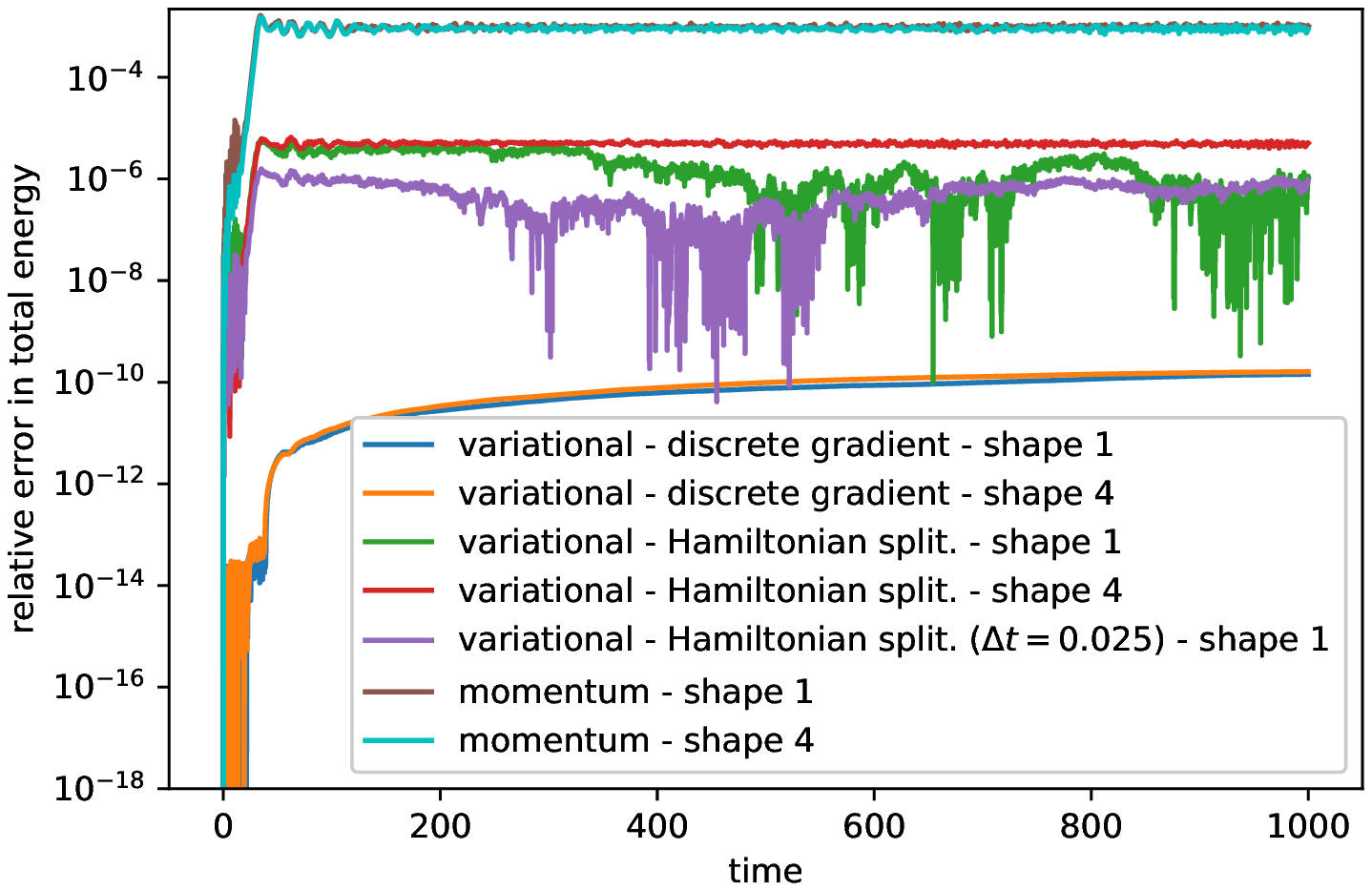}
  \caption{Two-stream instability: energy errors.}
  \label{fig:tsi_toten}
\end{subfigure}
  \begin{subfigure}{.5\textwidth}
  \centering
 \includegraphics[width=.8\linewidth]{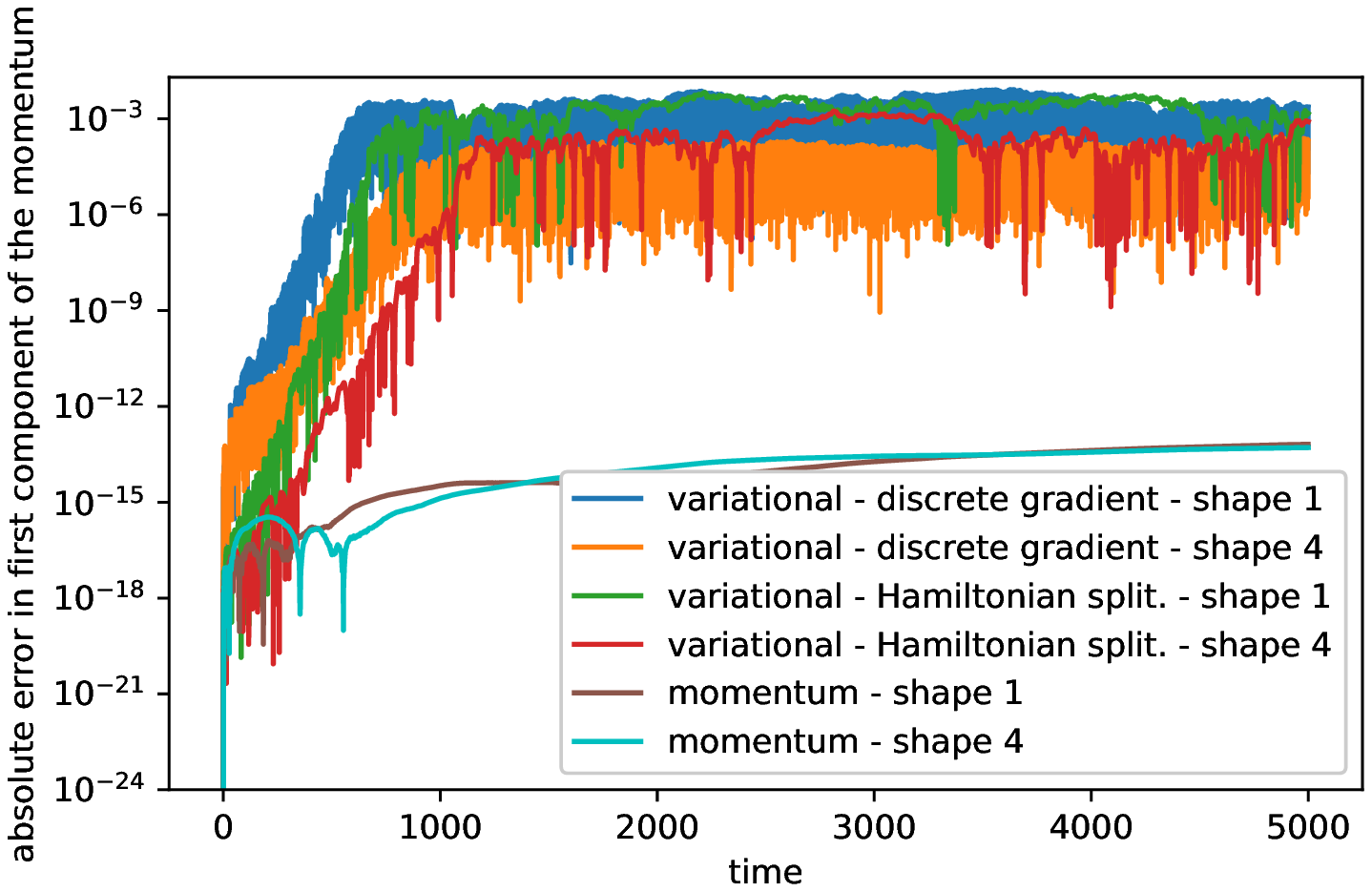}
  \caption{Weibel instability: momentum errors.}
  \label{fig:weibel_mom}
  \end{subfigure}
\begin{subfigure}{.5\textwidth}
  \centering
 \includegraphics[width=.8\linewidth]{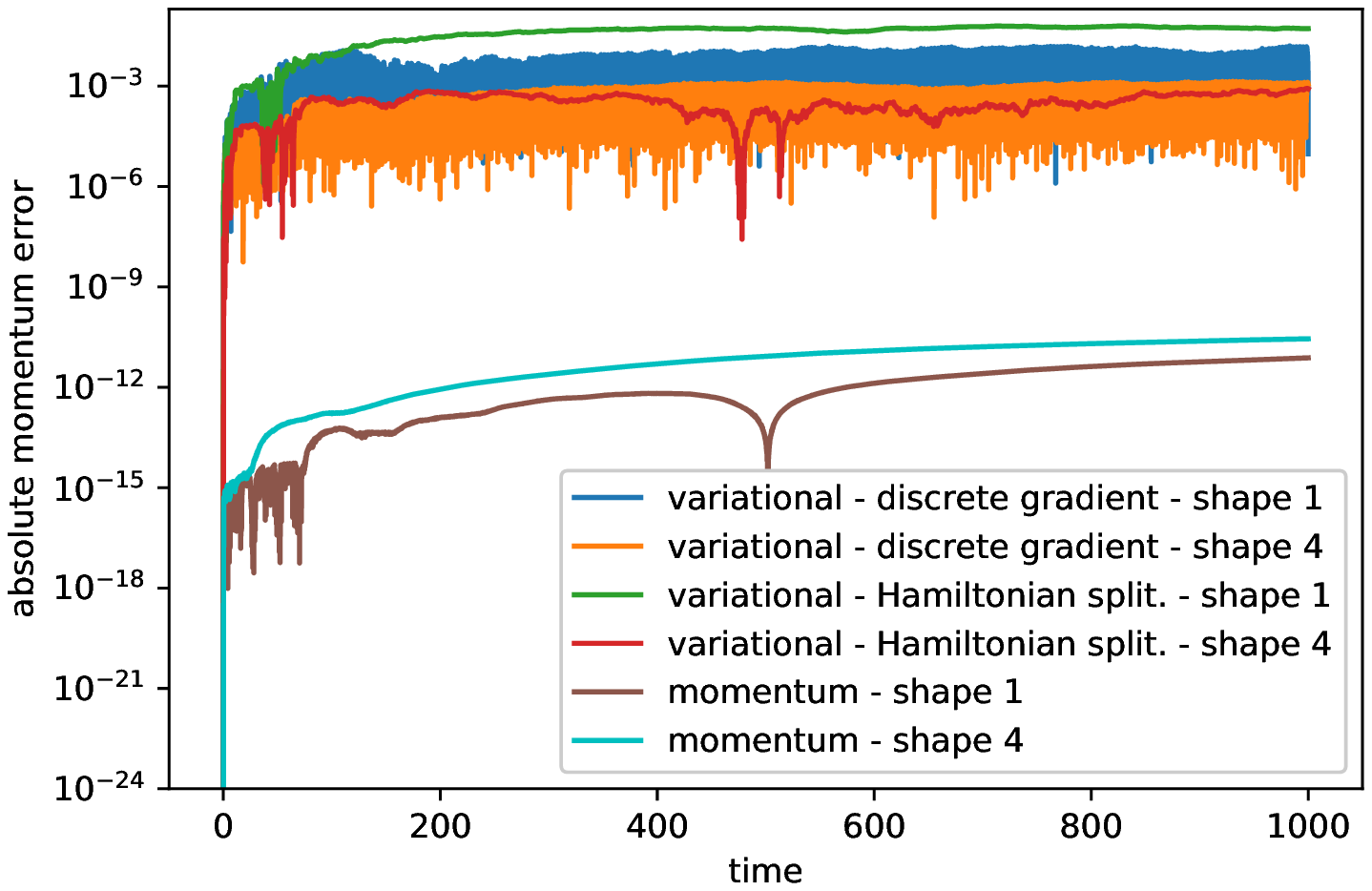}
  \caption{Two-stream instability: momentum errors.}
  \label{fig:tsi_mom}
\end{subfigure}
  \begin{subfigure}{.5\textwidth}
  \centering
 \includegraphics[width=.8\linewidth]{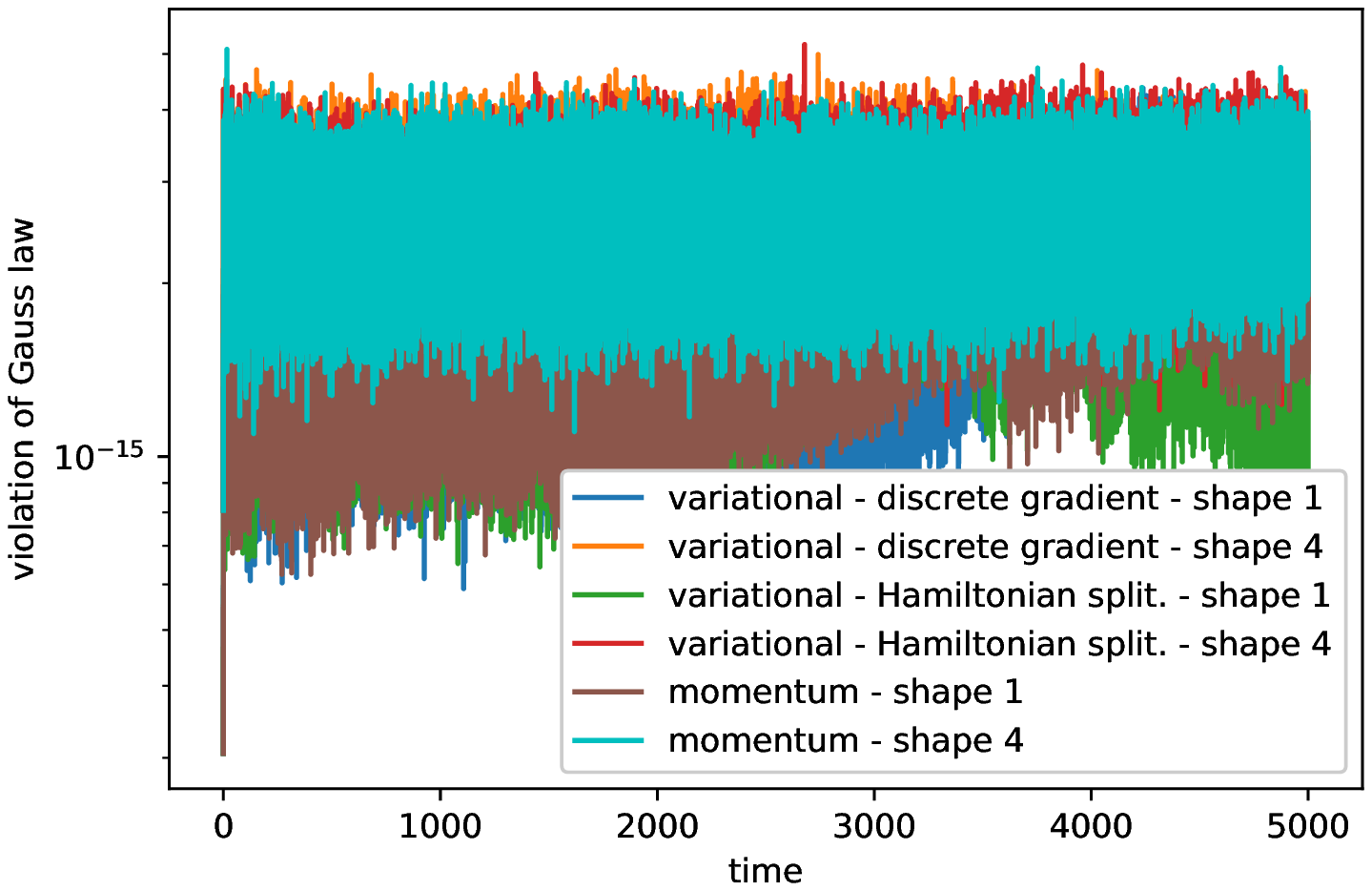}
  \caption{Weibel instability: Gauss' law errors.}
  \label{fig:weibel_gauss}
\end{subfigure}
\begin{subfigure}{.5\textwidth}
  \centering
 \includegraphics[width=.8\linewidth]{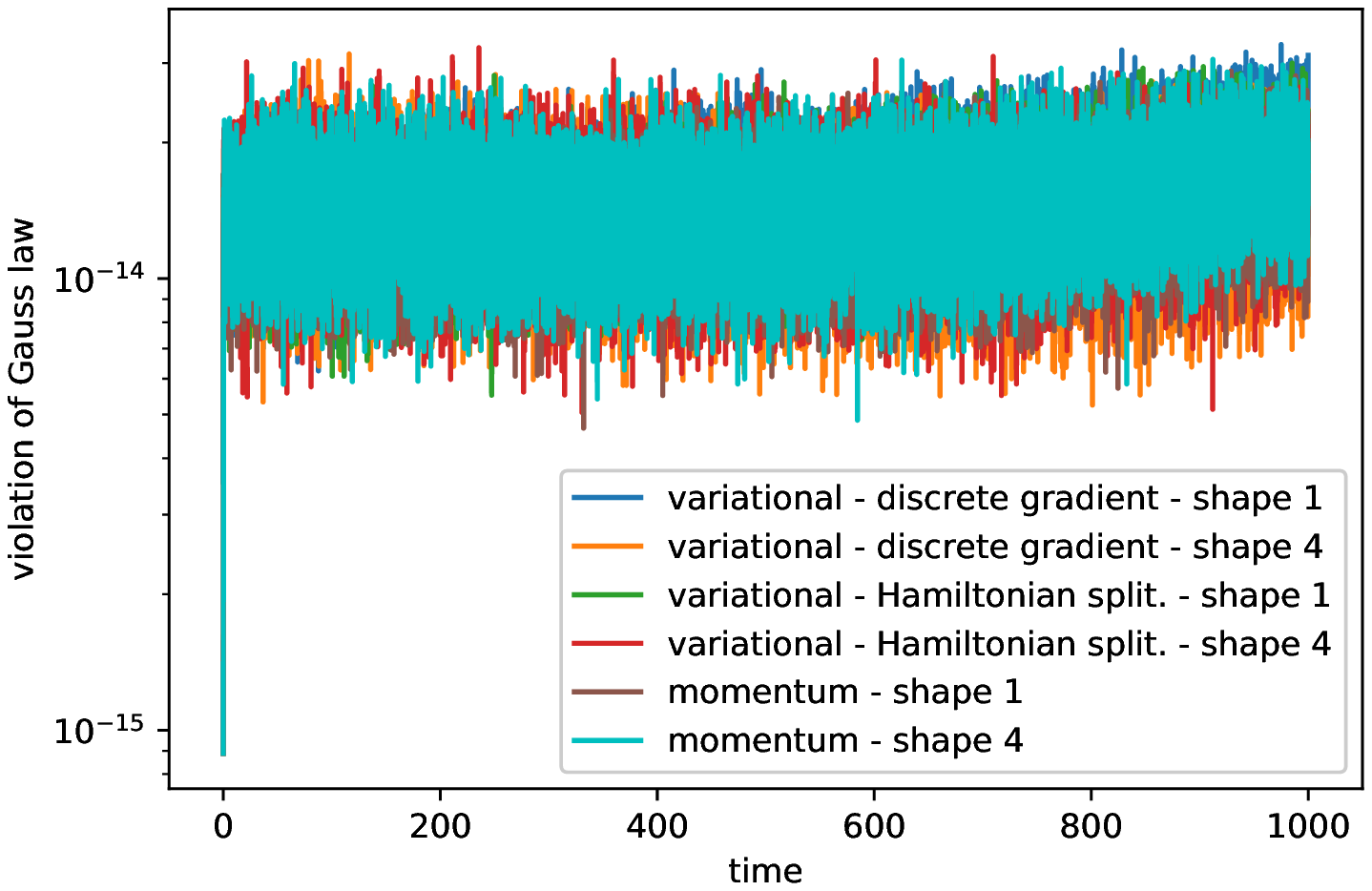}
  \caption{Two-stream instability: Gauss' law errors.}
  \label{fig:tsi_gauss}
\end{subfigure}
\caption{Conservation properties: energy, momentum and Gauss' law errors are shown for long-time simulations of the Weibel instability with $N=1000$ particles (left column) and of the two-stream instability with
$N = 48000$ particles (right column). In all the runs, a spectral solver with $K=7$ modes is used for the field.}
\label{fig:conservation}
\end{figure}

We now compare the conservation properties of the various methods. For this we consider long times simulations with both the variational and the momentum-preserving discretizations. For the time stepping, we consider in both cases a Hamiltonian splitting as before but we also provide the solution with an energy-conserving discrete gradient propagator for the variational scheme to show that the semi-discretization is indeed energy-conserving. Figures \ref{fig:weibel_toten} and \ref{fig:tsi_toten} show the relative error in energy conservation for the various runs. We can see that the energy is conserved up to the tolerance of the linear solvers for the variational scheme with an energy-conserving discrete gradient time propagator. If we use the Hamiltonian splitting instead, there is an energy error but its behavior is oscillatory and decreases with decreased time step. This is the typical behavior for such Poisson integrators. Finally, we see that the energy error is larger for the momentum-preserving scheme, in particular for the low order shape function with a low particle resolution. For the variational scheme, on the other hand, the energy error does not depend on the shape function.

Figure \ref{fig:weibel_mom} and \ref{fig:tsi_mom} show the error in momentum for the various methods. We can see that the momentum-preserving scheme indeed preserves momentum up to machine precision. On the other hand, for the variational scheme the error in momentum increases as soon as the nonlinear phase of the simulations starts and later flattens out at a certain level. As this error level seems to be rather independent of the propagator, and is smaller for higher order shape functions, we conjecture that it is dominated by the error in the spatial semi-discretization.

Finally the error in Gauss' law as a function of time is shown in Figures \ref{fig:weibel_gauss} and \ref{fig:tsi_gauss} for the two test cases, respectively. We can see that all scheme preserve Gauss' law to machine precision. 

\section{Acknowledgements}
This work has been carried out within the framework of the EUROfusion Consortium and has received funding from the Euratom research and training programme 2014-2018 and 2019-2020 under grant agreement No 633053. The views and opinions expressed herein do not necessarily reflect those of the European Commission.

\bibliographystyle{plainnat}
\bibliography{gempic_strong_ampere}

\end{document}